\documentclass[a4paper,10pt]{article}
\usepackage{amsfonts}
\usepackage{amssymb}
\usepackage{graphicx}
\usepackage{mathtools}
\usepackage{skak}
\usepackage{MnSymbol}

\usepackage{amscd}
\usepackage{MnSymbol}

\usepackage{wasysym}

\usepackage{ifsym}

\makeatletter

\newcommand{\rmnum}[1]{\romannumeral #1}
\newcommand{\Rmnum}[1]{\expandafter\@slowromancap\romannumeral #1@}
\makeatother

\title{ Dirac and Plateau Billiards in  Domains with Corners.}

\author{Misha Gromov}

\begin{document}

\maketitle \tableofcontents

\begin{abstract} Groping our way toward a theory of singular spaces with 
positive  scalar curvatures we look at  the Dirac operator and  a generalized Plateau problem 
in Riemannian manifolds with corners. 
  Using these, we prove that
 
  {\it the set of  $C^2$-smooth Riemannian  metrics $g$ on a smooth  manifold $X$, such that $scal_g(x)\geq \kappa(x)$, is closed under $C^0$-limits of Riemannian metrics  for all continuous functions $\kappa$ on $X$.}

Apart from that our progress is limited  but 
we formulate many conjectures. 
All along, we    emphasize   geometry, rather than topology of manifolds with their scalar curvatures bounded from below.

\end{abstract}

\section {Setting the Stage.}

A closed subset $P=P^n$ in a smooth $n$-manifold $X$
is called a {\it cornered} or  {\it curve-faced  polyhedral $n$-domain}  of depth $d=0, 1,...,n$ if   the 
 boundary of  $\partial  P$ of $P$ is decomposed
 into the union of  a countable  locally finite 
  (e.g. finite) family of (possibly disconnected) {\it $(n-1)$-faces} $Q_i=Q_i^{n-1}$  with a
   distinguished set of  {\it  adjacent pairs of faces}   $(Q_i, Q_j)$, such that 
   
 $\bullet $
 every face $Q_i$ is contained  in a smooth  hypersurface 
 $Y_i =Y_i^{n-1} \subset X$, where $Y_i $ is transversal to $Y_i$ for all adjacent pairs of
 $(n-1)$-faces;
 
 $\bullet \bullet$ the boundary $\partial Q_i$ of  each $Q_i \subset Y_i$ equals the union of the intersections 
 $Q_i\cap Q_j$ for all faces $Q_j$ that are adjacent to a given $Q_i$, where the corresponding decompositions 
 $\partial Q_i=\bigcup_j Q_i\cap Q_j$ give  polyhedral $(n-1)$-domain structures of depth $d-1$ to all $Q_i$.
 
 This defines the notion of a   polyhedral  domain structure by induction on $d$, where
  polyhedral domains  of depth zero are  non-empty closed subsets $P\subset X$ with empty boundaries, i.e. just smooth manifolds with no extra structures  and   domains of depth one are those bounded by smooth hypersurfaces.

 A polyhedral domain $P$ is called {\it cosimplicial } if the intersections of  all $k$-tuples of mutually non-equal  $(n-1)$-faces  satisfy
    $$dim (Q_{i_1}\cap Q_{i_2}\cap ... \cap Q_{i_k})\leq n-k+1.$$

 \vspace {1mm}
  
 {\it Polyhedra, Edges, Corners.} Our attention will be focused  on   the {\it  boundaries} of cornered domains   that are   unions of faces,
$\cup_iQ_i$.  We call these boundaries  {\it curve-faced polyhedra} or {\it polyhedral hypersurfaces}
and  may  occasionally denote them   $P$ rather than  $\partial P$.

The codimension two faces of $P$ that are (n-2)-dimensional intersections
 of faces $Q_i$ are called {\it edges of $P$} and the higher codimension faces are called {\it corners.} Thus, corners appear starting from $d=3$, such as the ordinary corners of a cube in the $3$-space.

 \vspace {2mm}

 A {\it cornered $n$-manifold} structure in a $P$ is defined by an  extension of $P$ to a smooth 
  manifold $X\supset P$ where $P$ makes a polyhedral domain. The tangent bundle $T(P)$
  equals, by definition, the tangent bundle of $X$ restricted to $P$.

  \vspace {1mm}
 
{\it A     Riemannian manifold with corners} is a  cornered $n$-manifold $P$ with a Riemannian metric $g=g_P$ on $T(P)$.
 
   \vspace {1mm}
  
  {\it Cartesian products} of cornered manifolds, 
  $P=P_1\times P_2\times...\times P_m,$  
 come with  natural (product) corner structures 
 where
 $ dim(P)= dim(P_1)+dim(P_2)+...+dim(P_m)$  and  
 $depth(P)=  depth(P_1)+depth(P_2)+...+depth(P_m).$ 
If $P_i$ are Riemannian, then the product carries the natural 
Riemannian Cartesian product structure as well.
   
 The simplest instance of this is the unit {\it Euclidean  $n$-cube}, 
 sometimes denoted   $\square$ or  $\square^{n}$ (instead of  more logical  $\square^{n/2}$) 
 that is the Cartesian $n$-th power $[0,1]^n$ of
 the segment $[0,1]$. 
 \vspace {1mm}

A cornered Riemannian manifold $P$ is called  {\it preconvex}, if 
   the {\it dihedral angles}, denoted $\angle[Q_i\pitchfork Q_j]$, that are continuous functions 
   on the  edges of $P$,  that are the  $(n-2)$-faces $ Q_i\cap Q_j$,  are  $<\pi$  for all pairs of adjacent 
   $(n-1)$-faces $Q_i$ and $Q_j$.
 
 Observe that  preconvexity, does not depend on the metric and that 
 intersections of finitely many domains
 with smooth mutually transversal boundaries in smooth manifolds $X$ are preconvex.

\subsection {Surgery at the Corners.}

The simplest instance of surgery is the following operation of 

{\it Adding  $1$-handle at the corners.} Let $P\subset X=X^n$ be a polyhedral domain
and let $A^1\subset X$ be a smooth curve, (future axis of the surgery) that   joins two vertices $b_i$ and $b_2$ vertices (i.e. codimension $n$ faces)  in $P$ and that does not meet $P$ apart from these points.
 
Suppose there exists a face preserving diffeomorphism $D$ between small neighbourhoods   $U_1$ and $U_2$ of $b_1$ and $b_2$ in $P$. (If $X$ is a Riemannian manifold one thinks of these neighbourhoods as intersections of the Riemannian $\varepsilon$-balls $B_{b_1}(\varepsilon)$ and $B_{b_2}(\varepsilon)$ with $P$.)

Notice that such a diffeomorphism does exist if $P$ is cosimple at these vertices.

The boundaries $\partial U_1$ and  $\partial U_2$, that are identified by $D$,  carry natural  structure of cornered manifold of dimension $n-1$  call it $Q^{n-1}$. Then the curve  $A^1$ can be thickened
in  $X$ to $A^1\times Q^{n-1}$ and glued to $P\setminus (U_1\cup U_2)$ at its two (thickened) ends, where the resulting domain, call  it $P'$ carries a natural cornered structure. Notice that $P'$ has  the same vertices as $P$, except that $b_1$ and $b_2$ are not there anymore.

{\it  Example: Doubling.} Let $P$ consist of two disjoint diffeomorphic polyhedral domains.
Then upon performing the above surgery over all pairs of corresponding vertices one obtains a new domain $P''$ of depth $n-1$, i.e. with no vertices at all.

{\it Remark.} The above works as stated only if $dim(X\geq 3)$ where one may assume  the curves joining the vertices do not have to intersect, but if $n=2$ then $P'$ can be a priori 
immersed rather than imbedded into $X$, that is, however,  rather immaterial for our purposes.

Now, in general let $B^{m-1}\subset  P$ be an $(m-1)$-dimensional face in $P$ that contains no lower dimensional faces. For instance, this may be a a closed $1$-face in the above $P''$.

Let $A^{m}\subset X$ be an $m$-submanifold with $\partial A^{m}=B^{m-1}$ 
that meets $P$ only at $B^{m-1}\subset P$.  Let the normal corner structure of  $B^{m-1}\subset  P$
is represented by a flat fibration $U\to B^{m-1}$    with the  fiber $C^{n-m+1} $, that is the cone over a cornered manifold  $Q^{n-m} $.

(Such a $C^{n-m} $ automatically exists if  $B^{m-1}$ is connected; if $B^{m-1}$ is simply connected then, moreover,   its small neighbourhood in  $   P$ splits as  $B^{m-1}\times C^{n-m+1} $.)
  
Let $W^n\to A^m$  be an extension of  the flat fibration $U$ from $B^{m-1}$ to $A^m\supset B^{m-1}$.  Then one can remove $U$ from $P$ and attach $W$ instead.

{\it  Example: Multi-doubling.} Take two disjoint copies of the above doubled $P''$
and join the unions of their curves of edges by a (possibly disconnected)  cylinder $A^2$. Then the result of the corresponding surgery, say $(P'')''$  will have depth $(n-2)$. If we repeat this $n$ times   we arrive at a manifold  $P^\ast$ of depth zero -- that is a slowed manifold with no corners of any kind.

{\it Face  Suppression.}  If one double $P^\ast$ along the boundary, one obtains  a closed manifold  with "pure edge singularities" (as in section  2.3). 
\vspace {1mm}

{\it Remark.}  The above kind of surgery applies to all spaces with normally conical singularities. Thus, for example, every pseudomanifold $P$  can  be turned to a  manifold $P^\ast$ 
that is "kind of cobordant" to $P$.

 \subsection {Mean Convexity and Dihedral Angles.}  
Call  a  preconvex  Riemannian manifold  $P$ with corners  
 {\it mean curvature convex}, and write   $mn.curv(\partial P)\geq 0$ 
  if all  $(n-1)$-faces $Q_i=Q^{n-1}_i$ 
  have (non-strictly) {\it positive}   
mean curvatures,  i.e. the variations of the  $(n-1)$-volumes of the faces $Q_i$ are {\it non-positive}  
under infinitesimal
{\it inward} deformations of $Q_i$. 
 
 {\it Remarks and Examples.}  (a)  The simplest instance of  a mean convex domain in $X$  is a full dimensional submanifold $U\subset X$ with  a {\it smooth} boundary that has positive mean curvature. Such domains are abundant in $X$. For instance every piecewise smooth subset $Z\subset X$ of codimension $\geq 2$ in $X$ admits an arbitrarily small smooth mean convex neighbourhood. 

Furthermore, if $U$ is smooth mean convex, then the union $U\cup Z$ admits a 
 smooth mean convex neighbourhood in $X$.

(b) The  intersection of finitely many  domains $U_i \subset X$ with mutually
  transversal smooth mean convex  (e.g. convex)  boundaries $\partial U_i$ is an    m.c. convex
   polyhedral domain in $X$ with the faces made of pieces of $\partial U_i$.

(c)  If $P$ is compact and   the faces have strictly positive mean curvatures, then (it is obvious, see \cite {plateau-stein} for related results) the boundary of a an $\varepsilon$-neighbourhood of $P$, call it $ \partial P_{+\varepsilon}$, is $C^{1,1,}$ smooth for small $\varepsilon >0$  and has  (discontinuous ) positive mean curvature. 

It easily  follows that  a mean convex $\partial P$ can be approximated by $C^2$-smooth hypersurfaces with 
 $mn.curv>0$, unless  (some connected component of) $\partial P$  consists of a single face with zero mean curvature.

But  this is not especially relevant in the present context: we are keen at keeping track of the {\it combinatorial pattern} of the corners of $P$ and of the {\it dihedral angles} between adjacent $(n-1)$-faces at the corners.

\vspace {1mm}

 {\it The combinatorial type/scheme ${\cal CT }={\cal CT}(P)$} of  a  manifold $P$ with corners refers to 
  the intersection  and the adjacency patterns between its  codimension one faces  $Q_i$.

Observe, that the combinatorial type is stable under Cartesian products of $P$ 
by  connected manifolds (without boundaries)  with their tautological corner structures of depth zero.

{\it Combinatorial equivalence.} We often call  two domains {\it combinatorially equivalent} if they are of the same combinatorial type. Notice that such domains do not even have to be of the same dimension.

{\it  Cubical Domains $P$ of Depth $d$ in $n$-Manifold $X$.}   Such a {\it $d$-cubical domain $P$} has $2d$  faces $Q$ of codimension one, where every face  $Q$ has a unique
{\it opposite face}, call it  $-Q$, which does not intersect $Q$.  Every cubical $P$ admits a continuous map
onto the unit $d$-cube,     $(P,\partial P) \to(\square^{d}, \partial\square^{d})$, with the  faces of $P$ being the pullbacks of the faces of  $\square^{d}$, where, such a map is uniquely, up to homotopy in the class of faces-to-faces maps, is determined of by to which faces of the cube the faces of $P$ go.

 An  $n$-cubical $P$ is called {\it essential} if the map   $(P,\partial P) \to(\square^{n}, \partial\square^{n})$ has {\it non-zero} degree, where the degree is understood $\mod 2$  if $P$ is not oriented.

\vspace {2mm}

 \textbf{ Combinatorial  Mean Curvature Convexity   Problem.  } When does    
   a Riemannian manifold $X$, e.g. $X=\mathbb R^n$, 
  contain a mean curvature convex polyhedron  $P=P_{\leq\alpha_{ij}}$
of a given combinatorial type where all its dihedral angles are  bounded  by
$\angle_{ij}(P) \leq \alpha_{ij}$ for given constants $\alpha_{ij}$? \vspace {2mm}

This is unknown even for $X=\mathbb R^3$, where the answer is available only for {\it "normal" mean  curvature convex} domains  $P\subset \mathbb R^3$   that are   combinatorially  equivalent to  {\it prisms},  where "normal" means that  the dihedral angles at the top and the bottom  of these "prisms" are equal  $\frac {\pi}{2}$ and where a simple argument 
(see 
 section 5.4) shows  that \vspace {1mm}

{\it the dihedral angles $\alpha_i$ between the side faces of these  $P$ satisfy
  $$\sum_i(\pi-\alpha_i)\leq 2\pi \leqno {(\ast)}$$
with the equality {\it only} for convex prisms  $P\subset \mathbb R^3$ with flat faces.}

 \vspace {1mm}

The inequality $(\ast) $ remans true for "prisms" $P$ with Riemannian metrics with 
{\it non-negative scalar curvatures} (see section 5.4) and it suggests a possibility of defining general spaces $X$ with non-negative scalar curvatures, at least in dimension $3$, 
via this $(\ast) $. But it is unclear what kind of singularities one may admit in such a definition. \vspace {1mm}

 The simplest case  where we do not know the answer is that of {\it continuous} Riemannian metrics on smooth $3$-manifolds  $X$.
 
  Here, the mean curvature convexity of hypersurfaces may be   defined as  in section 4.4 by the existence of "many" localized {\it inward  deformations}  of faces that are  {\it area decreasing}. (This defines,  as stated,   {\it strictly} positive mean curvature.)
  
 Alternatively (that is  not, in general, equivalent)   one may requite {\it all} small localized {\it outward} deformations of the faces to be area  {\it increasing. }
  
 Or, one may require both: increase of area for {\it all} outward deformation and decrees for {\it  "many"}  inward  deformations. \vspace {1mm}

{\it Test Question.}    Suppose that all  {\it sufficiently small} prisms  $P=P_x\subset X $ that contain  $x\in X$ satisfy  $(\ast) $ for {\it all}  $x\in X$.  Do then all $P\subset X$ satisfy $(\ast)$?
 
 (This question is motivated by a similar  localization  property of the {\it  comparison   inequalities}  for geodesic triangles  in Alexandrov's spaces with positive curvatures.)

\subsection{Mean Convexity under Surgery.}

Let $P$  be a strictly mean curvature convex polyhedral domain in a Riemannian $n$-manifold  $X$ with all  dihedral angles of $P$  bounded  by
$\angle_{ij}(P) <\alpha_{ij}$. Let $P'$ be obtained by a surgery along some union   $B^{m-1} \subset P$ of  $(m-1)$-faces that  themselves have no corners (see section 1.1).

Observe, that in general when $2m\geq n$, such a $P'$ is immersed rather than embedded into $X$.
\vspace {1mm}

{\it If $n-m\geq 2$, then $P'$, with its combinatorial structure coming along with the surgery,  can be arranged (i.e. immersed) in $X$
with strictly positive mean curvature of all its faces  and with all dihedral angles satisfying $\angle_{ij}(P') <\alpha_{ij}$} 
 \vspace {1mm}
 
 This is immediate with (a simplified version of) the  argument from  \cite {gromov-lawson2}  for a similar behavior of positive scalar curvature under surgery. 
 
  \vspace {1mm}
 
 {\it Corollary.} The non-strict version of  the above combinatorial  mean curvature convexity   problem, that is of the existence of  $P=P_{<\alpha_{ij}}$, is invariant under the codimension $\geq 2$-surgery.
 
{ \it Example.} Every mean curvature strictly convex   $P=P_{<\alpha_{ij}}$ can be transformed by multi-doubling to a   $P^{''...}$
 without corners, i.e. of depth $2$, where   \vspace {1mm}
 
\hspace {-4mm}$mean.curv(P)>0$   \& $\angle_{ij}(P) <\alpha_{ij} \Rightarrow$  $mean.curv(P^{''...})>0$ \& $\angle_{ij}(P^{''...}) <\alpha_{ij}$.
 \vspace {1mm}
 
 {\it Question.} Does the reverse implication also hold true? 
 
Namely, does the existence of a position (i.e. of an immersed) $P^{''...}$ in $X$ with strictly positive mean curvatures of the faces and all dihedral angles  strictly bounded by given
$\alpha_{ij}$ imply the existence of such a position for $P$ in $X$?

 \subsection { Mean Curvature Stability and Semistability Problems.} 
 
 Conjecturally, the existence of  the above $P_{<\alpha_{ij}}$,  is stable under {\it   smooth perturbations} of the Riemannian metric $g$ in $X\supset P$  that are $\varepsilon$-small   in the {\it $C^0$-topology.}
 
 More generally, 
 let $P\subset X$ be a compact  strictly preconvex  (i. e. all  $\alpha_{ij}<\pi$)  polyhedral domain in a smooth (meaning $C^\infty$)   manifold $X$ with a $C^2$-smooth Riemannian metric $g$.
 Let $X'=X_\varepsilon = (X_\varepsilon,   g_\varepsilon)$ be another Riemannian manifold with $C^2$-smooth metric $g'=g_\varepsilon$ and let 
 $f_\varepsilon:X' \to X$ be a continuous map.
 
 An essential   example is where $dim(X') =dim (X)$  and $f_\varepsilon$ is a $e^\varepsilon$- {\it bi-Lipschitz} homeomorphism. Also we  allow  $dim(X') >dim (X)$ where 
  $f_\varepsilon$ is   $\varepsilon$-close in some sense to a {\it Riemannian submersion.}
 
We seek   conditions on  $X'  $   and on  $f_\varepsilon$ that would guarantee 
the existence of another map, say $f'_\varepsilon: X'\to X$  such that \vspace {1mm}

$\bullet_{\to0} $ $f'_\varepsilon$ is {\it close} to  $f_\varepsilon$, where close may mean $dist_X(f'_\varepsilon(x'), f_\varepsilon(x'))\leq \delta(\varepsilon)\to 0$ for $\varepsilon\to 0 $ and all 
$x'\in X'$ or that the function  $ x'\mapsto dist_X(f'_\varepsilon(x'), f_\varepsilon(x'))$  is small in a weaker  (e.g some $L_p$)  topology;

$\bullet_{reg} $ the pullback 
$P'=P'_\varepsilon=(f'_\varepsilon)^{-1}(P)\subset X'= X_\varepsilon$ is a smooth polyhedral domain
with the faces being the $f'_\varepsilon$-pullbacks of those of $P$.\vspace {1mm}

In order to formulate the stability and the semistability conditions, we agree that
 the mean curvature of a polyhedral domain at an "edge point" $x$ i.e. at a point where {\it exactly two} faces meet signifies

  \hspace {10mm} $\pi$ { \it minus the dihedral angle $\angle_{ij}$ between these faces at $x$.}\vspace {1mm}

Now the stability  and the semistability conditions read:

$\bullet_{stbl} $ 
$$ |mn.curv_{x'}(P'_\varepsilon)-mn.curv_x(P)| \leq \kappa(\varepsilon)\to 0 \mbox  {   for  $\varepsilon\to 0$, }$$
for  all $x'$ at the faces and at the edges of $P'_\varepsilon$ and     $x=f'_\varepsilon(x')$;

$\bullet_{semistbl} $ $$ mn.curv_{x'}(P'_\varepsilon)-mn.curv_x(P) \geq - \kappa(\varepsilon)\to 0 \mbox  { for  $\varepsilon\to 0$ }$$
for  all $x'$ at the faces and at the edges of $P'_\varepsilon$ and    $x=f'_\varepsilon(x')$.\vspace {1mm}

A particular instance of the latter  is where $P$ is strictly mean curvature convex and 
the same is required of  $P'_\varepsilon$.\vspace {1mm}

{\it On Regularity at the Corners.} The  conditions $\bullet_{stbl} $ and $\bullet_{semistbl} $ depend on  the faces of  $P'_\varepsilon$ being  $C^2$-smooth away from the edges and $C^1$-smooth at the edges but no regularity at the {\it corners,}  i.e. at the codimension$\geq 2$ faces is formally needed. 

This suggests  modified versions of  the stability and semistability problems where instead of $\bullet_{reg}$ we require only  {\it $C^\alpha$-H\"older}  smoothness of  $ P'_\varepsilon$  at the corners   for some $\alpha>0$.

This relaxed regularity condition is easer to satisfy when we  construct   $P'=P'_\varepsilon$
by means of the geometric measure theory (as we do it in sections 3 and 4).  

For instance,  let, say a cosimplicial,  curve-faced polyhedron  $P'$ be constructed "face by face"   where each $(n'-1)$-face  $Q'_i$ of $P'$ for $n'=dim (X')$ comes as a solution of  a {\it Plateau bubble} problem  (as defined below in section 3)
with  free boundary contained in  the union  $Q'_{\hat i}$ of the remaining faces.  
(Technically, the solution $Q'_i$ of  this  Plateau problem must be kept  away from the boundary of $Q'_{\hat i}$, with the resulting $P'$ obtained by chopping away  a small neighbourhood  of  $\partial Q'_{\hat i} \subset Q'_{\hat i}$.)

Then, probably(?), one can not guarantee
the $C^2$-smoothness of   $Q'_i$ at the edge points of  the hypersurface   $Q_{\hat i}$ (although it is likely such $Q'_i$ are $C^1$ at all boundary  points where their tangent cones are flat, and possibly,   $C^{1,\alpha}$-smooth, even with   $\alpha=1$) but  H\"older can be sometimes obtained. For instance, a {\it Reifenberg flatness argument}  delivers such  H\"older   stability for   $f_\varepsilon$ being a $e^\varepsilon$- {\it bi-Lipschitz} homeomorphism with a sufficiently small, depending on  $X$, $g$ and $P$, positive $\varepsilon$
(see sections 4.6-4.8).\vspace {1mm}


Insufficient smoothness does not seem to  excessively harm essential geometric applications of minimal hypersurfaces as well as of higher codimensional subvarieties. 

Moreover, this apparently remans so for more general singular ambient spaces, e.g. for Alexandrov's  spaces with curvatures bounded from below. For instance, it seems likely  that Almgren's sharp isoperimetric inequality indicated on p. 475 in \cite {singularities} 
for smooth manifolds with non-negative sectional curvatures extends to singular  such spaces $X$ (with the conical spaces $X$ being extremal as pointed out in \cite {singularities}).

\subsection {Dihedral Extremality and Plateau-hedra.}

 A mean curvature convex  polyhedron 
$P \subset  X$, e.g.  a convex polyhedron in $
\mathbb R^n$, is called
{\it  dihedrally  extremal} if no "deformation"   $P \rightsquigarrow P'$  of $P$
can  diminish its dihedral angles while keeping the faces mean curvature convex.
(Compare \S5$\frac{4}{9}$ in  \cite {positive}.)

That is,  more precisely, if a 
mean curvature convex polyhedron $P' \subset X$  of 
the same combinatorial type 
as  $ P$  has (non-strictly) smaller suprema of its dihedral
angles along  {\it all} $(n-2)$-faces  $Q'_{i}\cap Q'_{j} \subset P'$ than $P$,
$$\sup_{Q'_i\cap Q'_j} \angle[Q'_{i}\pitchfork Q'_{j}] =  \alpha'_{ij} \leq \alpha_{ij}=  
\sup_{Q_i\cap Q_j} \angle[ Q_{i}\pitchfork Q_{j}],$$
then, necessarily, $\alpha'_{ij}= \alpha_{ij}$.

\vspace {1mm}

{\it PP-hedra.}  An  polyhedral domain $P \subset X$ is called a  { \it (poly) Plateau-hedron}, or
 {\it PP-hedron}  if all its $(n-1)$-faces 
have {\it zero   mean curvatures} and the
 dihedral angle  functions  
 $\angle_{ij}=\angle[Q_i\pitchfork Q_j]$  are {\it constant} on the edges  that  are
   $(n-2)$-faces $ Q_{ij}=Q_i\cap Q_j$ for the pairs of  $(n-1)$-faces $Q_i, Q_j\subset P$ with 
   $dim(  Q_i\cap Q_j)=n-2$.

Notice that we {\it do not require}  preconvexity of $P$ (which is  equivalent to mean curvature
convexity in this case) but we will be dealing mainly with preconvex PP-hedra.

 \vspace {1mm}
 
 Basic instances of Plateau-hedra are ordinary polyhedral domains   with flat
  (i.e. totally geodesic) faces in 
 manifolds $X$ of constant curvature, e.g. in $X=\mathbb R^n$. 

Also it is easy to see that   \vspace {1mm}

{\it  dihedrally  extremal mean curvature  convex polyhedra $P$ in 
 Riemannian 
 
 manifolds are  Plateau-hedra.}  

\vspace {1mm}

{\it Singularities: Cones and Corners.} The above does not imply, however,  that all combinatorial types of polyhedra contain 
dihedrally extremal representatives, since the corresponding existence/regularity theorem is unavailable.

On the other hand,  one may attempt a construction of  Plateau-hedra  in 
 a Riemannian manifold $X$  by a variational  argument where  a face $W=W_i$ of a desirable $P$ is obtained as a solutions  of a 
 Plateau  type problem with free boundary, i.e. the boundary of  $W$ must be contained in the union $W_{\bot} = W_{\hat i}$ of the remaining  faces, see section 3 and 4.

Yet,  such a   $W$ may have {\it singularities}, both in the interior and at boundary points  $x$ in  $ W$.

For instance, if  $n-1=dim(W)\geq 7$, then $W$ may have {\it quasi conical singularities} at some points $x\in W$  where, by definition, a tangent cone  is {\it non-flat}.  
 
But  if $x\in int (W)$, 
 or if $x\in \partial W\cap reg (W_{\bot} ) $, i.e. if $x$ lies    away from the {\it edges} of $W_{\bot} $ as well as  from {\it interior singularities of the  faces}  of $W_{\bot} $    and if a tangent cone of $W$ at $w$ is flat, then $W$ is smooth at $w$, see   \cite {allard} \cite {gruter}, \cite {nirenberg}, \cite {gruter-jost}.\footnote { This was explained to me by Fang-Hua Lin.  Also I had a useful communication with Frank Morgan and   Brian White concerning these problems.}

  But the behaviour of $W$ at the singular points of $W_{\bot}$, even at the regular corners, i.e. where $\partial W$ meets edges between  {\it non-singular} faces in $W_{\bot} $, may be  more complicated, e.g. see \cite{simon}, \cite {brakke}.\footnote{I must admit I have not truly studied   these papers.}
 
 In particular, one has the following 
 
 {\it  Perturbation Question.} Let $P_0\subset \mathbb R^n$ be a cosimpicial  convex polyhedron.  Does it admit an arbitrary small perturbation to a Plateau-hedron $P'$
 with  {\it non-flat} faces? 
 
A natural approach here would be via a solution of the linearized problem combined with the implicit  function theorem, but one can not guarantee regularity 
at the corners.\footnote{ The existence of singularities for linear boundary value problems in  domains with  edges and corners  was pointed out to me by  Jeff Cheeger and  the full extent of the difficulty of this problem was explained to me by  Vladimir Mazia.}

\vspace {1mm} \vspace {1mm}

Inevitability of singularities  suggests a more general definition of  Plateau-hedra and of cornered domains in general  and  such a concept is also needed in the (conjectural)  context   of the  theory of singular spaces with  scalar curvatures bounded from below.  But our understanding of singular polyhedra, in particular of  singular Plateau-hedra, remains  unsatisfactory.

{\it Why at the Corners?}  The   simplest instance of singularities at the corners is that of minimal surfaces in $ Y\subset \mathbb R^3$ contained in intersection $P$ of two subspaces  with free  boundary $\partial Y\subset  \partial P$. If the (dihedral) angle $\alpha$ between the half plane that make the boundary of $P$ is $\pi/k$ for an integer $k$
then $Y$ extends by reflections to a  minimal hypersurface, say  $2kY$  around the edge in $P$.  Consequently, $Y$ is smooth, actually real analytic, in this case.

Furthermore, because of the $2k$-th order symmetry,  $2kY$, and hence $Y$ as well,  are   flat of order $k$ at the corner point   $y_{\lefthalfcup}\in Y$, i.e. where $Y$ it meets the edge in $P$.

But if $\alpha$ is incommensurable with $\pi$, the same symmetry argument shows that
if smooth, then $Y$ must be flat of {\it infinite} order at $y_{\lefthalfcup}$. In particular $Y$ can not be real analytic unless it is flat and, probably,(?)  it cannot be even $C^2$.

On the other hand if a {\it curve-faced} $P$ has the dihedral angle $\pi/k$ one expects a reasonable smoothness at the corner.

\vspace {2mm}

\subsection {Euclidean Dihedral Extremality  Problem. } 
  
\hspace {23mm}   \textbf {   [?]   Conv}exity $\Leftrightarrow$  \textbf {D}ihedrally  \textbf {Ext}remality \textbf {[?]}

 \vspace {1mm}

Namely,
 
 \vspace {1mm}
 
 [\textbf {CONV}$\Rightarrow$\textbf {DEXT?}] \hspace {3mm}  Which {\it convex} polyhedra $P\subset \mathbb R^n$ are
{\it dihedrally extremal}?  

  \vspace {1mm}

  [\textbf {DEXT}$\Rightarrow$\textbf {CONV?}] \hspace {3mm}   Are there {\it  non-convex  
   dihedrally extremal} $P\subset \mathbb R^n$?

 \vspace {1mm}

 My guess is that  [\textbf {DEXT}$\Rightarrow$\textbf {CONV }] for {\it all} combinatorial types $\cal CT$ of $P$, 
 i.e. all m.c. dihedrally  extremal  m.c. convex polyhedral domains in $\mathbb R^n$ 
 are convex, even if we allow  the above mentioned singularities.
  
  On the contrary,    
 dihedral extremality   seems too good to be true for {\it all} convex 
 polyhedra  $P \subset \mathbb R^n$.  
  
  In fact,  even if $P' \subset \mathbb R^n$ is a 
{\it convex} (not just mean curvature  convex) polyhedron, combinatorially equivalent to $ P$,
it is unclear why   
the dihedral angles of $P'$ can not be all  {\it strictly smaller} 
than the corresponding angles in  $ P$. 
This {\it can not} happen for simplices $P=\Delta^n\subset \mathbb R^n$ 
by {\it Kirszbraun theorem} applied to the  {\it dual} simplices $(\Delta^n)^\ast$

 \vspace {1mm}

 Possibly, simplices and their Cartesian products are  dihedrally extremal.
 Also we have the following

\vspace {1mm}

\textbf {  $ \triangle \times \square  \times \pentagon\times 
\hexagon\times...$-Conjecture.} The Cartesian products 
$P\subset \mathbb R^{2k+m}$ of  $k$ convex polygones $\subset \mathbb R^2$ 
and an $m$-cube $\subset \mathbb R^m$ 
are m.c. dihedrally extremal: one {\it can not diminish} the dihedral angles by "deforming" such a $P$ 
{\it without} developing {\it negative} mean curvature in some of the faces.

\vspace {1mm}

\textbf{Theorem-Example:   $\triangle   \times\square\times \hexagon$-Extremality.} The simplest 
case where we settle this conjecture is for Cartesian products of $m$-cubes with 
Cartesian powers of  {\it regular} triangles and  regulate hexagons,
$$ P=  \square^{m_1}\times   \triangle^{m_2} \times \hexagon^{n_3} \subset
    \mathbb R^{2m_3+m}.$$
 In fact, (see "Gluing  around Edges" in section 2) this  
  $\triangle   \times\square\times \hexagon$--extremality   follows
  from the positive   
solution to the { \it Geroch conjecture} on  {\it non-existence}  of metrics of {\it positive scalar
curvatures} on tori $\mathbb T^n$.  

Strangely enough, there is no apparent {\it direct elementary} proof
of this apparently  {\it intrinsically   Euclidean} inequality/extremality.
\vspace {1mm}

There are, recall,  two approaches to the  Geroch conjecture.
The original one, due to Schoen and Yau, depends on smoothness of  minimal hypersurfaces
 $H^{n-1} \subset \mathbb T^n$ 
and applies  {\it only} to $n\leq 7$. Their argument easily extends to  $n=8$   by   {\it non-stability of singularities} of $7$-dimensional minimal hypersurfaces  \cite {smale}, while    a way
around  singularities  for $n\geq 9$,  found   
relatively recently  by Lohkamp \cite{lohkamp},   is rather intricate. (Possibly, one can prove the full  $ \triangle^{m_3}  \times \pentagon^{m_5}\times 
\hexagon^{m_6}...\times \square^{m/2}$-conjecture utilizing Lochkamp's techniques.)

Another proof (see \cite {gromov-lawson1}), that 
depends on the index theorem for  twisted Dirac  operators,  indiscriminately applies 
to all $n$, but it needs the {\it spin structure}. 
(This  causes no problem for $\mathbb T^n$ but becomes  a hurdle for non-spin manifolds 
$X$.)

 We shall interpret applications of  these methods to corned domains $P$  as   "billiard games" played  by  Dirac operators and minimal hypersurfaces in $P$ (see section 2.3).

\vspace {1mm}

{\it Extremality, Rigidity and Scalar Curvature.} The concepts of the dihedral extremality  as well as the extremality of manifolds with positive scalar curvature that was   studied in
\cite {positive} \cite {min-oo2} \cite{listing} \cite{goette}  are embraced by the following 

{\it Definition}. Let  $X =(X,g)$ be
a Riemannian $n$-manifold  $X$ with corners, let $f$ be a smooth map  of {\it non-zero degree} of another smooth  $n$-manifold $Y$ onto $X$ and let us endow $Y$ with a corner structure induced by $f$ from $X$. (One may admit $Y$ with $dim(Y)\geq dim(X)$ with a suitable concept of degree for maps $f:Y\to X$.)

$X$ is called {\it extremal} if every for every such map $f$  every Riemannian
Let $Y$ be endowed with a Riemannian metric $g_Y$   such that

 (1) $f$ is {\it distance decreasing} (sometimes one only requites $f$ being area decreasing);

(2) $f$ {\it decreases} the scalar curvature: 
$$scal_X(f(y))\leq scal_Y(y)\mbox {  for all $y\in Y$;}$$

(3) $f$ decreases the mean curvature of the $(n-1)$-faces: 
$$mean.curv(f(y)) \leq  mean.curv(f(y))\mbox {  for all $y$ in all  $(n-1)$ faces of $Y$;}$$

(4) $f$ decreases the dihedral angles at all edgers 
$$\angle_{f(y)} \geq \angle_y \mbox {  at all
$y$ in all $(n-2)$-faces of $Y$;}$$

Then $X$ is called {\it extremal}, if for {\it no}  such $Y$, $f$  and  $g_Y$ 
the  increases/decreases in the inequalities   (2)-(4) {\it can be strict}, {\it not  at a single point  $y\in Y$} while rigidity says, not quite precisely, that  all $Y$ satisfying (1)-(4) are obtained from $X$ by "obvious modifications".

Notice that (1)-(3) can be unified if the $scal_X(x)$ is understood as the mean curvature at
the points $x$ in the (interiors of) $(n-1)$-faces of $X$ as $\pi-\angle_x$ at $x$ in the 
$(n-2)$-faces.  

\vspace {1mm}

(A) Several extremality/rigidity results are  available  for  closed (i.e. of depth 0) manifolds with positive scalar curvature, in particular for most (all?) compact   Riemannian symmetric spaces $X$, see \cite {min-oo2}, \cite{listing}, \cite{goette},  \cite {goette2}, \cite{llarull}, \cite {positive} which is proved with Dirac operators. \vspace {1mm}

 {\it Can this  be proved by means of  minimal hypersurfaces or, rather, of $\phi$-bubbles  (see 3.1)?} \vspace {1mm}

 The extremality/rigidity of the round spheres \cite{llarull} implies  that  convex metric balls   in simply connected spaces  $X$ of  constant curvatures ("convex" is relevant if $curv(X)>0$, i.e. $X=S^n$)  are extremal  for $n=dim(X)\leq 8$. This follows by the wrapped product  argument from \S5$\frac{5}{6}$ in \cite {positive}, where the  case $n=8$   relies  on non-stability of (isolated) singularities of  $\phi$-bubbles in $8$-manifolds   (as well as of 
 minimal hypersurfaces \cite {smale})  and where, possibly, the extremality, but not, a priori,  rigidity may be obtained
 with \cite {lohkamp}. 
 
This rigidity is reminiscent of  the  generalized positive mass theorem  \cite {min-oo2}
and suggests a possibility of proving this extremality/rigidity by the Dirac operator method.

(B) The simplest examples of extremal/rigid $P$ are  convex  $k$-gons in surfaces of positive (not necessarily constant) curvatures,
where their extremality and rigidity follows from the Gauss-Bonnet theorem.  \vspace {1mm}

{\it Are Cartesian products of extremal/rigid manifolds, in particularly of those in the above examples  (A) and (B), extremal/rigid?} \vspace {1mm}

We prove in this regard  the rigidity of $3D$-prisms ($k$-gons $\bigtimes  { } [0,r]\subset \mathbb R^3$) in section 5.4.

Probably, {\it the rigidity of Cartesian   products of $k$-gons}  (at least in the spin case) follows by extending the methods of
\cite {goette2} to manifolds with singularities along codimension two (divisor-like)  subvarieties, where the relevant  examples are   Cartesian products of  surfaces with isolated conical singularities.  \vspace {1mm}

(C)  Some (but not all)  {\it warped products}  of Riemannian manifolds may be  extremal/rigid. 

For instance, let $X=Y\times [0,R]$, where  $Y\subset \mathbb R^{n-1}$ is  a convex polyhedral  domain and where the warped product  metric    $g(y,t)= a(t)^2 g(y)+dt^2$
 on $X$ has  constant negative curvature, namely,   $a(t)=e^t$. \vspace {1mm}
 
{\it Does    extremality/rigidity  of $Y$ imply  that of   $X$?} \vspace {1mm}

In fact, the rigidity of $X$ is proved in  \cite {positive} for $Y$ being a flat torus of dimension $n-1\leq 6$, where
the extremality (but not rigidity) for  $n-1=7$ extends with \cite {smale} and, probably,  with \cite {lohkamp} to all $n$.

The  $\phi$-bubble  argument from  \cite {positive}  can be combined with gluing around the edges (see 2.1) thus proving the extremality/rigidity of these warped product  $X$ for {\it  Euclidean reflection domains $Y$}, at least for $dim(Y)\leq 6$.

\vspace {1mm}
{\it Can one prove the rigidity/extremality of these warped products by a pure Dirac operator method in the spirit of \cite {min-oo1}?} \vspace {1mm}

(D)  Annuli $X=S^{n-1}\times [0,r]$  between concentric spheres in manifolds of constant curvature are instances of 
warped products which have geometric properties similar to but different from our extremality/rigidity. 

Such properties are proved in \cite {positive} by means of $\phi$-bubbles that limits the results to $n\leq 7$ (extended to   $n\leq 8$ with \cite {smale} and, possibly, to all $n$ with \cite {lohkamp}). 

Similar properties may be true for some (e.g. reflection, that is no big deal) domains  $P\subset S^{n-1}$ but the overall picture is far from clear. \vspace {1mm}

Observe, finally, that if a  certain space $X$  (without mean  convexity points at its boundary) of {\it negative} scalar curvature is  "extremal", this extremality must be opposite to what we saw above: when one enlarges such an $X$ its scalar curvature tends to increase rather then decrease. 
 
In other words,  the distance {\it decreasing} condition for maps $f:Y\to X$ is satisfactory restrictive  if $scal (X)\geq 0$ but  it seems more logical to require $f$ to be {\it distance increasing} (which needs to be properly defined for non-injective maps) at the points where $scal\leq 0$.

Ideally, one wants to prove  extremity and rigidity relative to maps  $f:Y\to X$ that  {\it decrease the  integrals} of the  scalar curvatures over some class of surfaces  in the two manifolds, something like
$\ \int _S scal(Y)ds \leq  \int _Tscal(Y)dt$ for $T\subset Y$ and $S=f(T)\subset X$ for the same kind of surfaces of $S$ and $T$    as in the semiintegral inequalities in section 5.4 of the present paper and  and in  0.5.C of  \cite {foliated}.

\subsection{ On  Acute Polyhedra.}  

Besides products  of $k$-gons, there is another class of 
"elementary" polyhedral domains where one may expect extremality/rigidity results.

 Call  a  Riemannian manifold $P$ with corners (non-strictly)   {\it acute}
if all its  dihedral angles are acute, i.e.  bounded by $\pi/2$,
$$\angle[Q_i\pitchfork Q_j] \leq \pi/2\mbox { for all pairs of adjacent $(n-1)$-faces $Q_i, Q_j \subset P$.}$$

{\it  Acute Spherical Polyhedra.} If a  convex   {\it spherical polyhedron} $P\subset S^n$ 
(i.e. an intersection of hemispheres) has {\it acute} dihedral angles then it 
is a {\it simplex}, or, in the degenerate case, the {\it spherical suspension} over a simplex
in $S^{n-i}\subset S^n$.

 Indeed,  the dual polyhedron, say $P^\perp_p \subset S^{n}$, 
    has all its edges longer than $\pi/2$.  Consequently, the distance between 
    {\it every two} vertices in  $P^\perp_p$ is $\geq \pi/2$; hence, there are  at most $n+1$ vertices
     in  $P^\perp_p$.    

It follows, that 

\vspace {1mm}

acute Riemannian manifolds $P$ with corners are  {\it simple} --  
there are exactly $n-i$   faces $Q_i^{n-1}$ transversally meeting  along every $i$-face
   $p \in \partial(P)$. 

\vspace {1mm}

Also observe that (non-strictly)  acute spherical triangles $\triangle \subset S^2$ have all their edges
bounded in length by $\pi/2$. It follows that all $m$-faces, $m=2,3,..., n-1$, of acute spherical
 $n$-simplices are acute.

\vspace {1mm}

   {\it   Acute Euclidean Polyhedra.} Cartesian products of {\it acute} simplices $\triangle^{n_i} \subset \mathbb R^{n_i}$,
   $$P=\triangle^{n_1}\times\triangle^{n_2}
   \times...\times \triangle^{n_k}  \subset \mathbb R^{n_1+n_2+...+n_k},$$
   are, obviously, {\it acute}.
  
   Conversely, 
 
 \vspace {1mm}
 
 {\it every acute  polyhedron $P\subset \mathbb R^n$  is a Cartesian product of simplices.}
  
   \vspace {1mm}

  This is easy and, certainly, has been known for ages.   But  I could not find this on the web  and  wrote down   a (few lines)  proof in \cite {hilbert}.
  
         \vspace {1mm}      
        
   { \it Questions}. Are all   acute  polyhedron $P\subset \mathbb R^n$ dihedrally rigid or at least external?

    What are possible combinatorial types of  mean curvature
    convex   Riemannian  cornered manifolds  $P$ with  $scal(P)\geq 0$ and
     acute dihedral angles? 
     
     Are there any constrains on the combinatorial types of mean curvature convex  
     $P\subset \mathbb R^n$ 
     and 
      $\angle[Q_i\pitchfork Q_j]\leq (\pi/2)+\alpha$ for a given  $0<\alpha <\pi/2$?
  
  (I  stated in my article {\sl Hilbert Volume in Metric Spaces}  that there are only finitely many combinatorial types of  convex polyhedra with 
  
  \hspace {26mm} $\angle[Q_i\pitchfork Q_j]\leq (\pi/2)+\alpha$,   $\alpha<\pi/2$, 
  
  \hspace {-6mm} but  
   Karim Adiprasito recently  showed me counterexamples starting from dimension $3$.)


On the other hand, the  scalar curvature can be made 
arbitrarily large by  multiplying  any (compact)   $P^{n-2}$ by a small $2$-sphere,
where $P^n= P^{n-2}\times  S^2(\rho)$  has the same  dihedral angles
and mean curvatures of the faces as $P^{n-2}$, but these   have "rather  degenerate" combinatorial 
types.

 One  may expect  that  some (most?)  manifolds  $P$ with corners (convex 
 polyhedra?)  support no
  metrics $g$ with $scal(g)\geq 0$, with     $mn.curv_g( Q_i)\geq 0$ and with
   $\angle_g[Q_i\pitchfork Q_j]\leq (\pi/2)+\alpha$ for  every $\alpha < \pi/6$
    and, less likely, for  $\alpha \geq \pi/6$.
  
  Yet,  finding a single such $P$ for  any $\alpha>0$ remains problematic.\vspace {1mm}
  
Also we can not solve the following 

{\it Simplex Problem.} Let $P\subset \mathbb R^n$ be a curve-faced  polyhedral domain that is combinatorially equivalent to the $n$-simplex and let $\alpha_{max}(P)$ denote the supremum of its dihedral angles at all edge points. Notice that if $n\geq 3$ then, obviously, $\alpha_{max}(P)>\pi/3$. \vspace {1mm}

What is the infimum of  $\alpha_{max}(P)$ over all {\it mean curvature convex} $P$?\vspace {1mm}

Conjecturally,   this $inf\alpha_{max}$ is assumed by (the dihedral angle  between a pair of faces of ) the {\it regular  $n$-simplex with flat faces}, but it is not even a priori  clear if  $inf\alpha_{max}$ is  {\it strictly} greater than $\pi/3$ for  $n\geq 3$.

\subsection{ $C^0$-Limits of Metrics with $scal\geq \kappa$.} 

Our study of mean curvature convex polyhedral domains in 
Riemannian manifolds $X$, even for $X=\mathbb R^n$, is
 intimately related to the scalar curvature. For example we shall prove the  extremality of the above $P$
 in the class of all {\it spin} manifolds $P'$ with corners which have   $scal(X)\geq 0$ 
 by utilizing minimal
   hypersurfaces along with
 Dirac operators. 
  
  We also achieve   this for {\it non-spin} manifolds 
   $X$ with
  $dim (X)\leq 9$, where  the singularities of minimal hypersurfaces  are at most
  $1$-dimensional (actually, 
  we shall need this for "Plateau bubbles" in $X$, see 5.3)  and, as  we will show,   they  do not "feel" spin obstruction
 that lives in dimension $2$. 
   Probably, the analysis
  of singularities developed  in \cite {lohkamp} would allow  a direct (with no use of spinors) proof
   for all $dim(X)$.  \vspace {1mm}

Most current  results on manifolds $X$ with  $scal(X)\geq 0$  rely on {\it global}
techniques and do not tell you much 
on the geometry of small (but not infinitesimally small) 
and moderately large regions
 $U\subset X$. For 
example, the Dirac operator  can be directly used (almost) exclusively on {\it complete } manifolds $X$ (an exception is Min-Oo rigidity theorem for the  hemisphere \cite {min-oo2})  while
 the  Schoen-Yau approach  depends on a presence of
  {\it closed/complete} (or  "quite large" as in  \cite {gromov-lawson3}) minimal hypersurfaces  $Y^{n-1}\subset X=X^n$, similarly  to  how 
 the proof  {\it Synge's theorem} for  manifolds $X$ with
 $sect.curv(X)>0$ uses {\it closed} geodesics in $X$.

  Sometimes, one can derive  semi-global results  from global ones,  either by extending
a  metric from  a manifold $X$ with a boundary (or such a domain $U\subset X$) to  a complete $X_+\supset X$ keeping
 $scal(X_+)\geq 0$ \cite {wang} or by exploiting 
   {\it Plateau "soap" bubbles}  $ Y^{n-1}\subset X$ \cite {positive} to which global techniques apply.

Yet, all this fails short of  Alexandrov's approach to spaces  with $sect.curv\geq 0$ 
(and more generally with $sect.curv \geq - \kappa^2$)  via {\it (comparison) inequalities for 
 angles of  geodesic triangles}
that  indiscriminately hold on {\it all scales} and provide a non-trivial
 information on the geometry of all domains  $U\subset X$, be they big or small.

Hopefully,  lower bounds on dihedral angles of extremal PP-hedra  $ U\subset X$ may
play  a similar  role for $scal(X)\geq 0$. This, in turn, 
points  toward an Alexandrov's type  theory of {\it singular} spaces $X$
with $scal(X)\geq 0$ and, possibly, with   $scal(X) \geq - \kappa^2$.

Notice that there is an  analytic approach to singular metrics with positive scalar curvature  understood in the distribution sense in \cite {grant} and somewhat similar in    \cite {lee} but these do not seem to apply to our situation.

We do not know what the theory of objects (spaces?) with positive scalar curvatures
understood in the distribution sense  should be but we  prove in section 4.9 the following
\vspace {1mm}

\textbf {$C^0$-Limit Theorem.} (Compare \S5$\frac{5}{6}$ in  \cite {positive}) {\it  Let a smooth Riemannian metric $g$ on a  Riemannian manifold $X$ equal the uniform limit of smooth metrics  $g_i$ 
with  $scal_{g_i}(x)\geq  \kappa(x)$ for a  continuous function $\kappa$ on $X$. 
Then $scal_g(x)\geq \kappa(x)$  as well..}

\vspace {2mm}
{\it Remarks, Questions, Speculations. }  (a)  The above is a local property of metrics  and the general case trivially follows from that where $\kappa$ is constant. 

The  $C^0$-limit property for $\kappa=0$ is derived  from the existence of  particular (small) strictly mean convex cubical domains  with acute dihedral angles in manifolds with $scal<0$ (see $(\smallsquare)$ in 4.9) and the solution to the Geroch conjecture on non-existence of metrics with $scal>0$ on tori, while 
 the cases $\kappa>0$ and $\kappa<0$ reduce to $\kappa=0$ as follows.

First, let $\kappa=n(n-1)$ for $n=dim(X)$, where observe this $\kappa$ equals the scalar curvature of the unit  Euclidean $n$-sphere $S^n$. 
 
Given a metric $g$ on $X$ let $\check X=(X\times\mathbb R_+, \check g)$ be the  standard Riemannian/Euclidean cone over $(X,g)$, that is $g(x,t)=   t^2 g(x) +dt^2$.

Notice that if  $scal_g(x)=n(n-1)$ at a point $x\in X$, then $scal_{\check g}(x,t)=0$ and if  $scal_g(x)<n(n-1)$, then  $scal_{\check g}(x,t)<0$ for all $t>0.$
 
 Thus, the $C^0$-limit theorem   for metrics on $X$ with
$\kappa=n(n-1)$,  hence, for all  $\kappa>0$, reduces to that for $\kappa=0$ on $X\times \mathbb R_+$. \vspace{1mm}

(If $X$ is a compact manifold with $scal>0$  then the  double $2\diamond \check X$
of the cone  $\check X =X\times \mathbb R_+$ at the vertex corresponding to t=0 admits a  complex metric of positive curvature that is conical at both ends of this double. This suggests that  geometric properties of such $X$ can be expressed and/or generalized in terms of asymptotic geometries of complete manifolds where   Witten and Min-Oo style spinor arguments may be applicable. Can one, for instance,  derive Llarull's sphere rigidity theorem along these lines?) \vspace{1mm}

Now  let $\kappa<0$ and proceed similarly albeit more artificially. Namely,
let $\hat X =(X\times [-\delta, \delta], \hat g) $ for $\hat  g(x.t) = C_\kappa t^2g(x)+dt^2$, where the constant $C_\kappa>0$ is chosen such that 
$scal_g(x)=\kappa\Rightarrow scal_{\hat g}(x,0)=0$ for a given $\kappa<0$.

Then the inequality $scal_g(x)<\kappa$ implies  $scal_{\hat g}(x,0)<0$ and  the $C^0$-approximation theorem on $X$ with $\kappa<0$ is thus reduced to that on $\hat X$ with $\kappa=0$.

(The natural  cone metric  in  $X\times \mathbb R_+$ where $\kappa(X)<0$
is a Lorentzian one to which our flat Riemannian argument does not(?) apply; yet,  \cite {min-oo1} suggests an approach to this metric.)

Conclude by noticing that   if one is willing to add {\it two} (or more)  extra dimensions one can reduce the case of  $\kappa\neq 0$ to that of $\kappa=0$  by taking the Riemannian product  $X\times D^2$ for 
 a disc $D^2$ with $scal(D^2)=-scal (X)$.

Thus, the relevant information concerning $X$ may be   seen in the geometry
of the torus $T^{n+2}$  that is obtained by  {\it reflection development} (orbicovering, see 2.1 ) of a small cube in 
$X\times D^2$  (see ($\smallsquare$) in section  4.9)  where the sections of this cube by $X\times d$, $d\in D^2$ make an $n$-dimensional foliation. 

Can one abstractly define and study this kind of foliations with no reference to any  external $D^2$?

\vspace {1mm}


(b) The $C^0$-limit theorem  contrasts with 
{\it Lokhamp's  $h$-principle}: \vspace {1mm}

{\it \hspace {4mm}every Riemannian metric  $g$ on a smooth manifold $X$ can be  $C^0$-approximated  

\hspace {4mm}by metrics $g_i$ with $scal(g_i)<0 $ and even with $Ricci(g_i)<0$.}
  \vspace {1mm}

 (c) If $n=dim(X)=3 $, then the $C^0$-limit theorem for all $\kappa$, be they positive or negative,
 follows, from the {\it Gauss-Bonnet prism inequality} in 5.4. and a version of $(\smallsquare)$ from 4.9 (where it is used for the case $\kappa=0$) adjusted to  $\kappa\neq 0$.  \vspace {1mm}
 
 (d) The $C^0$-limit theorem for  $\kappa<0$ can be also proven in a more natural fashion intrinsically in $X$ itself similarly to the case $\kappa=0$   with a version of $(\smallsquare)$ for bands around   (germs of) suitable {\it convex hypersurfaces} $Y$ in $(X,g) $  with induced metrics having scalar curvatures  zero (or, rather, close to zero)   by   reproducing   the argument presented at  the end of 
  \S5$\frac{5}{6}$  in \cite {positive} with {\it small  cubical} $P\subset Y$ instead of $(n-1)$-tori as in  \cite {positive}.
  (These $P$ are similar to hyperbolic $n$-dimensional "prisms" that are suspended over  reflection domains in $Y^{n-1}_{horo}$ discussed toward the end  of section 5.4.)  \vspace {1mm}

(e) The $C^0$-limit  theorem is reminiscent of  {\it Eliashberg's 
$C^0$-closeness theorem for  symplectic diffeomorhisms}.

Is there something in common between the two theorems besides superficial similarity?

Is there a scalar curvature counterpart of  {\it Hofer metric} between Hamiltonian diffeomorphisms and/or Lagrangian submanifolds?

Is there anything interesting  (besides a specific  $K$-area inequality from \S5$\frac{4}{5}$ in  \cite{positive})  in geometry of quasi-K\"ahlerian metrics with positive scalar curvatures on symplectic manifolds? \vspace{1mm}
 
Also, the $C^0$-closeness of $scal\geq \kappa$ resembles  {\it Novikov's theorem
on the topological invariance of Pontryagin classes}, but it is  equally unclear if there is something profound behind this similarity. 


\subsection { Rigidity Problems around  $scal\geq 0$.}

Let a {\it continuous}  Riemannian metric $g$ on a closed manifold $X$ admits
a $C^0$-approximation by smooth metrics $g_i$ with

\hspace {-6mm} $scal(g_i)\geq -\varepsilon_i\to 0$, $i\to\infty$.

$\bullet_1 $ {\it Does $X$ admit a smooth metric  $g'$ with $scal(g') \geq 0$?}

$\bullet_2$ {\it Suppose that $X$ admits a continuous map to the $n$-torus, $n=dim X$, of non-zero degree. Is then the metric $g$  itself necessarily  flat?}\vspace {1mm}

One asks  similar questions concerning Dirac operators $D=D_g$, say on spin manifolds $X$:

{\it Is the spectrum of $D^2=D^2_g$ semi-continuous  for the  $C^0$-topology in the space of  smooth Riemannian metrics $g$?}\vspace {1mm}

In particular, 
let $g_i\underset {C^0}\to g$, where $g_i$ are smooth metrics with positive squared  Dirac operators $D^2_{g_i}$. 

{\it If $g$ is smooth, is then $D^2_g$ also (non-strictly) positive? }\vspace {1mm}

{\it What are  classes of continuous metric $g$ on an $X$ where positivity of $D^2_{g}$
makes sense?}
(A particular   class of interest  is that of {\it piece-wise smooth metrics} on closed manifolds obtained by gluing  Plateau-hedra along their faces, where  the approach from \cite{miao},  \cite{mcferon}, \cite{lee} and  \cite{grant} may be relevant.)

\vspace {2mm}

Also one may raise   such questions for other geometric operators, e.g. for the {\it coarse Laplacians}
in the bundles associated to the tangent bundle, where  "spectral $C^0$-semi-continuity"  of such operators 
is related via the {\it Feynman-Kac formula} to the following  question.
 \vspace {1mm}

{\it $C^0$-Semi-continuity of Holonomy.} Let $g_i\underset {C^0}\to g$, where $g_i$ are smooth 
metrics on the $n$-dimensional manifold  $X$ with their holonomy groups  contained in a given closed subgroup $H\subset O(n)$
in the linear group acting on the tangent space $T_{x_0}(X)$. Is the holonomy group of $g$ also contained in $H$? 

This is, probably, not hard to prove for "classical $H$", e.g.  for $H=O(k)\times O(n-k)$,  where the corresponding manifolds $(X,g_i)$ split;  thus, carry many flat geodesic subspaces. One may  approach  in a similar fashion  holonomy groups of symmetric spaces and, possibly, of K\"ahler manifolds.

Is there a metric characterization of  {\it  K\"ahler manifolds} that is   stable under $C^0$-limits of metrics similar to that in \cite {kahler} but that, unlike \cite {kahler} would  apply to non-necessarily closed manifolds, e.g. to (small)  domains in projective algebraic
manifolds?

\subsection {Singular Spaces  with $scal\geq 0$ and Related Problems.}

A  potential pool of singular $X$ with $scal(X)\geq 0$ spreads  immensely wider 
than that or than  the class of the Alexandrov spaces with
$sect.curv \geq 0$. In particular, this "pool" must include:


$\bullet $  spaces partitioned  into "PP-hedral cells" with essentially conical singularities  where the local geometry
is similar to that
 for   $sect.curv \geq 0$,

$\bullet $ certain  spaces 
with "fractal singularities".  \vspace {1mm}

Here are two such examples.  

(1) Let  $scal(X)\geq 0$ and  $U \subset X$ be a  
domains $U \subset X$ where the boundary $\partial U$  of $U$, that is allowed to have 
singularities,  has  {\it non-negative} mean curvature, e.g.  
 where  $\partial U$ comes as  a minimal hypersurface in $X$.  Then the double of $U$ along 
 $\partial U$ must be regarded as a space with $scal.curv\geq 0$. In fact, such an $X$ is often
 (always?)  equals Hausdorff limit of smooth manifolds with  $scal\geq 0$

(2)  Let   $X_1$ and $X_2$ be  $n$-manifolds   and $Y_i\subset X_i$, 
$i=1,2$ be  submanifolds for which there exists a diffeomorphism $Y_1\leftrightarrow Y_2$
that induces an isomorphism of their respective normal bundles.

This delivers 
 a diffeomorphism between the boundaries  of
 the normal neighbourhoods $U^{nrm}_i \supset Y_i$, call it
 $\chi : \partial U^{nrm}_1  \leftrightarrow \partial U^{nrm}_2 $, and let 
 $X=X_1\#_\chi X_2$ be obtained by gluing $X_1\setminus U^{nrm}_1$ with $X_2\setminus U^{nrm}_2$ 
  according to $\chi$.

If $scal (X_i) >0$  and $codim(Y_i) \geq 3$, then this $X$ carries a  canonical class 
 of  metrics with $scal>0$ 
that  equal
the original ones on $X_i\setminus  U^{nrm}_i$.

This gluing operation can be repeated infinitely many times and the resulting limit spaces
 should be regarded as having $scal\geq 0$ with 
the simplest instance of this is as follows.

 Let $X_0=X_0^n$, $n\geq 3$, be a compact manifold with 
$scal(X_0) >0$ and $\{y_1,y_2, y_3\} \subset X$ be a 3-point subset. Let us attach to each of 
these points a copy of $\zeta\cdot X_0$, that is $X_0$ with the metric scaled by $\zeta>0$,  where the point in  
 $(\zeta\cdot X_0)_i$ that corresponds to $y_1 \in X_0$ is attached to $y_i\in X_0$, for each $i=1,2,3.$ 
The resulting manifold $X_1$ (that is the connected sum of $X_0$ 
and three copies of $\zeta\cdot X_0$) has six "free"
points  corresponding to the unused counterparts of $y_2$ and $y_3$ in the three $\zeta\cdot X_0$.  

Let   $X_2$ be obtained by  attaching a copy of $\zeta^2\cdot X_0$ to $X_1$ at  each of these points 
and then, similarly, we   get $X_3$,  $X_4$, etc.. 

 If $\zeta>1$, the the Hausdorff   limit, call it  $(1-\zeta)^{-1}\ast X$,
 of the resulting  sequence of spaces $X_0, X_1, X_2,...$ is
   a smooth complete non-compact  Riemannian  manifold. If   $\lambda<1$, this is a compact self-similar fractal  space $X$
  which behaves  in many respects as nicely as  Riemannian
 manifolds do (e.g. it may have  essentially Euclidean {\it filling inequalities} see section 4.7) 
 and it definitely must be regarded as having $scal(X)>0$.  

Also notice that if $\zeta^n <1/2$, $n=dim(X_0)$, then  the volume of $(1-\zeta)^{-1}\ast X$, that is $\approx \sum_i 2^i vol(\zeta ^i\cdot X_0)=vol(X_0)\sum_i 2^i \zeta ^{ni} $, is finite.

\vspace {1mm}

{\it Remarks.} (a)  The class of spaces with $scal\geq 0$, unlike that with $sect.curv\geq 0$,
 can not be stable under  Hausdorff limits, unless extra   strong "topological  non-collapsing" conditions are 
 imposed on the spaces involved.  
 These conditions can be enforced  in the above case where $X$ equals the
 Hausdorff-Lipschitz  {\it projective limit} of $X_i$  for 
 naturally defined  uniformly  Lipschitz maps
  $X_{i+1}\to X_i$. 
 
 (The classes of spaces considered in \cite {bernig}, \cite{miao},  \cite{mcferon}, \cite{lee} and  \cite{grant} do not seem to be stable under geometric limits.)

 (b) {\it On $scal\geq \sigma$.} 
   A lower bound on the scalar curvature by a  positive or negative constant  $\sigma\neq 0$   is not scale invariant
 and  geometric characteristic of the corresponding  manifolds must include       a bound on  their "size"
 and/or  a  lower  bound on the mean curvatures of their boundaries (if there are any)including   faces of  of  spherical and hyperbolic  polyhedra $P$, if we want to prove their extremality. 
 


 
 
 \vspace {1mm}

 (c) Besides  $scal>0$, there are  other classes of Riemannian metrics that are 
 stable under geometric  connected sums of manifolds. The most prominent among these
  are  {\it conformally flat} metrics  (where one may simultaneously keep  positivity of $scal$
   if one wishes) and metrics with {\it  positive isotropic curvature}. 
   This curvature is defined in terms of the complexified 
   tangent bundle of $X$ and its positivity may be expressed in writing by $K_\mathbb C (X)>0$,
    \cite {micallef-moore}.   
 
 Probably, suitable limits of such connected sums can be 
 embraced by a general theory that would allow singular spaces.

{\it Sample Question.} Is there a natural class of singular spaces $X$ with  $K_\mathbb C (X)>0$
that  would satisfy  (a suitable version of) the Micallef-Moore \cite {micallef-moore} and/or  La Nave \cite {nave}    bounds on indices and sizes of harmonic maps of surfaces into $X$? \vspace {1mm}

 {\it More Questions. }
What is the geometry of    piecewise smooth Riemannian metrics? 

Manifolds with such metrics (which are not, in general, Alexandrov as they may
 have minus infinite sectional curvatures at their $(n-1)$-faces) seem to support full-fledged 
 Dirac operators and allow solutions of the Plateau problem.

 \vspace {1mm} 

{\it  Limit/Closure Question}. 
 What are natural limits of such spaces retaining their essential  "nice" properties. 
 
The spaces we expect  to have among such limits  should include the doubles of 
domains $P$ with boundaries in smooth manifolds, where the boundaries $\partial P$
may have certain  singularities. For example, these $\partial P$ may be (stable?)
 minimal hypersurfaces with singularities.

  \vspace {1mm} {\it Approximation Question.} When does
   a Riemannian metric $g$ on $X$  admit an approximation by piece-wise smooth ones, 
   where

   $\bullet$ the interiors of all
   pieces are flat (or, more generally have a given constant curvature),

   $\bullet'$   the exterior  curvatures of all
   $(n-1)$-faces (looked from the interiors of the corresponding $n$-faces) satisfy some convexity condition,

  $\bullet''$   the sums of the dihedral angles between $(n-1)$-faces around all
   $(n-2)$-faces are bounded by $2\pi$?
 
 For example, let $Sc(g)\geq 0$. Can one approximate it with   $\bullet'$ signifying 
 positive mean curvatures?  
 
 This, possibly, can be done (at least for $Sc>0$)  by starting with a fine triangulation of $X$  into 
 fat simplices
 with almost flat faces, and then flattening the insides of these $n$-simplices and, at the same time,
 gaining  the mean curvature convexity of their $(n-1)$-faces while keeping the  total angles around the $(n-2)$
   faces equal $2\pi$)
 
 (If  $  sect.curv \geq 0$ , then one may ask for 
 convexity instead of mean convexity but this seems less realistic.)

\section { Reflection Domains.}

An $n$-manifold $P$ with 
corners  is called { a \it $\Gamma$-reflection domain} if it is 
represented as  a {\it fundamental domain}  of  a { \it discrete reflection group} 
$\Gamma$ which acts on  a topological space 
 $\tilde P \supset P$ that is seen as   {\it an orbifold covering} or  {\it reflection  development} of  $\tilde P$.

This means that $P\in \tilde P$ is a  {\it domain}, i.e. its topological boundary in $\tilde P$
 equals the boundary $\partial P$, and $\Gamma$ is generated by 
  {\it reflections} $R_i$ in the $(n-1)$ faces $Q_i \subset P$:

  \vspace {1mm}
  
 $\bullet$  every $R_i: \tilde P \to\tilde P$ is an involution, $R_i^2=id$; 
  
 $\bullet$ $R_i(P)\cap P=Q_i$;
  
   $\bullet$ $R_i$  fixes
 $Q_i$.
 
  \vspace {1mm}
 
 {\it Regularity.}  A  $\Gamma$-reflection domain $P\subset \tilde P$  is  called {\it regular}, 
 if $\tilde P$ is a manifold
 which
 admits a smooth Riemannian $\Gamma$-invariant metric.

 \vspace {1mm}

 Let  $2k_{ij}$ denote the number of $\Gamma$-transformed domains 
  $\gamma(P)\subset \tilde P$, including $P=id(P)$, that 
 contain the
 $(n-2)$-face $Q_i\cap Q_j \subset P$ and define {\it $\Gamma$-angles} of 
 $P \subset \tilde P$  as 
 $$\angle_\Gamma [Q_i\cap Q_j]=\pi/k_{ij}.$$
  
 Notice, that our $P$ and $\tilde P$ are not endowed with any metrics so far and 
 the  $\Gamma$-angles are purely topological/combinatorial invariants.
 But if $\tilde P$ {\it is}  a {\it smooth} Riemannian manifold and $\Gamma$ acts 
 by {\it isometries}, then the $\Gamma$-angles equal the 
 dihedral angles $\angle[Q_i\pitchfork Q_j]$.

  \vspace {1mm}
  
  {\it Example: Rectangular Domains.}  Let $P$ be a  {\it co-simple}  polyhedral domain $P$, i.e. where the intersections of  all $k$-tuples of  $(n-1)$-faces  satisfy
    $$dim (Q_1\cap Q_2\cap ... \cap Q_k)\leq n-k+1.$$
Then $P$ has a natural  regular reflection structure with all dihedral angles $\pi/2$.    

If one glues $P$ with its  $ R_i$-reflected copy along the corresponding face $Q_i$,  then $P'=R_i(P) \cup P$
carries again a rectangular reflection structure. Thus, by consecutively   applying  such reflections with gluing, one constructs  a manifold 
$\tilde P$ with the corresponding reflection group $\Gamma$  generated by $R_i$ acting on  $ \tilde P$.
  
 For instance, if $P$ is the $n$-cube then $\tilde P\subset
  \mathbb R^n$ with the reflection group $\Gamma$ being a finite extension of $\mathbb Z^n$.
  
  \vspace {2mm}

\subsection {Gluing around  Edges and  Dihedral Extremality Theorem.}

 Let a compact connected Riemannian  mean curvature convex 
$n$-manifold $P$ with corners be represented by a regular
$\Gamma$-reflection domain in $\tilde P \supset P$  such that the (geometric)
dihedral angles of $P$ are bounded by the corresponding (topological)  $\Gamma$-angles,

$$ \angle [Q_i\pitchfork Q_j]\leq  \angle_\Gamma [Q_i\cap Q_j] \mbox {
 for all $(n-2)$-faces $Q_i\cap Q_j \subset P$}.$$

 {\it If $scal(P)\geq 0$, then the manifold $\tilde P$ admits a smooth Riemannian
   $\Gamma$-invariant metric $\tilde g_{reg}$ with positive
    scalar curvature, $scal(\tilde g_{reg})>0$, unless 
 
 $\bullet$  all dihedral angles $ \angle [Q_i\pitchfork Q_j]$
  are constant and are equal to $\angle_\Gamma [Q_i\cap Q_j]$; 
 
  $\bullet$ the mean curvatures of all $(n-1)$-faces $Q_i\subset P$ equal zero;
 
  $\bullet$ the scalar curvature of $P$ is everywhere zero. }
 
 \vspace {1mm} 
 
{\it  Proof.}  Every Riemannian metric $g_P$ on $P\subset \tilde P$  
obviously extends to 
  a unique  $\Gamma$-invariant path
 metric  $\tilde g$ on $\tilde P$ but, typically, this $\tilde g$ is singular on the boundary of 
 $P\subset \tilde P$.
 
 However, the three  inequalities: 
 $$  \angle_\Gamma [Q_i\cap Q_j] -\angle [Q_i\pitchfork Q_j]\geq 0, \mbox {  } 
 mn.curv (Q_i)\geq 0 \mbox { and } scal (g_P)\geq 0 \leqno {[3_{\geq 0}]}$$
   say, in effect, that $scal(\tilde g)\geq 0$ in 
 some generalized sense. 
 
In fact, if there is no $(n-2)$-faces at all  and $ \tilde P$ equals the double of $P$ along a 
mean convex boundary,  this $\tilde g$ can be easily approximated by a smooth metric 
$\tilde g_{reg}$ with
 $scal (\tilde g_{reg})>0$ as is explained in \cite {gromov-lawson1} and in \cite{almeida} for $scal(g_P)>0$ (also see  "Gluing with Positive Scalar Curvature" below)
 and where   non-vanishing of   
 $scal(g_P)$ at a single point  inside $P$ actually suffices because the positivity  of $scal$ can 
 be "redistributed"
 over all  of $P$  from a single point by a simple perturbation argument. 
 (It is easier to make such metrics with $scal>0$ not on $\tilde P$ itself but on $\tilde P\times \mathbb T^N$ by the warping argument from section 12 in \cite {gromov-lawson3}.)
 
This takes the care of rectangular domains, where an essential example is that of $P$ being cubical.
 
  In general,   a  close look at the "double-gluing/smoothing"  argument also  shows that it 
goes well along with
  $\angle_\Gamma [Q_i\cap Q_j] \geq \angle [Q_i\pitchfork Q_j]\geq 0$: 
  if $P$ is connected and {\it at least one}
   of the  three inequalities    $[3_{\geq 0}]$ is {\it non-strict} at some point, then
   the metric $\tilde g$ on $\tilde P\supset P$ 
    admits a smooth  $ \Gamma$-invariant approximation $\tilde g_{reg}$  with 
    $scal(\tilde g_{reg})>0$.
  
  This is obvious
 for $n=2$ and the general case is not difficult. 

 \vspace {1mm}
 
\textbf{  Corollary: Dihedral Extremality Theorem.}  {\it  Let $P\subset \mathbb R^n$ be a convex  polyhedron  where all dihedral angles are integer fraction of $\pi$, that are   $\pi/k$  for some integers $k$.}

(An essential instance of this is 
 the $n$-cube $P$  with the dihedral angles equal $\pi/2$  at all  $(n-2)$-faces and where a simple instance of a reflection domain with different angles is 
  the Cartesian  product of an $(n-2)$-cube with a   regular triangle where there are $\pi/2$ and $\pi/3$ angles.)

{\it Then $P$ is dihedrally extremal: there is no curve-faced  domain   $P'\subset \mathbb R^n$  combinatorially equivalent to $P$, such that the
 faces of   $P$ have strictly positive mean curvatures and all dihedral angles are bounded by the values of  the dihedral angles at the corresponding $(n-2)$-faces of $P$.}  \vspace{1mm}

{ \it Proof }  Since every Euclidean 
 reflection group $\Gamma$ contains $\mathbb Z^n$ of finite index, the 
 mean curvature convex  dihedrally extremality  of these $P$   follows from the Geroch conjecture for the 
torus $\mathbb R^n/\mathbb Z^n$.\vspace {1mm}
 
 {\it Generalizations with  $scal\geq 0$.} The above does not take  much of the Euclidean geometry of $P$, but  rather applies to  general corned Riemannian manifolds  $P$ with $scal(P)\geq 0$.  For instance,

 {\it no {\it essential} $n$-cubical $P$ with $scal(P)\geq 0$ can have acute dihedral angles and strictly mean curvature convex faces,  provided $P$ is  spin
 or  $n=dim(P)\leq 7$.} 
 
 (The non-spin cases for $n=8,9$ are settled   in section 5.3  and \cite {lohkamp}  allows all $n$.)

 \vspace {1mm}

 {\it  On Irregular Reflection Domains $P$. }   An orbifold covering or  "reflection  development" $\tilde P$ of a cornered manifold $P$ may have topological  singularities issuing from finite reflection (sub)groups acting  at the corners of $P$. These singularities  may be avoided if we replace 
 $P$ by its multi-double $P^{''....}$ without corners (see section 1.1).

 Since   $P^{''....}$ 
inherits from $P$ the (strict)  inequalities $scal>0$, $mean.curv>0$ and $\angle_{ij}>\alpha_{ij}$,    one can derive lower bounds on $\alpha_{ij}$ whenever the topology of 
 $\tilde P^{''...}$, allows no metric of positive scalar curvature  invariant under  the  reflection group acting on $\tilde P^{''...}$.

Yet, this does not help unless $\alpha_{ij}=\pi/k$.\vspace{1mm}

 {\it Gluing with Positive Scalar Curvature.}  Let  $X_1$ and $X_2$ be smooth Riemannian manifolds 
 and let $\gamma: \partial X_1 \leftrightarrow \partial X_1$ be an isometry between their boundaries. 
 
The manifold $X^{1\cup 2} =X_1\cup_\gamma X_2$  obtained by gluing the two along the boundaries
carries a natural {\it continuous} Riemannian metric $g^{1\cup 2}$. 

 If the  {\it shape operators} (corresponding to the second fundamental forms)  $A^\ast_1$ and $A^\ast_2$ of the two boundaries "match" , i.e.   $A^\ast_1=_\gamma A^\ast_2$ for $\partial X_1$ cooriented outward and $\partial X_2$ inward,   then $g^{1\cup 2}$ is $C^1$-smooth. It follows, that

 $(\ast)$    {\it if the scalar curvatures of $X_1$ and $X_2$ are strictly positive, then $g^{1\cup 2}$ can be smoothed to a metric that also has
$scal>0$}.

Indeed,  $g\mapsto scal(g)$  is a differential  operator on Riemannian metrics $g$ on $X^{1\cup 2}$ that is  {\it linear} in the second derivatives of $g$ and so the  smoothing with any standard smoothing kernel does the job.

The gluing construction from \cite {gromov-lawson1}  delivers, in effect, a deformation of a metric  $g$ on a manifold $X$ with $scal(X)>0$ and $mn.curv(\partial X)>0$  that  keeps  $scal.curv >0$, that does not change the restriction of $g$ to  $\partial X$  and that makes  the second fundamental form zero. Then the above gives one a metric with $scal>0$ on the double of $X$.\vspace {1mm}

 More generally,  let $Y_1$ and $Y_2$ be, say closed, Riemannian manifolds with metrics $g_1$ and $g_2$ that are also are "decorated" by quadratic differential forms $A_1$ and $A_2$ and let us look at 
   compact (complete?) smooth Riemannian  manifolds $X=X_{12} =(X,g)$ such that\vspace {1mm}

$\bullet$ the boundary of such an $X$ is decomposed  into a disjoint union,  $\partial X=\partial_1X\sqcup \partial_2(X)$, where  
$\partial X_1$ is cooriented by an {\it  inward} vector field and $ \partial_2(X)$, by an  {\it outward} field; 

$\bullet$ there are  
isometries $I_i:Y_{i}\to \partial_iX$, $i=1,2$, that induce the metrics $g_i$ from $g$ and  send
$A_i$ to the second quadratic (exterior curvature) forms  of $\partial_iX\subset X$where these forms are evaluated with given coorientations.  \vspace {1mm}

For instance, $Y_1$ and $Y_2$ may be  be two concentric spheres in $\mathbb R^n$ with outward coorientations. If $Y_1$ is contained in the ball bounded by $Y_2$, then these spheres    serve as the boundary of the annulus $X$ between them,
where, according to our   convention, both (second quadratic)  forms,  $A_1$ and $A_2$ --   on the concave  interior  $Y_1$-boundary with the inward coorientation  and on the convex exterior $Y_2$-one cooriented by an outward field --  are positive definite.\vspace {1mm}

These $X=X_{12}$ may be seen (almost) as morphisms between "decorated" Riemannian manifolds $Y$: if  we  glue $X_{23}$ to $X_{12}$ along $Y_2$ the resulting metric is 
$C^1$-smooth and the true $X_{13}$ is obtained by smoothing this metric.(One could avoid smoothing if working with $C^{1,1}$-metrics.) 

The metrics and quadratic form serve for defining several more interesting  smaller categories such as

 (1) ${\cal B}_{sc>0}^n$, $n=dim(X)$, the subcategory of the above category  where  the manifolds $X$  have {\it strictly positive} 
scalar curvatures, $sc(X)>0$,  (the category  ${\cal B}_{sc\geq0}^n$ is equally interesting but slightly 
harder to handle);

 (2) the subcategory  ${\cal B}_{sc>0}^n$ where manifolds $X$ are cobordisms;

(3) the subcategory made by   those $X$ where the distance function $x\mapsto dist_X(x,\partial_1X)$ is smooth
with $\partial_2X$ being a constant level of this function. 
 
 (These  three categories naturally extend to $\infty$-categories with $d$-morphisms 
being represented by $d$-cubical cornered manifolds, that suggests a  "topological field theory" for $scal>0$. Also it is amusing to think of reflection groups as of "enhanced" $\infty$-categories.)
 
  \vspace {1mm}
  
  Let us focus our at attention on an, "infinitesimal $\varepsilon$-subcategory" of (3) where
the distance between $\partial_1X$ and  $\partial_2X$ equals $\varepsilon\to 0$ and where
the $C^2$-distance between the metrics $g_1=g_{|\partial_1X}$ and  $g_2=g_{|\partial_2X}$
is also $\leq \varepsilon$ when the two metrics brought to the same manifold, say to ${\partial_1X}$ via the normal projections ${\partial_1X}\leftrightarrow {\partial_2X}$, and
answer the following question.\vspace {1mm}

{\it When can a  quadratic form $A_2$ on $Y_2=\partial X$  be obtained by 
an equidistant $\varepsilon$-deformation of $Y_1=\partial_1X$ with a given quadratic form $A_1$ on $Y_1$.}.\vspace {1mm}

To  effortlessly compute the curvatures, etc. of equidistant deformations  $Y_t$ of hypersurfaces $Y\subset X$, 
recall  that, in general, \vspace {1mm}

 {\it the first derivative of a Riemannian metric $g$ on $X$ restricted to  $Y_t\subset X$ under the normal equidistant deformation $Y_t$  equals the  second fundamental form 
$A_t$ of $Y_t$},  where both forms  $g_t=g_{|Y_t}$ and $A_t$ are brought to $Y$ by the normal projection $Y\leftrightarrow Y_t$, \vspace {1mm}
 $$\frac {d}{dt} g_t=A_t,  \leqno {[\frac{d}{dt}]} $$ 
 while the derivative of the corresponding shape operator $A_t^\ast$, that is  defined by $g_t(A_t^\ast(\tau),\tau)=A_t(\tau, \tau)$, for the tangent vectors $\tau\in T(Y_t)$,
and then brought to $Y$ by $Y\leftrightarrow Y_t$, is expressed by the \vspace {1mm}
 
  {\it   Riemannian   Hermann Weyl Tube Formula}\footnote{This is a most useful formula in Riemannian geometry that {\it directly} leads to geometrically significant results, e.g. to basic {\it comparison theorems} (see \cite {sign}), unlike   still persistent      roundabout computations with curvature tensors and  Jacobi fields used for this purpose.}:
$$\frac {d}{dt} A_t^\ast=-(A_t^\ast)^2 +B_t \leqno {[\frac{d^2}{dt^2}]}$$ 
where the operators $B_t$ are defined via the sectional curvatures $K$ of $(X,g)$ on the $2$-planes 
$\sigma\subset T_y(X)$, $y\in Y_t$, that are normal to the tangent spaces  $T_y(Y_t)$, as follows,
 $$g(B_t(\tau), \tau)= -K(\sigma),\leqno {[K]}$$
where $\tau$ is a unit vector in the line $ \sigma\cap T_y(Y_t)$. (We use here the notation from \cite {sign}, p 43.)

Now, let $U_\varepsilon=Y \times [0,\varepsilon]$ for $Y=Y\times 0$,
 and let $g_\varepsilon$ be the following metric on $U_\varepsilon$ defined with
 given  smooth metric $g_0$ and two smooth quadratic forms $A_0$  and $A_+$ on $Y$,
   $$g_\varepsilon(y,t)=g_0(y) + t A_0(y) + \frac {t^2}{2\varepsilon}( A_{+}(y)-A_{0}(y)), \mbox { } 0\leq t\leq \varepsilon.\leqno {(++)}$$

When $\varepsilon \to 0$, then

$\bullet_g $  the metrics $(g_\varepsilon)_{|Y\times t}$ converge to $g_0$ in the $C^2$-topology;

$\bullet_A $ the second  quadratic forms of $Y=Y\times \varepsilon \subset U_\varepsilon$ similarly
$C^2$-converge to $A_+$;

$\bullet_{scale} $  if $trace (A_+(y)) < trace (A_0(y))$ at all $y\in Y$, then  $scal(U_\varepsilon)\to +\infty$ at all $x\in U_\varepsilon.$   \vspace {1mm}

In other words, an infinitesimal positive scalar  curvature  cobordism/morphism  can transform  a quadratic form $A_0$ to a given $A_+$, whenever the mean curvature (trace) of the latter is strictly smaller than that of the former.

 (It is helpful to visualize this by thinking of   
 $ Y_\varepsilon=Y\times \varepsilon$ as  the equidistant $\varepsilon$-deformation  of an equatorial $(n-1)$-subsphere $Y_0$ in the $n$-sphere of radius $r$, with 
 $U_\varepsilon$ being the annulus pinched between these two spheres, where $scal(U_\varepsilon)= n(n-1)/r^2$ and  the second fundamental form $A_0$ of $Y_0$ is zero. 
 
 A relevant example here is where
 $r\to 0$ and $\varepsilon =const\cdot  r^2$ and  where the second fundamental form $A_+$ of $Y_\varepsilon$ is, negative,  being, roughly,    $-const\cdot \mathbf 1$.)

  \vspace {1mm}

{\it Proof.} The claim   $\bullet_g $ is obvious,  $\bullet_A $ follows from  $[\frac{d}{dt}]$, while 
  the asymptotic  estimate $scal ( U_\varepsilon) \sim  const\cdot \varepsilon^{-1}(trace (A_0) - trace (A_+))$ is seen with $[\frac{d^2}{dt^2}]$, $[K]$    and  {\it Gauss' theorema egregium.} 

 \vspace {1mm}

{\it Gluing Corollary.} Let $X_1$ and $X_2$ be manifolds with strictly positive scalar curvatures and let $\gamma:\partial X_1\leftrightarrow \partial X_2$ be an isometry.

{\it If the mean curvatures of the two boundaries satisfy
 $$mn.curv_{\partial X_1}(y_1)+ mn.curv_{\partial X_2}(y_2)\geq 0 \mbox { for all
 $y_1 \in \partial X_1$  and  $y_2=\gamma(y_1)\in \partial X_2$,}$$  
 then 
the metric  $g^{1\cup 2}$ on the manifold  $X^{1\cup 2}$ obtained by gluing the two along their boundaries  can be perturbed in a neighbourhood  of the glued  boundaries to a metric of positive scalar curvature 
 on  $X^{1\cup 2}$.} 
 \vspace {1mm}  
 
{\it Proof.}  If we apply $(++)$ to $g'_0(y)=g_0(y)+\delta(y)$ instead of $g_0(y)$  for 
$\delta(y)= g_0(y)-g_\varepsilon(y,\varepsilon)$  we end up with $g'_\varepsilon(y,\varepsilon)=g_0(y)$. Thus, an arbitrary small $C^2$-perturbation of $g_0$ allows us to achieve $g_\varepsilon(y,\varepsilon)=g_0$.

Now we can modify  the metric on one of the two manifolds, say on $X_1$,  such that

$\bullet$ the modified metric equals the original one away from an arbitrary small neighbourhood (that is our  $U_\varepsilon \supset Y=\partial X_1$) of the  boundary of $X_1$;

$\bullet$ the restriction of the modified metric to the boundary remains equal the original metric (corresponding to the above $g_0$);
 
$\bullet$ the modified metric has   $scal>0$;

$\bullet$ the second fundamental form of the boundary with respect to the  modified metric has the second fundamental form  opposite to that of the boundary of $X_2.$
 
Then, by the above $(\ast)$, the manifold $X_1$ with the modified metric can be glued to $X_2$  and the proof follows.

 \vspace {1mm}

{\it  Gluing with  $scal>\kappa\neq 0$.} The above equally applies to manifolds with scalar curvatures bounded from below by any constant $\kappa$, not necessarily $\kappa=0$:

 \vspace {1mm}
 
 {\it If  manifolds $X_1$ and $X_2$  with their scalar curvatures bounded from below by a constant $\kappa$  are glued by an isometry  $\gamma$ between their boundaries the mean curvatures of which  satisfy the above positivity } ($\geq 0$) {\it of their sum  inequality,  then the  metric $g^{1\cup 2}$ on the resulting manifold   $X^{1\cup 2}=X_1\cup_\gamma X_2$ can be  approximated by smooth metrics  $g_{app}$ on  $X^{1\cup 2}$  with $scal(g_{app})>\kappa$  where these   metrics $g_{app}$ can be chosen
 equal  $g^{1\cup 2}$  away from an arbitrarily small neighborhood of the glued boundaries $\partial X_1\underset{\gamma}= \partial X_2\subset X^{1\cup 2}$.}
 \vspace {1mm}

 \vspace {1mm}
{\it Remarks.}  (a) The above style local  "gluing+ smoothing"  appears in different forms in 
\cite{gromov-lawson1}, \cite {almeida}, \cite {miao} and non-local  smoothing with the  Ricci flow is suggested in \cite {mcferon}.\vspace {1mm}

(b) Probably, ideas from \cite{miao},  \cite{mcferon}, \cite{lee} and  \cite{grant} may help to establish a  version of this   under  the  {\it non-strict} assumption $scal  (X_i)\geq \kappa$, $i=1,2$, and with the corresponding  non-strict conclusion   $scal (g_{app})\geq \kappa$  for most (all) manifolds $X_1$ and $X_2$. \vspace {1mm}

(c) If
two given smooth Riemannian metrics on $Y_1$ and on $Y_2=Y_1$ can be  included in a {\it  continuous family }  of metrics with {\it positive scalar curvatures},  then, by combining the above with \cite {gromov-lawson2}, one sees that   the existence of a  {\it non-infinitesimal} cobordism as in the above (3) with given  second quadratic forms  $A_1$ on $Y_1$ and $A_2$ on $Y_1$ and $Y_2$  {\it does not need} the  assumption  $trace (A_1)>trace (A_2)$.
 Let us spell it out in detail. \vspace {1mm}

Let $A_1$ and $A_2$ be smooth quadratic forms on a manifold $Y$  an let $g$ be a metric on $X=Y\times [1,2]$, such that

$(+)$ the metrics  $g_{|Y\times t}$ on $Y\times t=Y$ have  {\it strictly positive scalar curvatures} for all $t\in [1,2]$. \vspace {1mm}

Then  there exists a homotopy $g_\tau$ of the metric  $g=g_{\tau=0}$ on $X$, such that \vspace {1mm}

$\bullet $ the metrics  $g_{\tau{|Y\times t}}$ on $Y=Y\times t$ have positive scalar curvatures for all $t\in [1,2]$ and $\tau\in[0,1]$;

$\bullet $ the result of this homotopy --  the metric $g_{\tau=1}$ on $X$  has strictly positive scalar curvature;

$\bullet $ the homotopy is constant on $Y\times 1$ and $Y\times 2$,
$$g_{\tau|Y\times 1}=g_{|Y\times 1}\mbox {  and } g_{\tau|Y\times 2}=g_{|Y\times 2}\mbox { for all $\tau \in [0,1]$};$$

$\bullet $  the second quadratic forms of $Y\times 1$ and $Y\times 2$ in $(X, g_{\tau=1})$ equal $A_1$ and $A_2$ correspondingly.

$\bullet $ the submanifolds $Y\times t\subset X$ are equidistant to $Y\times 1$ as well as to $Y\times 2$ for all $t\in [1,2]$ (as in the above (3))
with respect to $g_{\tau=1}$; moreover, if $g$ has this equidistance property, then one
can have  all $g_\tau$ with this property as well.
\vspace {1mm}

{\it Question.} When does a closed subset  $Z$ in a Riemannian manifold $X$ equal the intersection of an decreasing sequence of  smooth domains $U_i\subset X$ where the induced metric $g_i$ on the boundaries $Y_i=\partial U_i$ have $scal(g_i)\to +\infty$ for $i\to \infty$?

Is there  a sufficient condition representable by an inequality $dim_{?}(Z)< dim(X)-2$  for some  notion of dimension
as it is the case for 
   {\it piecewise smooth} polyhedral  subsets $Z\subset X$ of codimension $>2$  by the argument from \cite {gromov-lawson2}.

A related question (we reiterate it in section 3) is that of finding "nice" functions $\phi_i$ on $X\setminus Z$ that blow up at $Z$ and such that
the intersection of  certain $\phi_i$-bubbles equals $Z$.

\vspace {1mm}

(d) There are   global PDE  constructions  of metric with positive  scalar curvatures on
"glued manifolds" like  the above  $X^{1\cup 2}$ under
 integral rather than  point-wise  assumptions on the mean curvatures (e.g. see \cite {wang} and references therein) but the available results of this kind  apply  so far only  to rather special metrics.

 \subsection {Dihedral Rigidity Conjecture.}  The above does not say {\it what  are} dihedrally  extremal mean curvature domains and, more generally,  what  are the above cornered  Riemannian  $n$-manifold $P$
where
  $$  \angle [Q_i\pitchfork Q_j]  = \angle_\Gamma [Q_i\cap Q_j], \mbox {  } 
mn.curv (Q_i)= 0 \mbox { and } scal (g_P)= 0.$$

 Probably, they are all isometric to convex  Euclidean polyhedra.  In particular, 

 { \it bounded Euclidean polyhedral reflection domains  $P\subset \mathbb R^n$  are, conjecturally,  dihedrally  rigid.}\vspace {1mm} 

Namely, \vspace {1mm}

 let a  curve-faced  $P'\subset \mathbb R^n $ have not necessarily strictly positive  mean curvatures of all their faces $\geq 0$ and all dihedral angles bounded by the corresponding angles of $P$.

{\it  Then,
conjecturally, all faces of $P'$ are flat; moreover, $P'$ is obtained from $P$ by parallel translations of its faces followed by an isometry.}


 
\vspace {1mm}


\vspace {1mm}

{\it  Five Incomplete    Proofs.} 
 (1) The most transparent case of the problem is where $P $ is a curve-faced  cubical polyhedron
 in the Euclidean $3$-space $\mathbb R^3$.
 
  If $P$ is extremal,
 then all $2$-faces of it are minimal surfaces meeting each other at the angles $\pi/2$.
 If such a face $Q$ can be slightly moved inside $P$ with a strict decrease of its area, then it can be perturbed to strict mean convexity and if every such a move   strictly increase the area of $Q$, then one could make  $Q$ strictly mean curvature concave by such a move. 
 
 Thus, we may assume that 
 
 {\it $Q$
 includes in a continuous family $Q_t \subset P$ of minimal faces normal to the rest of the boundary of $P$.} 
Now,  following  Schoen-Yau  \cite {schoen-yau1}, (compare    \cite {burago-toponogov})   we observe that 
  the second variation integral for  $area(Q)$ equals  the integral of the Gauss curvature of $Q$  plus  
 the boundary term that is the integral over the curve $\partial Q$ of the difference between the mean curvatures of $\partial Q$ in $Q$ and of the surface of $P$ normal to $Q$. 

Thus,  the (non-strict) positivity of the second variation implies that the integral of the Gauss curvature of $Q$ plus the integral of the curvatures of the four edges of its boundary is non-negative. On the other 
hand, the four vertices contribute $2\pi$ to the Gauss-Bonnet integral; hence, $Q$ must be flat. 

The only unsettled  point in this proof is hating that nothing bad happens at the corners of $P$  but this does not seem to be difficult.
 In fact, the Gauss-Bonnet prism inequality (see 5.4)
does work in this case and implies rigidity of all, not necessarily (reflection) convex Euclidean prisms.
\vspace {1mm}

If  $P$ is  a curve-faced  cubical polyhedron 
 in  $\mathbb R^n$ for $n\geq 4$, one may apply the   Schoen-Yau dimension reduction argument \cite {schoen-yau2}
or rather the  warped product  version of it from \cite {gromov-lawson3}, where, one needs doing it only once, and where the dimension restriction $n\leq 6$ is unnecessary since all we  need after all is a certain {\it perturbation} of  a {\it smooth} face $Q$. \vspace {1mm}

 (2)   Let us turn now to general reflection domains and recall that if $X$ is a compact manifold with zero scalar curvature $g_0$,
then, according to \cite {kazdan-warner}, $g_0$ can be perturbed to $g_1$ with $scal (g_1)>0$ unless $g$ is Ricci flat.  Probably, a similar perturbation is possible for  the above cornered  Riemannian  $n$-manifold $P$
with
  $$  \angle [Q_i\pitchfork Q_j]  = \angle_\Gamma [Q_i\cap Q_j], \mbox {  } 
mn.curv (Q_i)= 0 \mbox { and } scal (g_P)= 0.$$

This suggest the following \vspace {1mm}

{\it Question.} Let $P\subset \mathbb R^n$ be a preconvex  Plateau-hedron where at least one face is {\it non-flat}. When  does $P$  admit a  perturbation  to a  {\it strictly} mean curvature convex
polyhedron $P'$ with all dihedral angles $\angle_{ij}(P')\leq \angle_{ij}(P)$?
(A general  positive answer would settle the rigidity problem.)

{\it Example.} Let $P\subset \mathbb R^3$ be bounded by a catenoid and a pair of hyperplanes normal to its axes. (This $P$ has depth two as it has no corners.)
It seems easy to decide if    $P$ can be perturbed to a  strictly mean convex polyhedron 
with a decrease of its two (circular) dihedral angles.
\vspace {1mm}

 (3)  Let   $\tilde P$  be the above  manifold where  reflection group $\Gamma$ acts with $ P\subset \tilde P$ being  a $\Gamma$-reflection domain and  let $\tilde g$ on $\tilde P\supset P$ be the continuous metric coming from $P$. 
 
 The above gluing argument says in this case that $\tilde g$ 
    admits a smooth  $ \Gamma$-invariant approximation  by metrics $\tilde g_{reg}$  with 
    $scal(\tilde g_{reg})>-\varepsilon $ for arbitrarily small $\varepsilon>0$. Then the positive solution  to  the  $C^0$-rigidity problem  $\bullet_2$ in 1.9 would  show that the metric  $\tilde g$  is flat. \vspace {1mm}

(4)
\textbf { Torus Conjecture.} Let $\tilde X$ be the universal covering of the $n$-torus $X$  with a smooth non-flat Riemannian metric. 

{\it Then, conjecturally, $\tilde X$ can be exhausted by cubical strictly  mean curvature convex  cornered domains $\tilde P\subset \tilde X$  with all dihedral angles bounded by $\pi/2$.}\vspace {1mm}

 {\it Half-Proof}.  Such domains   $\tilde P_k$ with {\it singular} faces are constructed    as follows. 
 
Let $ X^{{\sim}_ k} =X/k\mathbb Z^n$  for $Z^n=\pi_1(X)$ and let $Y_1 \subset X^{{\sim}_ k} $ be  a minimal hypersurface homologous to an $(n-1)$-subtorus   in $ X^{\sim_ k}$. Let $ X_1^{{\sim}_ k}$ be (possibly singular for $n\geq 8$) space obtained 
by cutting   $ X^{{\sim}_ k}$ along $Y_1$ where the  boundary of $ X_1^{{\sim}_ k}$,  denoted  $\partial X^{{\sim}_ k}=Y_{\pm 1}$,  consists of two copies of $Y_1$. 

Take an $(n-1)$-volume minimizing  hypersurface $Y_2\subset X_1^{{\sim}_ k}$  with $\partial Y_2\subset \partial  X_1^{{\sim}_ k} = Y_{\pm 1}$ that represent the relative homology class "suspending" the class of some "$(n-2)$-subtorus"   in $Y_1 $. 

  Cut $ X_1^{{\sim}_ k}$ along  $Y_2$  and take the result in space for  $ X_2^{{\sim}_ k}$. This  $ X_2^{{\sim}_ k}$ is a cornered $2$-cubical space where the boundary $\partial X_2^{{\sim}_ k}$ consists of    two pairs mutually orthogonal faces  that are $Y^{cut}_{\pm 1}$  and 
$Y_{\pm 2}$.

 Cut $X_2^{{\sim}_ k}$ along a  a volume minimizing hypersurface  $Y_3$ in  $\partial X_2^{{\sim}_ k}$ with boundary in $\partial X_2^{{\sim}_ k}$ and thus obtain a $3$-cubical $X_2^{{\sim}_ k}$ with norma n dihedral angles where these are defined.
 
 Keep doing this unless you arrive a singular mean curvature  $n$-cubical  cornered space $X_n^{{\sim}_ k}$ with normal faces that lifts to
 a  singular cubical domain   $ \tilde P_k\tilde X$, where $\tilde X$   can be exhausted by such  $\tilde P_k$ for  $k\to \infty$.

What remains is to  smooth the  singularities in these $\tilde P_k$. Prior to smoothing
one has to  modify the construction by taking {\it $\varepsilon$-bubbles} with small $\varepsilon>0$ (defined below) instead of  {\it minimal }  hypersurfaces $Y_i\subset X_{i-1}^{{\sim}_ k}$  (corresponding to $\varepsilon=0$). 

Such bubbles do exist if $X$ is  {\it non-flat}  (this seems obvious  but a proof will not hurt) and then the resulting polyhedra $\tilde P_k\subset \tilde X$
are {\it strictly} mean curvature positive. Probably -- and this is what we {\it do not prove} -- such polyhedral domains $\tilde P$ can be always approximated by  face-wise smooth strictly mean curvature convex domains with $\pi/2$ dihedral angles .\vspace {1mm}

A construction of such  strictly  mean curvature convex cubical  domains $Q$  in  manifolds  $\tilde X$ with a {\it  piecewise} smooth metrics  $\tilde g$ would imply   dihedral rigidity of reflection domains, since one could construct a metric of positive  scalar curvature on $X$ (or rather on some $X'$ admitting a  degree $1$ map to the $n$-torus) by "gluing  $Q$  around the corners".\vspace {1mm}

(5) Probably, one can make sense of the Dirac $D_{\tilde g}$ operator being (non-strictly) positive on $(\tilde P)$, confront this with existence of harmonic spinors twisted with flat bundles as in   \cite  {gromov-lawson1},
and use a piecewise smooth  version of Bourguignon theorem on parallel spinors.

\vspace {1mm}

\vspace {1mm}
{\it  Topological Mean Convex   Exhaustion Problem.} Let $\cal P$ be a combinatorial  class of (potential) cornered domains $P$ with  numbers $\alpha_i$ attached  to the edges of $P\in \cal P$. 

Let $\tilde X$ be a smooth manifold acted upon by a discrete cocompact group $\Gamma$. For instance,  $\tilde X$  may be  a universal covering of a closed manifold $X$ with $\pi_1(X)=\Gamma$.

 If $\tilde X$ is endowed with a smooth $\Gamma$-invariant RIemannian metrics $\tilde g$, we say that $\tilde X=(\tilde X, \tilde g):\Gamma $ is  {\it exhaustively  dominated} by $({\cal P},\alpha_i)$ if 
 $\tilde X $ can be exhausted by mean curvature convex domains $\tilde P\in \cal P$  where all dihedral angles in  all $P$ are bounded by $\alpha_i$.  (The  definition of combinatorial equivalence allows  maps $\tilde P \to P\in {\cal P}$ of "positive degrees" where, possibly, $dim(\tilde P)>dim(P)$.)

Say that  $\tilde X:\Gamma$ is {\it topologically exhaustively  dominated} by $({\cal P},\alpha_i)$ if $(\tilde X, \tilde g) $ is  exhaustively  dominated by $({\cal P},\alpha_i)$ for all  smooth $\Gamma$-invariant RIemannian metrics $\tilde g$ on $\tilde X$.

For instance, the torus conjecture claims that the topological $n$-torus  $T^n\mathbb= R^n/\mathbb Z^n$ is so dominated by  the  (combinatorial class of the)  cube  $(\square^n, \alpha_i=\pi/2)$.

In general, given  $\tilde X$ acted upon by $\Gamma$ and  a combinatorial class  $\cal P$  with the edge set $I=\{i\}$, we denote by $A_{\cal P}(\tilde X:\Gamma)\subset \mathbb R^I$ the set of vectors $\{\alpha_i\}\in \mathbb R^I$ for which  $({\cal P},\alpha_i)$
topologically exhaustively  dominate $\tilde X:\Gamma$.

This  set  $A_{\cal P}$ is a  {\it topological invariant} of $(\tilde X:\Gamma)$ (or equivalently of the quotient manifold $\tilde X/\Gamma$ for free actions) and the problem is to evaluate it in particular cases.\vspace {1mm}

 {\it Examples}.  (\rmnum{1})   Let  $P\subset \mathbb R^m$ be a convex polyhedron with dihedral angles $\alpha_i$.

Let $X$ be a closed smooth $n$-manifold that admits a continuous map   $X\to \mathbb T^m$ of "non-zero degree", i.e. such that the fundamental cohomology class of the torus goes to a {\it non-zero class} in $H^m(X,\mathbb Q)$ and let $\tilde X$
the Galois $\mathbb Z^m$-covering of $X$ induced by the universal covering $\mathbb R^m\to \mathbb T^m$.

Is   $\tilde X:\Gamma$ topologically exhaustively  dominated by $({\cal P},\alpha_i)$
for the combinatorial class  ${\cal P}$ of this $P\subset \mathbb R^m$?

  (\rmnum{2}) Does  this remain so, for $\Gamma\neq  \mathbb Z^m$, if  $\tilde X$ admits a smooth proper {\it distance decreasing}   map  $\tilde X \to \mathbb R^n$  of {\it non-zero degree}?
  
  Here  "non-zero degree" means that the pullback of a generic point  is $\mathbb Q$-non-homologous to zero in $\tilde X$ and "distance decreasing" is understood relative to some $\Gamma$-invariant metric in $\tilde X$.

  (\rmnum{3}) Are the universal coverings $\tilde X$  of manifolds $X$ with  {\it infinite $K$-area} (defined in  \cite {positive})    topologically exhaustively  dominated by these $({\cal P},\alpha_i)$?

\vspace {1mm}

  (\rmnum{4})   Let     the hyperbolic $n$-space $H^n$ of constant curvature $-1$ be exhausted by  curve-faced polyhedral domains $P$
  of a certain combinatorial type $\cal P$ with 
  umbilical faces  with positive mean curvatures    and acute dihedral angles.
  
  Is then $H^n:\Gamma$    {\it topologically} exhaustively  dominated by $({\cal P}, \alpha_i=\pi/2)$?
  
Does a similar property hold true for exhaustions of $H^n$ by  domains  $P_1\subset P_2\subset .... \subset H^n$  
of variable combinatorial types?
  
 What are combinatorial classes $\cal P$, such that $({\cal P}, \alpha_i=\pi/2)$ topologically exhaustively dominate (the universal coverings of)  {\it $\mathbb Q$-essential } closed manifolds 
$X$ with other  "large" fundamental  groups $\Gamma$?

(Recall that "$\mathbb Q$-essential" means that the fundamental cohomology class  $[X]_{\mathbb Q}\in H^n(X;\mathbb Q)$, $n=dim X$, comes from the cohomology of $\Gamma$.)

Particular instances of interesting "large" groups are cocompact lattices in Lie Groups and  Cartesian projects of
{\it world hyperbolic}  groups.

 Probably, there  are significantly of more such $\cal P$   for hyperbolic groups than
for products of these and than   for groups  co-compactly acting on non-hyperbolic  symmetric
 spaces. 
  
 \vspace {2mm}

Let is limit the $\Gamma$ invariant metrics $\tilde g$  on $\tilde X$ to  those where
$scal(\tilde G) \geq -1$.  Then, besides  bounds on  the dihedral angles of domains $P_k$ exhausting $\tilde X$  
  by given numbers $\alpha_i$  one may requite {\it lower bounds on  the  mean curvatures} of the faces of these  domains.
 
 What are realizable (by some exhaustions)   possibilities for such bounds if, for instance,

$(\bullet$) $\tilde X$ is   symmetric space with non-positive curvature;

$(\bullet$)  $\tilde X$  is acted upon by an isometry group $\Gamma$ (for a metric $\tilde g$ with $scal(\tilde g) \geq -1$) and it  admits a proper  equivariant map of non-zero degree onto a symmetric space that is  isometrically  and co-compactly acted upon by $\Gamma$. \vspace {2mm}

$(\ast$)  {\it What are extremal/rigid corned domains in spaces of constant (positive or negative) curvatures  with given bounds on  the dihedral angles and lower bounds on the mean curvatures of their faces?}\vspace {1mm}

 For instance, let  $P\subset S^n\subset \mathbb R^{n+1}$ be a reflection domain,
 e.g. the spherical  simplex  $\Delta$  with all dihedral angles equal $\pi/2$
 and let  $P$ be a  mean convex cornered manifold  with all dihedral angles $\leq \pi/2$ that admits a {\it $1$-Lipschitz} combinatorial equivalence  $f:P\to\Delta$.
 
  Then,
 if $P$ is spin,  the above argument combined with  Llarull's theorem \cite{llarull} shows that

  \vspace{1mm}
  
  {\it there is a point $p\in P$, such that $scal_p(P)\leq scal_{f(p)}(S^n)$.}  \vspace{1mm}
 
 But it is unclear if the equality  $scal_p(P)= scal_{f(p)}(S^n)$ at all $p\in P$
 implies that  $f$ is an isometry.
 
 \subsection{  Billiards, Pure Edges  and Ramified Coverings.} Most (if not all) of our understanding 
 of mean curvature convex  cornered manifolds $P$ is derived from the geometry of minimal varieties and/or the Dirac operator on orbifold coverings (reflection developments) $\tilde P$ of $P$.

 For example if $P\subset X=X^n$ is a $n$-cubical polyhedral domain, then the quotient  $P_\circ=\tilde P/\Gamma_\circ$,
 for  a subgroup $\Gamma_\circ = \mathbb Z^n \subset \Gamma$ of finite index in the corresponding reflection group, $\Gamma$ acting on $\tilde P$, is a
  $\mathbb T^n$-essential manifold, i.e. it comes with a map of positive degree $P_\circ\to \mathbb T^n$.

 The relevant minimal subvarieties in   $P_\circ$ are those   representing the $(n-1)$-homology classes that come as pullbacks  from $(n-1)$ subtori in $\mathbb T^n$.
 
There are infinitely many of these classes, but only $n$ of them, the ones that correspond to the coordinate subtori,  have a simple representation
 in $P$, namely by (eventually minimal) hypersurfaces separating pairs of opposite faces. 

The remaining ones are similar to  multiply reflected periodic orbits of billiards in polygonal  domains.

{\it Question.} Is there a counterpart of these multiply  reflective minimal  hypersurfaces for  polyhedral domains $P$ that are not reflection domains? \vspace {1mm}
   
\vspace {1mm}

A similar issue arises for the Dirac operator. \vspace {1mm}

How can one descent Dirac operator $D$  proofs of non-existence of metrics with $scl>0$  from 
$\tilde P$ to $P$?

One problem is  the discontinuity of $D$ for the natural (only continuous but not smooth) extension of the Riemannian metric from $P$ to $\tilde P\supset P$. 

Even more serious difficulty stems from the fact that it is not $D$ itself is used but $D$ twisted with (almost flat) vector bundles $V$ over $\tilde P$,  where these $V$ is bu no means $\Gamma$-invariant. \vspace{1mm}

An essential difference between ordinary billiards and what we have here is that the dynamics and geometry of  billiards are shaped by interactions of the orbits with  the {\it faces} of $P$ while the geometries of minimal hypersurfaces and of Dirac 
crucially depend on what happens at the {\it edges} of $P$. 

Below is an attempt to isolate  the edge geometry. 

\vspace {1mm}

{\it Pure Edges without Faces.}  Let $X=X^n$ be a closed manifold and $Z=Z^{n-2}\subset X$ a closed submanifold
of codimension $2$, e.g. a knot in the $3$-sphere.

Consider all metrics on $X$ that are smooth with non-negative scalar curvatures away from 
 $Z$ and such that the geometry near $Z$ is {\it corner singular with angle $\alpha$},  i.e. a neighbourhood of $Z$ in $X$ is isometric to   $(Z, g_Z)\times C_\alpha$ where  $g_Z$ is a smooth Riemannian metric on $Z$ and  $C_\alpha$ is a surface with a rotationally symmetric  Riemannian metric that is  singular at a single point $c_0\in C_\alpha$ where its tangent cone has total angle $\alpha$. For example, if $\alpha=2\pi$
 then $C_\alpha$ is non-singular.

Denote by $\alpha_{min} (X,Z)$ the infimum of these angles of all above metrics.

This is  a {\it topological invariant} of the pair that can be bounded from below
by looking at the ramified coverings of $X$. Namely, if there is such a covering $\tilde X$
with ramification of order $\leq k$ that admits {\it no smooth metric with
positive scalar curvature} then, clearly, $\alpha_{min} (X,Z)\geq 2\pi/k$. 
\vspace {1mm}

{\it Question.} Are there pairs $(X, Z)$ where $\alpha_{min} (X,Z)$ is finite but yet not of the kind  $ 2\pi/k$ for any integer $k$?
\vspace {1mm}

The simplest instance of where such bounds are available  is where $X=S^2\times S^1$ and  where $Z$ equals a union of "coordinate circles",  i.e. $X=Z_0\times S^1$ for a finite subset $Z_0$ in the sphere $S^2$ and where the geometries near these circles $Z_i=z_i\times S^1 \subset X=S^2\times S^1$, $z_i\in Z_0$, are corner singular with (not necessarily mutually equal)   angles $\alpha_i$.

On can show, as we do it  in section 5.4 that, albeit singular, $X$ contains a {\it minimal surface} $Y$ in the homology class of the $2$-sphere $S^2\times z_0\subset X$ to which   the Gauss-Bonnet theorem applies and yields the (sharp!)  inequality 
 $$\sum_i (2\pi-\alpha_i) \leq 4\pi. $$  

{\it Remark.} Some  geometry questions on  general  cornered manifolds $P$ with $scal(P)>0$
reduce to those about the above  "pure edged"   spaces with the multi-doubling procedure from section 1.1. \vspace {1mm}

{\it Problem.} What are complete singular Riemannian  spaces   that are locally isometric to Cartesian   products of   flat manifolds and $2$-dimensional Riemannian cones?

Does, for instance, every stably parallelizable manifold admit such a singular Riemannian metric?

\section  { $\mu$-Area and $\mu$-Bubbles.}

 The  
 variational approach  indicated in section 1.4  for construction of Plateau-hedra
 also applies to more general curve-faced  polyhedra with prescribed mean curvatures of the faces and given dihedral angles at the edges as follows. \vspace {1mm}

An  open or closed subset $U\subset X$  is a {\it domain} if its boundary  also serves as the topological boundary of the complement to the closure of $U$, that is 
$$\partial U=\partial( int( U))=\partial( clos(U))=\partial
(X\setminus U) =\partial (X\setminus clos(U))=\partial (X\setminus int(U)).$$

  Each  component $Y$ of the boundary of a domain $U$ admits two coorientations  represent by arrows  directed toward interior or exterior of $U$. Accordingly,  we denote by    $[...\overset {in} \leftarrowtail Y] \subset X$ the germ of the intersection of $U$ with an arbitrarily small neigbourhood $W_\varepsilon$ of $Y$ in $X$, where we agree that $Y\subset [...\overset {in} \leftarrowtail Y] \subset X$  and  we denote by $[Y \overset {out} \leftarrowtail...]$ such a germ of the exterior of $Y$ that is the complement of the interior  of $U$ in $W_\varepsilon$.

Observe that this notation makes sense even if a hypresurface $Y$ does not bound anything in $X$, but only divide its {\it small} neighbourhood into  two parts, one regarded as interior and the other as exterior of $Y$.

A Borel measure $\mu$ on $X$ defines a  closed 1-cochain   on cooriented hypersurfaces $Y$ written is 
 $\mu[...\overset {in} \leftarrowtail Y]$, where the function  $Y\mapsto \mu[...\overset {in} \leftarrowtail Y] $  changes by an additive  constant
   under a  change  of the representative of the germ 
$[...\overset {in} \leftarrowtail Y] \subset X$.

If    $X$ is also endowed with a Riemannian metric  we define the {\it $\mu$-area} of $Y$ by
  $$area_{\mu}(Y)=_{def} vol_{n-1}(Y)-\mu [...\overset {in} \leftarrowtail Y] $$
 where $vol_{n-1}$ stands for the $(n-1)$-dimensional Hausdorff measure.

 Call a  hypersurface $Y \subset X$ a {\it $\mu$-bubble} 
  if it  {\it locally minimizes}  
 the function $Y\mapsto area_{\mu}(Y)$.

If $\mu$ is given by a continuous density function
 $\phi(x)$,  $x\in X$, i.e.  $ \mu =\phi dx$, then the mean curvature  of a
 $\mu$-bubble $Y \subset X$, obviously, 
satisfies $mn.curv(y)=\phi(y)$.  These {\it $\mu$-bubbles} are also called  $\phi$-bubbles.

In particular, {\it $C$-bubbles} with constant mean curvature equal  $C$  corresponds to   $\mu$ being proportional to the  Riemannian $n$-volume.
\vspace {1mm}

 {\it Example.} If  $X=\mathbb R^n$, $n\geq 3$, and $\phi(x)=(n-1)||x||^{-1}$, then the $R$-spheres defined by $||x||=R$ are (non-strictly) locally minimizing $\phi$-bubbles.

{\it Questions.}  Can a  (stable?)  {\it minimal} submanifold $Z$ of codimension $\geq 2$ in a Riemannian $n$-manifold$ \subset X$  be "surrounded" by   $C$-bubbles (with an arbitrarily large positive constant  $C$)  that are small {\it perturbations} of the levels of the  function 
   $(n-1)dist(x,Z)^{-2}$?  (This, possibly, can done by arranging suitable traps for such bubbles, see section 4.2.)

When, in general,  does a closed subset $Z\subset X$ admit arbitrarily small neighbourhoods $U\supset Z$ with (smooth?)  boundaries $\partial U$ of  (almost?) constant mean curvatures?

For instance, when can a subset  $Z\subset X$  be surrounded 
by  $\phi$-bubbles that approximate levels of a function  $\phi\sim (n-1)dist(x,Z)^{-2}$?

Is this possible  for  {\it smooth non-minimal} submanifolds and/or for   {\it singular minimal}  subvarieties $Z\subset  X$ with $codim(Z)\geq 2$?

What happens in this regard to  piecewise smooth subpolyhedra and to real algebraic subsets  $Z$?

Is there anything of this kind for  nice(?),  possibly   fractal, subsets  $Z\subset X$,   e.g. for singular loci of minimal subvarieties in $X$?

Do  constructions of minimal and mean curvature convex  hypersurfaces  from \cite {plateau-stein} extend to  $\phi$-bubbles in non-complete
manifolds such as $X\setminus Z$  with suitable functions $\phi(x)$  that have  "pole-like singularities" 
on $Z$?

 \vspace {2mm} 
   
   If $X$ is a manifold with a smooth boundary $\partial X$ and   $\mu=\mu_\partial$ is
 given by a continuous density function on $\partial X$, say by $\psi_\partial(x')$, $x' \in \partial X$, then 
  $\mu_\partial$-bubbles $Y\subset X$ with $\partial Y \subset \partial X$  meet $\partial X$ at the (dihedral) angle $\angle =\angle (Y\pitchfork \partial X)$ that,
     obviously, satisfies 
$\cos\angle (y')=\psi_\partial(y')$ for all $y' \in Y\cap \partial X.$
Moreover, this  equality remans true for the $\mu'_\partial$-bubbles
where $\mu'=\mu_\partial +\phi\cdot vol_n$ for  
 $\phi=\phi(x)$ being  an arbitrary continuous function on $X$.

\vspace {1mm} 

{\it Question}. What is the minimal regularity of a measure $\mu$ needed 
for the existence of
$\mu$-bubbles and  their regularity comparable to that of
 minimal hypersurfaces? Does the condition $\mu\leq const\cdot Hau_{n-1}$
for the $(n-1)$-Hausdorff measure  $Hau_{n-1}$ suffice?

 \vspace {1mm}

  {\it Remarks.}   (a) If $X$ has trivial $(n-1)$-homology,  $H_{n-1}(X)=0$, then 
   minimal $\varepsilon$-bubbles are associated   
to the {\it  supporting lines $ a- \varepsilon v=const$} of the ({\it convex hull} of the) 
 {\it isoperimetric profile}  
of $X$ in the positive quadrant of the $(a,v)$-plane, where
$profile_{isop}(X)\subset \mathbb R_{++}^2$ is defined as   the    image of the
 map from the space of compact domains $P$
  to the plane
given by 
$$P\mapsto \big(a=vol_{n-1}(\partial P), v=vol_n (P)\big),$$ 
where, observe, the boundary of $profile_{isop}(X)$ is contained in the critical set (curve)
of this map.

(If  $X$ is a Galois coverings of compact manifold,  one usually works with the  {\it Foelner-Vershik} profile
that is the convex hull of the logarithmic map $P\mapsto \big(\log vol_{n-1}(\partial P), \log vol_n (P)\big)$.)

\vspace {1mm}

\vspace {1mm}
(b)   $\mu$-Bubbles  are well defined for closed hypersurfaces
  $Y$ which  bound 
 (non-compact) domains with {\it infinite} $\mu$-measures, like
  $ Y\times 0\subset X=Y\times \mathbb R$, provided 
 the $\mu$-measures of the regions between $Y$ and hypersurfaces
 $Y'$ homologous to $Y$ are finite.
 
  .

\subsection{Poly-Bubble-Hedra.}

A Riemannian $n$-manifold $P$ with corners, e.g.  a polyhedral domain 
in an ambient Riemannian manifold $X=X^n$, is called a  { \it PB-hedron}  
if all its $(n-1)$-faces $Q_i$
have {\it constant} (possibly mutually non-equal) mean curvatures $m_i$
 and the
 dihedral angles $\angle_{ij}=\angle[Q_i\pitchfork Q_j]$ 
  are constant on all $(n-2)$-faces
 of $P$.

For example,  the domains $P$ in spaces of constant curvature which are  bounded by 
(convex or concave)
 {\it umbilical} hypersurfaces (where all principal curvatures are constant
  and mutually equal,   e.g.  as it is for  spheres)  are PB-hedra.
  
  The combinatorially simplest  PB-hedra are {\it di-B-hedra} with two $(n-1)$-faces $Q_1$ 
  and $Q_2$
  meeting across a single $(n-2)$-face $Q_{12}=Q_1\cap Q_2$. Probably,
  the space of {\it di-B-hedra} $P\subset X$ with $ Q_{12}$ contained in a  hypersurface $H=H^{n-1}\subset X$ is Fredholm, i.e.   it locally has finite,  positive or negative, 
  "virtual dimension" $d$:   if $\cal H$ 
 is a generic $N$-dimensional family of hypersurfaces in $X$ with large $N$, then the space $\cal P$ 
 of  {\it di-B-hedra}
  $P\subset X$ which are close to a given $P_0$ and have $Q_{12}$ contained in some 
 $H\in \cal H$, satisfies $dim({\cal P})=N+d$. (Possibly, this $d$ may depend on whether
a  hypersurface $H$ separates $Q_1$ and $Q_2$ or not.)
  
  \vspace {1mm}

  {\it The  geometric  type} of  a $P$ is, by definition, the totality 
 of the numbers  $\{M_i,  \alpha_{ij}\}$  for $M_i= mn.curv (Q_i)$ and 
 $  \alpha_{ij}=\angle{[Q_i\pitchfork Q_j]}$
 associated to the $(n-1)$-faces and $(n-2)$-faces of (the combinatorial scheme of)  $P$. 

\vspace {1mm}

  \textbf {${\cal PB}$-Problems.} The space   ${\cal PB}(X)_{\cal GT}$ of  $PB$-hedra of 
  a given geometric type $\cal GT$
   (which includes the combinatorial type) in a given Riemannian 
manifold  $X=X^n$ is similar in many respects to the 
space  ${\cal MIN}(X)$ of closed minimal hypersurfaces in $X$, albeit the spaces
  ${\cal PB}(X)_{\cal GT}$ are  infinite  dimensional at
   certain $P\in {\cal PB}(X)_{\cal GT}$  and it may be 
 hard to decide, for example,
if   ${\cal PB}(X)_{\cal GT}$ is non-empty for given $X$ and   $\cal GT$.

Also, the compactness properties of (subspaces in)    ${\cal PB}(X)_{\cal GT}$ 
are less apparent than these in ${\cal MIN}(X)$, since sequences of PB-hedra 
$P_i\in  {\cal PB}(X)_{\cal GT}$ may Hausdorff 
converge to subsets $Z \subset X$ with $dim(Z)\leq n-1$, where the picture is not fully clear
even for decreasing sequences $P_0\supset P_1\supset P_i\supset...$
where the intersection $Z=\bigcap_iP_i$ may be(?) rather complicated, say for sequences 
of acute Plateau-hedra  of the combinatorial type of the $n$-cube. For example,

\vspace {1mm} 

What are smooth $k$-dimensional submanifolds $Z\subset X$ 
with boundaries which  can be represented as such intersections $\bigcap_iP_i$?

\vspace {1mm}

Another new feature of Plateau-hedra (and of  more general PB-hedra)
 is their dependence 
on the underlying combinatorial scheme $\cal GT$:   the space ${\cal PB}(X)$
of all Plateau-hedra is stratified by various  ${\cal PB}(X)_{\cal GT}$ according to
their combinatorial types, (similarly to the space of convex Euclidean polyhedra), but 
 the topology/geometry of this stratification is far from transparent.

  \vspace {1mm}
  
 The concepts of  {\it dihedral  extremality and rigidity} we met earlier   obviously  extend (as in  the problem $(\ast)$ stated at the end of  section 2.2) to  PB-hedra $P$  in general  Riemannian    manifolds $X$  where they seem particularly interesting in manifolds with constant curvature and where they  (partly) generalize    
   rigidity phenomena for hyperbolic warped products  \cite {positive}.

  \subsection 
{Multi-bubble Description  and Construction of PB-hedra.} 
 Let ${\cal P}_{\cal CT}$ be the space of 
polyhedral domains    $P$ in $X$ with $m$ faces meeting according  
a given combinatorial  pattern (scheme/type) $\cal CT$.

Denote
by  $v_i =v_i(P)$, $P \in {\cal P}_{\cal CT}$, the $(n-1)$-volume form (measure) 
on the face $Q_i$ of $P$, $i=1,...m$, regarded
 as the measure on $X$. 
 
Let $-\infty < \varepsilon_i < \infty$ and $ -1 < c_{ij} < 1$ be given 
constants,  let
 $$\mu_i=\varepsilon_i \cdot vol_n +\sum _jc_{ij}v_j,$$ 
 where $vol_n$ is the Riemannian volume (measure) in $P$ and
 where the sum is taken over the faces $Q_j$  adjacent to $Q_i$.

It is clear that:

\vspace {1mm}

 {\it if, for a given $i$, the face $Q_i$ is a local $\mu_i$-bubble, i.e. 
 if $Y_i$ locally minimizes the $area_{-\mu_i}(Y_i)$ 
then $Q_i$ has constant mean curvature $ =\varepsilon_i$ and the  dihedral
 angles between  $Q_i$ and $Q_j$  equal 
 $\arccos (c_{ij})$ for all $Q_j$ adjacent to $Q_i$,

$$  mn.curv(Q_i)=\varepsilon _i   \mbox  {  and   } 
\cos\left (\angle[Q_i\pitchfork Q_j]\right )=c_{ij}.$$} 

\vspace {1mm}

{\it Minimal  PB-hedra.} Let $P^\circ$ be a Riemannian manifold with 
corners, e.g. a polyhedral domain in a manifold $X=X^n$, such that
  the mean curvatures of the $(n-1)$-faces  $Q^\circ_i \subset P^\circ$ and the  dihedral angles between them satisfy
   $$ mn.curv(Q^\circ_i)\geq \varepsilon_i \mbox { and } 
    \angle(Q^\circ_i\pitchfork Q^\circ_j)\leq \alpha_{ij} \mbox
   { for given }  \varepsilon_i\geq 0 \mbox { and } \alpha_{ij}\leq \pi/2$$
   and let
   $$\mu_i=\varepsilon_i \cdot vol_n +\sum _j (\cos  \alpha_{ij})v_j.$$
    
 Now,  take some $i_1$ and  minimize the the $\mu_{i_1}$-area of of the $i_1$-th 
 $(n-1)$-face
  by varying  
 $Q^\circ _{i_1} \rightsquigarrow Q_{i_1}$ {\it within} $P^\circ$ and with the
 boundary $\partial Q_{i_1}$ contained in the union of
  the $(n-1)$-faces adjacent to $Q^\circ_{i_1}$
while keeping all $(n-1)$-faces $Q^\circ_{i\neq {i_1}}$ unchanged. 
 
  In other words,   consider all closed subsets $P\subset P^\circ$ with 
  m.c. convex boundaries  $\partial P$ that
  contain all faces of $P^\circ$ which do not intersect $Q_i^\circ$, 
  and then minimize the $\mu_{i_1}$-area of the
    new part of $\partial P $, that is $Q_{i_1}=\partial P \setminus \partial P^\circ.$
  
 Thus we obtain some subset, say $P^{[{i_1}]}\subset P^\circ$ which we regard 
 as m.c. polyhedron with a minimal $i$-th face $Q[i_1]_{i_1}=Q_i^{min}$
  and with new faces
 $Q_i^{[i_1]}=  P^{[i_1]}\cap Q_j$ for $i\neq i_1$.

 \vspace {1mm}

This $P^{[i]}$ can be more singular than it is  allowed by our definition 
of "polyhedral domain",
but we pretend {\it it is} such a domain  and 
then apply the same minimization process to $P^{[i_1]}$ with respect to some 
face $Q^{[i_1]}_{i_2}$ for $i_2\neq i_1$.

By  continuing this process  with a given sequence $i_1,i_2,...$ one arrives at
 a decreasing family of mean curvature convex "polyhedra"

$$...  \subset P^{[i_1,i_2, i_3]} \subset  P^{[i_1,i_2]}  \subset P^{[i_1]}  \subset P^\circ$$

If this family stabilizes at some $ P_{min}(k)=P^{[i_1,i_2, i_3,...i_k]}$  and this $P_{min}(k) \subset X$ qualifies
 as
a polyhedral domain, it serves as a minimal PB-hedron inside $P^\circ$.

\vspace{1mm}

  \section {   Plateau Traps, $\phi$-Convexity, Qiasiregularity and Regularization.}
 
Let $P$ be a cubical polyhedral domain (or a manifold with corners) let
 $\pm Q \subset P$ be 
a pair of opposite faces and $X= P\setminus (Q\cup-Q)$ be the cornered manifold with the boundary
$ \partial X= \partial P\setminus (Q\cup-Q)$ that is the union of all $(n-1)$-faces except $\pm Q$.

If $P$ has  strictly acute dihedral angles and strictly  mean (curvature) convex faces,   then the minimum of the $(n-1)$-volumes of hypersurfaces $Y\subset X$
with $\partial Y\subset  \partial X$ is assumed by some $Y_{min}$ {\it away} from these faces {\it inside} $P'$.
 
 One may say that the relative homology class  $[Y]\in H_{n-1}(X,\partial X)$ is {\it  trapped 
 between $Q$ and $-Q$}. 
 
 The trapping feature, if shared  by  all pairs of opposite faces in $P$, is {\it equivalent} to
  the [acute angles + mean convexity] property of $P$. On the other hand, traps are quite robust. In particular the definition of traps  needs significantly   less
    regularity than that of mean curvature,  be it the  topology of $P$, its    Riemannian metric
 or smoothness assumption on the faces of $P$. 

In what follows we review simple standard properties of traps;   we limit ourselves  for the most part of the exposition to
{\it closed} hypersurfaces representing {\it absolute} $(n-1)$-homology classes in order not to overburden our  notation.

\vspace {1mm}

 \subsection {  Directed Homology. }  Let $X$ be an $n$-manifold, possibly with a    boundary.
 A {\it direction}   $ \overset{ \leftarrowtail}C_{n-1}$ in $X$ is, by definition, 
 a ({\it directing})   homology class 
 which is   representable  by a closed  {\it cooriented} hypersurface $Y\subset X\setminus \partial X$
which divide $X$ into two closed  subsets denoted 
$U_{in} \subset X$ and $U_{out} \subset X$
where 
$$ \partial  U_{in}=\partial U_{out}=Y= 
 U_{in}\cap U_{out}.$$

 {\it Remark.}  If $X$ is oriented, then  $ \overset{ \leftarrowtail}C_{n-1}$  is an 
ordinary homology class, but, in general, it is,  strictly speaking, a {\it $1$-cohomology} class.
We say "homology" to emphasize  geometric representations of  $ \overset{ \leftarrowtail}C_{n-1}$ by hypersurfaces even 
for non-orientable $X$,  e.g. where $X=Y\times \mathbb R$ and $Y$ is non-orientable.

\vspace {1mm}
{\it Relative Case.} Let $\partial_\circ X \subset  \partial X$ be an open subset in the boundary of $X$, e.g. $\partial_\circ X = \partial X$. Then there is an obvious generalization of the above  to hypersurfaces $Y \subset X$ with $\partial Y \subset \partial_\circ X$.  

The directing homology class is represented in the general relative  case by non-closed cooriented hypersurfaces $Y$ that are still are closed as subset in $X$ and that may have boundaries
contained in the   $\partial_\circ X$-region of the boundary of $X$. 

If all infinity of our $Y$ is contained in  $\partial_\circ X$, everything  can be reduced to the  absolute case   by taking the double of $X$ along    $\partial_\circ X$. On the other hand, such doubling does not (quite) apply to a more general setting, e.g. where $X$
has no boundary at all and where   $\partial_\circ X$ is a  "virtual subset at infinity" represented by
a descending family of open subsets in $X$,  say  $ U_1 \supset U_2 \supset...\supset U_i \supset... $, such that the intersection of every compact subset  $ K\subset X$ with $U_i$ is empty for all $i\geq i(K)$.

The  relevant hypersurfaces $Y\subset X$ here are those where the difference sets
$Y\setminus U_i$ are compact for all $i$.

  \vspace {1mm}
  
  Given a Riemannian metric on $X$ and  a Borel measure  $\mu$, 
  define as earlier 
$$ area_{\mu}(Y)= vol_{n-1}(Y)-\mu [...\overset{in} \leftarrowtail Y] \mbox { for } n=dim(X),$$
for all cooriented  $Y$ in the class $ \overset{ \leftarrowtail}C_{n-1}$ as we did in the previous section.

\subsection { Traps and  Walls}. 

   A domain
 $U\subset X$ ,
is called a {\it $\mu$-trap  (or well) for $ \overset{ \leftarrowtail}C_{n-1}$ in $X$}   and its boundary is called a {\it $\mu$-wall }
if   $ \overset{ \leftarrowtail}C_{n-1}$ can be represented by a smooth cooriented hypersurface in the interior $int(U)$ and  if  every
 smooth cooriented hypersurface 
$Y \subset U$  which represents  $ \overset{ \leftarrowtail}C_{n-1}$, i.e. $[Y]= \overset{ \leftarrowtail}C_{n-1}$,  (such a $Y$ separates the two $ \overset{ \leftarrowtail}C_{n-1}$-ends in $X$)
can be  
  "moved" to $ Y'\subset int(U)$, such that
 $$ area_{-\mu}(Y') \leq area_{-\mu}(Y)$$
where "moved" signifies that $Y'$ is a smooth cooriented closed hypersurface {\it homologous} to $Y$, 
i.e. $[Y']=[Y]= \overset{ \leftarrowtail}C_{n-1}$ and  
  this inequality is sharp 
unless $Y\subset int(U)$ to start with.

If $\mu=0$ these are called {\it Plateau traps/walls}  or just traps and walls;  if $\mu=\phi dx$ we speak of  {\it $\phi$-traps} and {\it $\phi$-walls} for functions $\phi$ on $X$.

\vspace {1mm}

{\it Remarks.} (a) The topologically  simplest situation, and this is the one we mostly deal with, is where $X$ is compact 
with two boundary components, e.g. a cylinder  $X=Y\times [0,1]$,  and where our
 hypersurfaces in $ X$ separate these components.

(b) More general singular measures $\mu$,  e.g. those supported on the boundary $\partial X$ are  relevant for PB-hedra and  the corresponding traps for the {\it relative} homology $H_{n-1}(X, \partial_\circ X)$ for some  subset $\partial_\circ X \subset \partial X$.

\vspace {1mm}

 {\it Smooth Mean Convex Traps.} If  the boundary of a domain $U$ has  strictly positive  mean curvature  then  $U$ it traps all classes  $ \overset{ \leftarrowtail}C_{n-1}$  that have representative cycles in $U$. More generally (and equally obviously), if $U$ is preconvex and the mean curvatures of the faces are  strictly bounded from below by a function $\phi(x)$
 then $U$ is  $\phi$-traps. \vspace{1mm}

{\it Normal Traps.} Let $Y\subset X$ be a closed  
smooth cooriented hypersurface that represent a non-zero directed homology class in $X$   and $\phi(x)$ be a $C^1$-smooth function
such that $\phi(y)=mn.curv_y(Y)$ for all $y\in Y$. 

 If {\it the inward normal derivative $\frac {d\phi(y)}{d\nu^{in}}$ on $Y$,} (for   $\nu^{in}$ being the inward looking unit normal vector  along $Y$) {\it  is sufficiently large} , say
$$\frac {d\phi(y)}{d\nu^{in}_y}>  curv_y^2( Y)  +|Ricci_X(\nu^{in}_y, \nu^{in}_y)|\mbox { for all } y \in Y,$$
where $curv^2$ denotes the sum of squares of the principal curvatures of $Y$ and where, observe, $Ricci(\nu^{in}, \nu^{in})=Ricci(\nu^{out}, \nu^{out})$,  then

\vspace {1mm}
{\it the homology class  $\overset{ \leftarrowtail}C_{n-1} =[Y]$ in $U_0$ is $\phi$-trapped in all neighbourhoods 

$U \subset X$  of $Y$.}

\vspace {1mm}

Indeed, let $d(x)$ be a smooth function in a neighborhood of $Y$ without critical points that is negative  inside $Y$, i.e. on $[...\overset {in} \leftarrowtail Y]$,   and positive on  $[Y \overset {out} { \leftarrowtail}...]$. If the mean curvatures of the levels of this function  cooriented as $Y$  strictly minorize  $\phi(x)$ strictly inside $Y$ i.e.  for   $x\in [...\overset { \leftarrowtail}Y]\setminus Y$ and strictly majorize   $\phi(x)$   strictly outside $Y$, then, obviously the homology class of $Y$ is trapped 
in the domain $d^{-1}[-\varepsilon, \varepsilon] \subset X$ for all (arbitrarily small) $\varepsilon >0$.

  If this applies to the signed distance function  $d(x)=\pm(x)dist(x, Y)$ negative inside $Y$ and positive outside, then the sufficiency of the  above lower bound on $\frac {d\phi(y)}{d\nu^{in}_y}$ follows from  the {\it second variation formula} for $vol_{n-1}(Y)$.

\vspace {1mm}

    {\it  Locally Trapped Hypersurfaces with Boundary.}  Let $Y$ be a smooth compact hypersurface with a boundary $\partial Y \subset \partial X$, let $\phi_\circ : \partial X \to [-1+\beta,1-\beta]$, $\beta>0$, be a  $C^1$-smooth function
and let $\mu =\phi \cdot vol_n + \phi_\circ vol_{n-1}$ for  $vol_{n-1}$ referring here to the Riemannian volume (measure) on the boundary $\partial X$. 
 
 Let the (dihedral) angle $\angle_y=\angle_y (Y \pitchfork\partial X)$ satisfy 
 $$\cos \angle_y=\phi_\circ(y)\mbox  {  for  all $y\in \partial Y$} $$

 \vspace {1mm}
 {\it If the inward normal derivative $\frac {d\phi(y)}{d\nu^{in}}$
 is  sufficiently large, where
  the lower bound depends, besides $curv^2( Y)$  and $Ricci(X)$, on
   $curv^2(\partial Y)$, on $( 1-|\phi_\circ|)^{-1} \leq \beta^{-1}$, 
  and on the normal derivative of $\phi_\circ$ in $\partial Y$, then  $ Y$   
  minimizes  $area_\mu$ in every small neighbourhood $U_0 \supset Y$.}

\subsection {Smoothing and Doubling.} Let  $U$  be  a compact domain in a smooth Riemannian $n$-manifold $X\supset U$ that traps a directed homology class   $ \overset{ \leftarrowtail}C_{n-1}$. Then    \vspace {1mm}

{\it there exists a subdomain $U'\subset U$ with $C^2$-smooth boundary $Y'$ that also $\phi(x)$-traps
  $ \overset{ \leftarrowtail}C_{n-1}$ and such that the mean curvatures of $Y'$ satisfy
  $$mn.curv_x(Y')>0\mbox {  for all } x\subset Y'.$$}
 
 {\it Proof.} Since the class $ \overset{ \leftarrowtail}C_{n-1}$ is trapped in $U$ it can be realized by, a priori singular,  cooriented hypersurface $Y_{min} \subset int(U)$.
 
This   $Y_{min}$ may be non-unique; in this case we take the union of all these minimal $Y_{min}$, denoted  $Y_{MIN}=\cup Y_{min}$ that is a compact subset in the interior of $U$. 

Clearly, the class $ \overset{ \leftarrowtail}C_{n-1}$ is trapped in {\it all}  neighborhoods $V\subset U$ of  $Y_{MIN}$, and,
by continuity,  the class class $ \overset{ \leftarrowtail}C_{n-1}$ is  $\varepsilon$-trapped  in $V$ for all, positive and negative $\varepsilon =\varepsilon(V)$ with sufficiently small absolute values $|\varepsilon|>0$. Therefore,  $Y_{MIN}$ is pinched between two  $\varepsilon$-bubbles, $\varepsilon>0$, i.e.  it is contained in a domain $V_\varepsilon'\subset V$ bounded by these bubbles.

We conclude the proof by smoothing these bubbles keeping their mean curvatures positive (see \cite {plateau-stein} and thus, approximate $V'$ by the required smooth domain $U'_\varepsilon$. QED \vspace {1mm}

 \textbf {Scalar Curvature Corollary.} Let a domain $U$ that traps a homology class   $ \overset{ \leftarrowtail}C_{n-1}$  admit a continuous map to the $(n-1)$-torus, say $ f: U\to\mathbb T^{n-1}$, such that    $ \overset{ \leftarrowtail}C_{n-1}$ goes to a {\it non-zero} multiple of the fundamental class of the torus. 
 If $U$ is spin, then \vspace {1mm}
 
 {\it the scalar curvature of $X$ is strictly negative at some point $x\in U$.} \vspace {1mm}

{\it Proof.}  As one knows (see  \cite {gromov-lawson1}  that if $scal(U_\varepsilon')\geq 0$, then   the  double, say  $U'_\varepsilon{^+}'$ of $U'_\varepsilon$  is a closed manifold  that admits a smooth  metric 
 with   positive  scalar curvature. Then  the proof follows from the validity of the Geroch conjecture in the spin case, \cite {gromov-lawson1}. \vspace {1mm}

{\it Remarks.} (a)  The above applies whenever  the topology of (the closed manifold)
 $U'_\varepsilon{^+}'$ does not allow  metrics on it  with $scal>0$ on it   that is known in a variety of  cases. In particular, the Schoen-Yau theorem allows one  to suppress the spin condition for $n\leq 7$.

 It was conjectured by Brian White (a private communication about 20-25 years ago) 
that \vspace {1mm}

\hspace {12mm}  {\it  singularities of minimal hypersurfaces are unstable.} \vspace {1mm}

This would allow an extension of the Schoen-Yau method to all $n$ (with some problems remaining in proving rigidity results).

 White's  conjecture was confirmed in  \cite {smale} for  $n=8$, where, observe, singularities are isolated.

The case $n\geq 9$ remains open, but we shall  see  in section 5.3 how  to extend  Schoen-Yau results to  {\it non-spin} manifolds of dimension  $n=9$ where the singularities  $\Sigma$  of minimal hypersurfaces are most $1$-dimensional. 

Furthermore,   Lohkamp's method of "going around singularities"    applies to all $n$,  but I have not studied it  in depth and can not   apply it to the problems at hand.  \vspace {1mm}

(b)  Generically    $Y_{min}$ is unique and then the (smooth mean curvature convex) domain $U_\varepsilon '\supset Y_{min} $ can be chosen arbitrarily close to $Y_{min}$.
 
 The  geometry of such  $U'_\varepsilon$  is 
similar to the warped product metric on $Y\times \mathbb R$ from \cite {gromov-lawson3}
Apparently, the former converges in a suitable sense to
 the latter  for $\varepsilon \to 0$  as  the neighborhood $U'_\varepsilon$, that converges to  $ Y_{min}$, becomes "infinitesimally narrow".

In fact the arguments using the warped product from \cite {gromov-lawson3} can be trivially adjusted to 
$U'_\varepsilon$ and/or to its double  $U'_\varepsilon{^+}'$. 

The advantage of this  over  the warped product is that it make sense for {\it singular} $Y_{min}$ but  we have only managed  a limited use of it (such as handling  dimensions $n=8,9$ for non-spin manifold).

Also  $U'_\varepsilon$  it is harder (but, probably, possible) to use it with bubbles $Y$
in spaces with negative scalar curvature as we did in \cite {positive} with a help of warped products. \vspace {1mm}

 {\it Smoothing  at the Edges.} Let $P$ be a  {\it singular} polyhedral domain in a Riemannian manifold $X$,
 i.e. where the faces may have singularities of the kind minimal hypersurfaces have. 
 
 For instance, Let  $P$ be bounded by a pair of  hypersurfaces, say $W_0$ and $W_1$ where
 $W_0$ is  a minimal hypersurface  and $W_1$ is a $\phi$-bubble  with free boundary $\partial W_1\subset W_0$ for a   continuous function  $\phi$ on $X$.
 
 The dihedral angle between $W_0$ and $W_1$ equals $\pi/2$ at the points $x$
 at the edge  $\partial W_1=W_0\cap W_1$ where both $W_0$ and $W_1$ are regular.
 
 But when it comes to singular points it is not even clear if the concept of dihedral angle is  symmetric, i.e. if  $ \angle(W_0,W_1)= \angle(W_1,W_0)$.

Both $W_0$ and $W_1$ can be approximated by smooth hypersurfaces with a minor decrease of their mean curvatures   \cite {plateau-stein}
 (in fact, with an increase if $Ricci(X)\geq 0$). Such a  smoothing starts with a small  equidistant inward deformation of $W_0$ and $W_1$.  If $W_0$ and $W_1$ are smooth
 at the edge points, then the dihedral angle also changes little but  the dihedral angle may, a priori, uncontrollable increase near the singularities in the edge.
 
{\it Question.}  Is it possible to smooth the faces, with  at most   $\varepsilon$-decrease of their mean curvatures and with  at most  $\varepsilon$-increase of the dihedral angle(s)
 for  an arbitrarily small $\varepsilon>0$?

 \subsection{$\Delta$-Stable Mean Convexity.}

We want to 
define a concept of  {\it mean curvature bounded from below}, in particular, of (strictly and non-strictly) 
positive mean curvature for non-smooth domains $U\subset X$ and  want this positivity to be stable under small continuous perturbations of $U$ as well as under $C^0$-perturbations of the Riemannian metric in $X$.  Eventually,  
we want a concept   adaptable to    singular spaces. There are several candidates for such a stable mean convexity.   Below is  a instance of such a  definition that depends on Almgren's concept of 
\vspace {1mm}

{\it $\varepsilon $-Minimization.}     Let $H\subset X$ be a smooth cooriented hypersurface, let  $H_0\subset H$ a compact domain in $H$ and $H_1=H\setminus H_0$.

  {\it A (smooth)   $\varepsilon $-minimization of $H$  supported in $H_0$}, denoted $H\rightsquigarrow H'$, 
  is  a replacement of $H$ by another smooth cooriented hypersurface  $H'$  such that $H'=H'_0\cup H'_1$, where
  $H'_1 =H_1$, where the  union $H_0\cup H'_0$ makes a closed hypersurface that {\it bounds} in $X$ 
 (and so $H'$ is homologous to $H$)  and such that
  $$ vol_{n-1}(H_0')\leq e^{- \varepsilon}vol_{n-1}(H_0),
  $$
 where we assume, to avoid immaterial   complications, that the boundary of $H_0$ in $H$ has Hausdorff dimension at most $n-2$.
 
  Clearly, \vspace {1mm}
  
  {\it composition of     $\varepsilon $-minimizations
  $$H\rightsquigarrow H'\rightsquigarrow H''$$  is an $\varepsilon $-minimization
  $$H\rightsquigarrow H''$$
the support of which equals the union of the supports of its factors}.\vspace {1mm}

Also observe that this minimization is stable under small perturbations of the Riemannian metric $g$ on $X$:

{\it if 
$$\lambda^{-2} g  < g_1 <    \lambda^2 g\mbox   { for } \lambda^{n-2}\leq e^\varepsilon   $$  
then $H\rightsquigarrow H' $ is  an $\varepsilon_1$-minimization in $(X, g_1)$ for $\varepsilon_1=\varepsilon- (n-2)\log\lambda.$}

 \vspace {1mm}

{ \it $\Delta$-Stable Mean Convexity.} Let 
$U$ be  a closed domain in $  X$ with compact  boundary $Y=\partial U$, let $V_1, V_2\subset U$ be closed subsets in $U$ that contain   $Y$ and  such that $V_1$ is contained in the interior of $V_2$; let $\Delta=V_2\setminus V_1$.  


Say that $Y$ is ({\it multiplicatively}) {\it $(\varepsilon, \Delta)$-stably  mean convex} if there exists a (small) neighbourhood $W\subset V_2$ of $Y\subset V_2$ such that
every 
smooth cooriented hypersurface $H \subset V_2$ that intersect $V_1$   admits an 
$\varepsilon $-minimization 
$H\rightsquigarrow H'$, such that 

$\bullet_1 $   $H'$  does not intersect $W$;

$\bullet_2 $ the support $H_0\subset H$ of the $\varepsilon $-minimization $H\rightsquigarrow H'$  is contained in $V_2$;

 $\bullet_3$ the hypersurfaces $H$ and $H'$ coincide outside $V_2$, i.e. the "new part" 
 
 $H_0'\subset H'$ of $H'$ is contained in $V_2$.




\vspace {1mm}
 
{ \it Remark.} What is most essential here is that the volume of  $H'_0$   is smaller  than that of $H_0$ by a  {\it definite   amount}, (roughly by  $\varepsilon\cdot vol_{n-1}(H_0)$) that is   independent of $H$: this makes this mean convexity stable under small perturbation of the metric in $X$.

On the other hand,  keeping  $H'$ away from $Y$  is a minor issue since   boundaries of {\it domains}, as we define them,  are  large.  (Nothing like   $U=X\setminus Z$ with 
$codim (Z)\geq 2$ is allowed.) Yet, this is needed to make this convexity stable under small perturbation of $Y$.

\vspace {1mm}
 
($\ast$) {\it Mean Convex $\Rightarrow$ $\Delta$-Stable Mean Convex.}  Let $U$ be a domain in a  complete Riemannian manifold $X$,    denote by  $U_{-\delta}\subset U$  the subset  of those  $u\in U$  where $dist(u,Y)\geq \delta$ and let $Y_{-\delta}=\partial U_{-\delta}$ be the corresponding equidistant hypersurfaces.

Let  $Y$ has strictly positive mean curvature and observe  that  $Y_{-\delta}=\partial U_{-\delta}$ also have positive mean curvatures for small $\delta< \delta_0=\delta_0(Y)>0$. \vspace {1mm}

{\it There exists an $\varepsilon=\varepsilon(U, \delta) >0$ for all positive  $\delta <\delta_0$ 
such that $U$ is $(\varepsilon, \Delta)$-stably mean convex for   
 for  $\Delta=U_{-\delta_1}\setminus U_{-\delta_2}$ and all
$0< \delta_2<\delta_1\leq\delta $.}
 
 \vspace {1mm}
 
 {\it Proof.} One may assume (this is easy, compare \cite {almgren},  \cite {filling})   that $H\cap V_2$ is {\it  $\varepsilon$-minimal} that is  no   $\varepsilon$-minimization (for a suitably small $\varepsilon>0$) that satisfies the above $\bullet_2$ and $\bullet_3$ exists. This property  provides a {\it lower bound} on the $(n-1)$-volumes  of $H$ intersected with
 $\rho$-balls $B_h(\rho)\subset \Delta$, see  \cite {almgren}.   It follows, that  there exists $\delta'$ in the interval $ \delta_2 <\delta'< \delta_1$ such that the region $H'_0 \subset Y_{-\delta'}$ bounded in $Y_{-\delta'}$ by $H\cap Y_{-\delta'}$ satisfies:
 $$ vol_{n-2}(H'_0)\leq e^{-\varepsilon} vol_{n-2}(H_0)\mbox { for  } H_0=H\cap (U\setminus U_{-\delta'}).$$
QED. \vspace {1mm}

($\ast\ast$)  {\it   $\Delta$-Stable Mean Convexity $\Rightarrow$  Mean Convexity.}  Let a closed domain  $U\subset X$ be $(\varepsilon, \Delta)$-stably  mean convex for some $\Delta=V_2\setminus V_1\subset U$ and $\varepsilon>0$. \vspace {1mm}
 
{\it  Then there exists a  smooth mean curvature convex  subdomain $U_1$ in $U$  such that  $V_1\subset U_1\subset V_2$; thus,  the (compact smooth  mean curvature convex)  boundary $Y_1=\partial U_1$ is contained in 
 $\Delta$.}
 
 \vspace {1mm}
 
 {\it Proof.} Let  $\phi(x)$  be   a non-negative  continuous function   that is very small on $V_1$   and that  fast grows  in $\Delta$ as $x$ approaches the complement $U\setminus V_2$. Then, the $\Delta$-stable mean convexity of $U$ implies   the existence of  a $\phi$-bubble trapped inside $\Delta$ and the required $Y_1$ is obtained with a smooth approximation of this bubble, see \cite {plateau-stein}. 
 
    \vspace {1mm} 
 
  {\it Remark/Question.} This arguments relies on the basic regularity theorems of the geometric measure theory. Is there an elementary proof of the  implication ($\ast\ast$)?   \vspace {2mm}

  {\it On Relative $\Delta$-Stable Mean Convexity.}   The  definition  of $\Delta$-stable mean convexity obviously generalizes to the relative case  where (possibly non-smooth) hypersurfces $Y$ have  boundaries  that must be  contained in a given closed subset $ W\subset X$.  The so defined notion of  "mean convexity" is stable under small perturbations of $Y$ and of the metric $g$ on $X$, but... not under $C^0$-perturbation of $W$.
 
For instance, if $W$ is a smooth  closed  hypersurface in $X$ that bounds a domain $U\subset X$ and $W_\varepsilon$ is obtained by a smooth diffeotopy
that keeps $W_\varepsilon$  all along in a small  $\varepsilon$-neighbourhood of $W$ in $X$, 
then the local geometry of $X$ near
$W_\varepsilon $ is {\it not necessarily} close in any sense to the local geometry of $X$ near $W$. Because of this, even if $Y\subset U$ is locally isolated minimal submanifold with (free) boundary 
$\partial Y\subset W=\partial U$, one can not guarantee  that
the minimal $Y_\varepsilon \subset U_\varepsilon$  with $\partial Y_\varepsilon\subset W_\varepsilon$ will  be, in general, geometrically close to $Y$.

Of course, the bulk of $Y_\varepsilon$ away from $W_\varepsilon$ and also  for most of the boundary of 
$Y_\varepsilon$ will be close to $Y$ but if the local relative filling profile of $U_\varepsilon$ near $W_\varepsilon$ becomes  badly non-Euclidean, then
 $Y_\varepsilon$ may have develop "long narrow fingers"  in vicinity of $W_\varepsilon$ protruding somewhere   at the boundary of $Y_\varepsilon$ and spearing in  $U_\varepsilon$ along  its boundary  $W_\varepsilon$.    \vspace {1mm}

However,  we shall see presently that 
 if  the hypersurface   $W_\varepsilon $  serves as a bubble  for a $C^0$-small 
$\varepsilon$-perturbation $g_\varepsilon$  of the original Riemannian metric $g$ in $X$, then 

  {\it the domain $U_\varepsilon\subset X$ bounded
by $W_\varepsilon$ satisfies  almost Euclidean  filling

inequalities for relative $(k-1)$-cycles $C\subset U_\varepsilon$ for all $k=2,3,..., dim(X)-1$.} \vspace {1mm}

  This  rules out   "fingers" and guarantees the existence of a relative  $g$-bubble  $Y_\varepsilon  \subset U_\varepsilon $ with $\partial Y_\varepsilon\subset W_\varepsilon=\partial U_\varepsilon$ that approximate a given submanifold  $Y\subset U$ with $\partial Y\subset \partial U$.
 
 Furthermore, if $U'_\varepsilon \subset U_\varepsilon$ equals a part of  $U'_\varepsilon$ bounded by 
  $Y_\varepsilon$ and a "half" of the boundary of $ U_\varepsilon$,  then one can, in some cases,  construct bubbles   $ Y'_\varepsilon \subset U'_\varepsilon$ with
  boundaries in  $\partial U'_\varepsilon$, where, observe  the boundary  $\partial U'_\varepsilon$ has a corner along the set  where  $Y_\varepsilon$ meets $W_\varepsilon=\partial U_\varepsilon$ and this corner may necessitate a presence of a singularity in  $ Y'_\varepsilon$ at this corner.
  
 The construction of such bubbles  $ Y_\varepsilon$, $Y'_\varepsilon$,  $Y''_\varepsilon$,... will be used   for   proving that

 \hspace {5mm}{\it $C^0$-limits $g$ of metrics $g_\varepsilon$ with $scal(g_\varepsilon )\geq 0$ also have $scal(g)\geq 0$. }\vspace {1mm}

 We shall do it in section 4  where we shall explain how to  compensate   for possible singularities of the solutions of the Plateau problem with free boundary at the corners. 
 
\subsection { $C^0$-Stability of the  Mean Curvature and  Reifenberg Flatness.}

The above implications   $(\ast)$  and  $(\ast \ast)$ show that strictly mean convex hypersurfaces are $C^0$-stable but,   in fact, the following more general property holds \cite {hilbert}.

   Let $Y$ be a   smooth closed cooriented hypersurface in  a smooth Riemannian
manifold $X=(X, g_0)$. This $Y$ can be seen as a $\phi$-bubble for a smooth function  $\phi(x)$ that extends the function  $ y\mapsto mn.curv_y(Y)$; moreover, if $\phi$ has a sufficiently large normal derivative on $Y$, then, as we know, the  homology class  $[Y]$ is $\phi$-trapped near $Y$.

Therefore, if a metric  $g_\varepsilon$ is sufficiently $C^0$ close to $g_0$ then a small neighborhood of $Y$ in $X$ contains a $\phi$-bubble, say $Y_\varepsilon \subset X$ for the metric
  $g_\varepsilon$. 
  
 The $g_0$-volume of this bubble is close to that of $Y$, since, clearly
 $$vol_{n-1}(  Y_\varepsilon)\to vol_{n-1}(Y)\mbox   {  for $\varepsilon  \to 0$;} $$
moreover,  the volumes of  $Y_\varepsilon$ within all   $R$-balls in $X=(X,g_0)$ are also close to the volume of the balls in $Y$,
$$|vol_{n-1}(B_y(R)\cap Y_{\varepsilon})- vol_{n-1}(B_y(R)\cap Y)|\leq \delta=
\delta(\varepsilon)\to 0\mbox {  for $\varepsilon \to 0$}$$ 
 for all $y\in Y_{\varepsilon}$ and all $R>0$.

 Since the metrics $g_\varepsilon$ are close to a fixed continuous (smooth in our case)  metric $g_0$,
 they satisfy almost Euclidean {\it filling inequalities} by Almgren's theorem {\cite {almgreniso}.
 Consequently, small balls in the bubbles  $Y_\varepsilon$ in $ (X, g_\varepsilon)$ have {\it an almost Euclidean lower bound on the volume growth. }
 
On the other hands, small balls in $Y$ have roughly Euclidean volumes and 
since the metrics $g_0$ and $g_\varepsilon$ are  mutually $\varepsilon$-close, the volumes of $B_y(\rho)\cap Y_{\varepsilon}$ are essentially the same for these metrics and the corresponding  metric balls $B_y(\rho).$

Then it  follows by all of the above   that  the $g_\varepsilon$-volumes of  $Y_\varepsilon$ within small   balls in $X=(X,g_\varepsilon)$ are {\it multiplicatively} bounded  by the volumes of the Euclidean balls
  $$vol_{n-1}(B_y(\rho)\cap Y_{\varepsilon})\leq (1+\delta')vol_{n-1}(B_{Eucl}^{n-1}(\rho)),\mbox { where }\delta'=\delta'( \varepsilon, \rho)\to 0  \mbox {  for } \varepsilon +\rho\to 0. $$

  By Allard gap/regularity theorem, \cite {allard} this bounds implies {\it regularity} of $Y_\varepsilon$ for small $\varepsilon>0$.
  
  Since this equally applies to the metrics $g_t=(1-t)g_0+tg_\varepsilon$, $0\leq t\leq 1$, one can construct a {\it diffeotopy of smooth(!) bubbles}   $Y_t \subset (X,g_t)$ between  $Y$ and $Y_\varepsilon$.\vspace {1mm}

{\it Warning.} The bubble $Y_t$ may be non-unique due to possible  bifurcations   at certain $t$ (as critical points of one parameter families of smooth functions do) and one can not guarantee smoothness of the  family $Y_t$ in $t$.

  {\it Remark.} Since $Y_\varepsilon$ comes by a deformation process of {\it smooth} bubbles, the standard elliptic estimates   suffice for the existence of $Y_\varepsilon$ and the use of the  geometric measure theory can be avoided at this point.\vspace {1mm}

  {\it Reifenberg's Flatness.}  Besides being $C^1$-smooth, the family  $Y_{\varepsilon\to 0}$ is {\it Reifenberg flat:}   \vspace {1mm}

   \vspace {1mm}
  {\it the   Hausdorff  tangent cones of  $Y_{\varepsilon\to 0}$   are isometric to $\mathbb R^{n-1} $ where the (abstract)  Hausdorff distance  of $\lambda\cdot Y_\varepsilon= (Y_\varepsilon, \lambda\cdot g_0)$ to $\mathbb R^{n-1} $ is uniformly small  at all points $y$
  in   $Y_\varepsilon$ and all $\varepsilon \leq \varepsilon_0>0$,}   \vspace {1mm}
  
 $$ dist_{Hau} (Y_\varepsilon\cap B_y(\rho), B^{n-1}_{Eucl}(\rho))\leq \delta''\rho,\mbox {where  $\delta''=\delta''(\varepsilon_0, \rho)\to 0$ for $\varepsilon_0+ \rho \to 0$.}  $$
 
 {\it Proof.} The Hausdorff limits  of  $\lambda\cdot (Y_\varepsilon\cap B_y(\rho))$ for $\lambda\to \infty$ and $\varepsilon\to 0$ are minimal hypersurfaces in $\mathbb R^n=\lim_{Hau}( X,\lambda g_\varepsilon)$ for $\lambda\to \infty$. 
 
 By the above, these hypersurfaces  satisfy the  Euclidean  bound on their volume growth; hence they are flat. QED. \vspace {1mm}
  
 {\it Remarks.} (a) The Reifenberg flatness of submanifolds   $Y_\varepsilon$ {\it does not} make them   $C^1$-close to $Y_0$ for small $\varepsilon$. In fact, the  $g_0$-normal   projections $Y_\varepsilon\to Y$ may be {\it many}-to-one  for all $\varepsilon>0$.
 \vspace {1mm}

(b)  Not all "$g_\varepsilon$-bubbles"  close to  $Y_0$ have their mean curvatures close to those of $Y_0$. For instance,  there (obviously) exists   an  arbitrarily small $C^0$-perturbation $g_\varepsilon$ of $g_0$ such that $Y_0$ becomes a {\it totally geodesic   submanifold} in $(X, g_\varepsilon)$, where moreover,   a small neighborhood $U$ of $Y$  {\it metrically splits} $(U,  g_\varepsilon)=  (Y, g_\varepsilon |_{Y})\times [-\delta_\varepsilon,\delta_\varepsilon]$.  (Of course, the above $\phi$-bubbles  $Y_\varepsilon$, albeit being close to $Y_0$ 
will not, in general, lie in such neighbourhood.) \vspace {1mm}

(c) The construction of the  "bubble perturbation"  $Y_\varepsilon$ of $Y_0$   partly generalizes to  smooth $k$-dimensional closed  submanifolds $Y_0\subset X$   for all $k<n-1$.
 
 For instance if $Y_0$ is minimal and moreover if it is an {\it isolated} local minimum for the function $Y\mapsto vol_k(Y)$, then, clearly,   a small neighbourhood of  $Y$ in $X$ contains  a $g_\varepsilon$-minimal subvariety that is non-singular and is diffeotopic to $Y_0$ as well as  Reifenberg flat.   
 
 In general, if  $Y_0$ is  not locally minimizing, one minimizes  
the function 
$$Y\mapsto vol_k(Y)-C\int_Y dist_X(y, Y_0),$$
where $vol_k$ is taken with $g_\varepsilon$ and $dist_X$ with $g_0$.

It is not hard to see that, if $g_\varepsilon$ is $C^0$-close to $g$, then this function has a local minimum
realized by a  subvariety $Y_{min}\subset X$ that is  contained  in a small neighbourhood of $Y_0$  where it is homologous to $Y_0$. 

Moreover, 
 Almgren's regularity result \cite {almgreniso}   seems to imply (as it does  happen for $k=n-1$) that such a minimal $Y_{min}$
is necessarily $C^1$-smooth;  hence, diffeotopic to $Y_0$. (I am not certain about the mean curvature of this   $Y_{min}$.) \vspace {1mm}

 (d) Reifenberg flatness of (families of  subsets) $Y_\varepsilon $ in smooth Riemannian  manifolds  $X$ implies that these are topological (actually H\"older) submanifolds by {\it Reifenberg's topological 
 disk theorem} \cite {reifenberg  }that was extended to  abstract metric space $Y$  in \cite{cheeger-colding}.
 
 An easy result in this respect  that implies the homotopy version of  the  disk theorem  for Reifenberg (sufficiently) flat  $Y_\varepsilon \subset X$   is the existence of smooth  approximation $d(x)$ to the distance function $x\mapsto dist (x,Y)$  such  that the function   $d(x) $  vanishes on $Y$ and has no critical points $x\in  X\setminus Y$  close to $Y$  \cite {hilbert}.

  \vspace {1mm}

  {\it On Piecewise Smooth  $C^0$-Mean Stability.}  What we  needs for the proof of the  $C^0$-stability of the inequality $scal(g)<0$  is a generalization of the above  to {\it piecewise } smooth $Y\subset  X$ that are boundaries of 
 mean convex polyhedral domains $U$ in $X$, where approximating $Y_\varepsilon$ must be also {\it piecewise } smooth and  have
 the mean curvatures of their   smooth pieces, that are  $(n-1)$-faces for $Y=\partial P$,    being close 
 the mean curvatures of the corresponding pieces of $Y$ while keeping
 the (dihedral) angles between these pieces  of $Y_\varepsilon$   close to the corresponding angles in $Y$.
 
We construct such $Y_\varepsilon$ by consecutively adding (hopefully)   smooth faces one by one, solving at each step  {\it Plateau's free boundary problem}  where its solvability follows from  Reifenberg kind of flatness that  implies  {\it rough filling bounds} needed for  "cutting off" undesirable "long narrow  fingers" ,  as we shall see  in the next section.

But  the solutions to such Plateau problem, say  for $Y \subset U $ with boundary $\partial Y \subset \partial U$, may have singularities  at the  corners in $\partial U$ even if all tangent cones to $Y$ at these points are smooth, where we prove this smoothness in our case by the means of {\it sharp Euclidean  filling bounds}  

  These  singularities may be due to the failure of Allard's type regularity, which seems unlikely,  or, which is more  probable, they may come  from the linearized   Plateau.
 
We shall manage to  "go around singularities" in our special   case by a rather artificial argument in section 4.9   while the general  problem of

\hspace {15mm} $C^0$-stability of {\it mean curvature} \& {\it dihedral angles} 

 \hspace {15mm} for piece-wise smooth hypersurfaces 

 \hspace {-6mm}  remains open.

  \subsection {Reifenberg Flatness and   Filling Profile.} 

{\it $(\delta, \lambda)$-Flatness.} A  hypersurface $W$ in  Riemannian manifold $X$ is called  {\it $(\delta, \lambda)$-flat  at a point $w\in W$ on the scale $\rho$}, for given positive numbers  $\delta>0, \lambda\geq 1$ and $\rho>0$, if there exists  a "flattening" $\lambda$-bi-Lipschitz homeomorphism $L$ from the ball $B=B_w(\rho)\subset X$ to the Euclidean space  $\mathbb R^n$, $n=dim (X)$, such that 

$\bullet$ $L(w)=0\in \mathbb R^n$;

$\bullet$ the image $L(W\cap B)\subset \mathbb R^n $ is contained in the $\delta \rho$-neighbourhood of the hyperplane $\mathbb R^{n-1}  \subset \mathbb R^n$.

$\bullet$ the ball $B_0^{n-1}(\lambda^{-1}\rho)\subset  \mathbb R^{n-1} $ is contained in the $\delta \rho$-neighbourhood of the image $L(W\cap B)\subset \mathbb R^n $.

  \vspace {1mm}
 
{\it Filling Volume.}  Let $U\subset X$ be a domain with boundary $W=\partial U$ and let
$C \subset U$ be a  relative $(k-1)$-cycle  with $\partial C\subset W$
 
 Denote by  $Filvol_k(C)|_U$ the infimum of the $k$-volumes of the   chains $D\subset U$ bounded by $C$, i.e. such that  $\partial D\setminus Z= C\setminus Z$ and write 
 $Filvol_k(C)|_{U\cap B}$ for this infimum taken over the chains $D$ that are  contained in the intersection of $U$ with a given ball $B=B_u(r)\subset X$.
 
{\it Remark.}  A specific local geometry of these chains and cycles is rather irrelevant for our purpose. One may think at this stage  of these being realized by piecewise smooth subvarieties in $X$.

 \vspace {1mm}

Let us show that if $\delta>0$ is sufficiently small depending on $\lambda\geq 1$ and if the boundary $W$ of $U$ is $(\delta, \lambda)$-flat  on the scales $\leq \delta_0$ for a given $\rho_0$, then {\it $U$ satisfies rough Euclidean filling inequality }
 $$Filvol_k(C)|_U  \leq const_{k,\lambda} vol_{k-1}(C)^\frac {k}{k-1},$$ 
for all $k$-cycles $C$ of diameters $\leq 0.1\rho_0$, all $k=2,3,...n-1$. 
A  precise formulation of this is as follows. \vspace {1mm}

\textbf { Rough  Filling Inequality. } There exist a continuous function $\delta_n( \lambda)>0$
and a constant $const_k\leq (10k)^{10k^2}$ with the following properties.

 Let  $\rho_0>0 $  be given and let $U\subset  X$ be a domain  with boundary $W$
 and let $C\subset U$ be a relative $k$-cycle   that is contained in the intersection of $U$ with the Riemannian ball $B_{u_0}(\rho_0/10)\subset X$, $u_0\in U$.  \vspace {1mm}

 {\it If $W$ is $(\delta, \lambda)$-flat for  some
  $\delta\leq \delta_n\cdot  (10\lambda)^{-(10 n)^{10n}}$  on all scales $\rho \leq \rho_0$ at all points $w\in W$ within distance $10\lambda \rho_0$ from $u_0$, where $\delta_n>0$ is a constant, that, in fact,  may be  assumed  $\geq (10n)^{-10n}.$
 Then  
$$Filvol_k(C)|_{U\cap B}  \leq const_k\cdot \lambda^{2k} vol_{k-1}(C)^\frac {k}{k-1}\leqno (1)$$ 
 for the ball $B=B_{u_0} (10\lambda \rho_0)\subset X$.
 
 Moreover, if $dist(c, W)\leq d$ for a given $d\geq 0$ and  all $c\in C$, then 
 $$Filvol_k(C)|_{U\cap B}  \leq const'_k\cdot \lambda^{2k}\cdot d\cdot  vol_{k-1}(C)\leqno (2)$$ }  
 for another universal constant $const'_k$. \vspace {1mm}

\vspace {1mm}

\vspace {1mm}

{\it Remarks.} (a) Filling inequalities for relative cycles  $C$ in $U$ with $\partial C\subset \partial U$ are equivalent to such inequalities for absolute cycles 
 in the {\it double} of $U$, where these cycles are symmetric under the obvious involution of the double.

(b) If $W=\partial U$ is everywhere $(\delta, \lambda)$-flat then $W$ admits a {\it collar } in $U$ that is a neighbourhood homeomorphic to $W\times [0,1]$. This   is seen with  the    smoothed distance function $u\mapsto dist(u, W)$ that has no critical points away from $W$.

In general, there is {\it  no Lipschitz collar}. On the other hand, 
 Reifenberg's disc theorem   says in the present context that  $W$  is a H\"older $(n-1)$-manifold. (Possibly, the  gradient flow of the    smoothed distance function $u\mapsto dist(u, Z=\partial U)$ may lead to a simple poof of this.)

\vspace {1mm}

{\it Proof of (1) and (2)}. We    combine the   filling  argument by induction on $k$  from \cite {filling}, \cite {wenger} with  Reifenberg's style  multi-scale  iteration process where the latter is amounts to proving the inequalities (1) and (2) on the scale $\rho$, i.e. for cycles contained in $\rho$-balls,
provided that such inequalities with slightly different constants are valid  on the scale $\delta\rho$ for some moderately   small positive  $\delta <1$.

Induction starts with $k=2$.  One  think of relative  1-cycles as 
 arcs  $C \subset U$ with the  ends in $W=\partial U$ and 
 multi-scaling is
especially simple. 

 To see this 
let $C\subset U\cap B$ where $B$ is a ball of radius $\rho\leq \rho_0/10$ and let us "chop away" the part of $C$ that is $d$-far from $W$  
for $d=2\delta_n\lambda\cdot\rho $. 

This done with the  "flattening"  homeomorphism $ L:B\to \mathbb R^n$, where may assume that it sends $B\cap U$ to the half space $\mathbb R^n_+$. 

Let  $\mathbb R^{n-1}_{-\varepsilon} \subset \mathbb R^n_+$ be the hyperplane parallel to  
$\mathbb R^{n-1}=\partial \mathbb R^n_+$ 
obtained by moving 
 $\mathbb R^{n-1}$ inside $\mathbb R^n_+$ by $\varepsilon$. We cut $C$ by the $L$-pullback 
 $H=L^{-1}(\mathbb R^{n-1}_{-\varepsilon})$ for $\varepsilon=\delta_n\rho$ and let  $C_{outH} \subset $ be the pullback,
  $$ C_{outH} =C\cap  L^{-1}((\mathbb R_+^{n})_{-\varepsilon})$$
for  $(\mathbb R_+^{n})_{-\varepsilon}\subset  \mathbb R_+$ being the half-space  bounded by $\mathbb R^{n-1}_{-\varepsilon}$.
 
 Observe that the 0-chain $H\cap C$  bound a $1$-chain $C' \subset U$, such that
 \vspace {1mm}

$\bullet$ $length(C') \leq \lambda length(C)$;

$\bullet$  $C'$ is contained in the $d$-neighbourhood of $W$ for $d\leq \lambda \varepsilon.$

$\bullet$ the sum $C^{\Rightcircle}=C'+  C_{outH} $ of $C'$ is an absolute cycle 
that bounds a $2$-chain $D^{\Rightcircle}\subset U\cap \lambda\cdot B$ with 
$$area(D^{\Rightcircle})\leq 4\lambda^3\cdot  length(C)^2,$$
where $\lambda\cdot  B$ denotes for a ball $B$ of radius $r$ denotes 
the concentric ball of radius $\lambda r$. \vspace{1mm}

\hspace{19mm}{\it Thus, the inequality (1) is reduced to (2) with 

\hspace{19mm} a controlled  "worsening" of the  constants.}

Now let us  derive (2) on the $\rho$-scale from (1) on a significantly smaller scale $3d$  as follows.

Let  an arc $C$ of length $l$  lie within distance $d$ from $W$.  
 Subdivide $C$ into $k\leq (l+d)/d$ segments of length $\leq d$ and connect the ends of these segments with nearest points in $W$ by curves of length $\leq d$.  
 
 Thus, we decomposed the chain $C$
into the sum of $k$  chains $C_i$, $i=1,2,...,k$ of length $\leq 3d$ and summing up inequalities (1) for $C_i$ we obtain (2) for $C$.\vspace {1mm}

Then we iterate this process with fillings of $C_i$ reduced to those for even smaller $C_{ij}$ etc.

Since the bound (2) in the filling volume (area for $k=2$)  is significantly stronger then (1)
for small $d$, the  infinite iteration of this process produces a finite total sum  of filling 
  areas that satisfies the inequality (1). QED. \vspace {1mm}

{\it Remark.} Even for smooth arks $C$,  the final chain $D$ filling $C$ may became "infinitely complicated" in the vicinity of $W$. But if $W$ is smooth, then piecewise smooth $C$
will be filled by piecewise smooth surfaces $D$.

In order to apply a similar argument to  higher dimensional cycles $k> 2$  we observe the following.  \vspace {1mm}

(I)   Cutting away the part of $C$ that is $d$-far from $W$ is same for all $k$ and causes no additional problem.

\vspace {1mm}

(II) Let $C \subset U$ be an {\it absolute} $(k-2)$-cycle, i.e. $\partial C=\emptyset$. If $Diam(C)\leq \rho_0$ and  $C$ lies within distance $d$ from $W=\partial U$,  then $C$ bounds a {\it relative} $(k-1)$-chain $D$ such that 
   $$vol_{k-1} (D)\leq10  \lambda^{2k-1}(d+\delta\cdot Diam (C)) vol_{k-2}(C).$$
If  $k=2$, this corresponds to moving ends of arks to $W$ where the Reifenberg's flatness of  $W$ was unneeded. v

(III) A compilation  for $k>2$ arises when we try to subdivide a chain $C$ of diameter $\rho$ that lies $d$-close to $W$ for $d<<\rho$ 
into pieces  of diameters  $\approx d$ and with $(k-1)$-volumes $\approx d^{k-1}$. This may be impossible if $C$ contains a significant thin part, where the intersections of $C$
with the balls $B_c(d)\subset X$, $c\in C$, have their volume much smaller than $d_{k-1}$. 

However, the inductive filling argument from \cite {filling} shows that the thin part can be filled in  on the  scale $\leq d$. Thus one may assume there is no thin part and $C$ can be subdivided 
into pieces of diameters about $d$, with $(k-1)$-volumes about $d^{k-1}$  and with $(k-2)$-volumes of the boundaries of these pieces about $d^{n-2}$. The latter together with (I) allows 
a decomposition $C=\sigma_iC_i$ as for $k=2$ and the validity of (1); hence of (2), disestablished for all $k$. QED.
\vspace {1mm}
 
{\it Remark.} The above argument   makes only an {\it rough outline} of a proof. But  the details are rather trivial if seen in the filling context of \cite {filling}. In fact, the above argument perfectly works for general Banach spaces instead of $\mathbb R^n$ with (properly readjusted) constants  depending on $k$ but not on $n$. \vspace {1mm}

{\it Corollary:  $C^0$-Stability of Relative Bubbles.} Let  $X =(X,g_0)$ be a smooth Riemannian manifold, $U\subset X$ be a compact domain with smooth boundary $W=\partial U$ and let $Y\subset U$ be a smooth hypersurface with boundary in $W$ that is everywhere normal to $W$.
Let $g_\varepsilon$, $\varepsilon>0$,  be a family of smooth Riemannian metrics that $C^0$-converge to $g_0$ for $\varepsilon \to 0$. \vspace {1mm}

{\it Then   there exist families of smooth domains $U_\varepsilon$ bounded by  hypersurfaces
$W_\varepsilon$ and of smooth hypersurfaces $Y_\varepsilon \subset U_\varepsilon$  with 
$\partial Y_\varepsilon \subset W_\varepsilon=\partial U_\varepsilon$ that are everywhere normal to $W_\varepsilon$ and where $W_\varepsilon$ converge to $W$ in the Hausdorff metric 
while $Y_\varepsilon$ Hausdorff converge to $Y$ for $\varepsilon \to 0$ and  such that the mean curvatures of  $W_\varepsilon$ and of  $Y_\varepsilon$ converge to the mean curvatures
of $W$ and $Y$ respectively.

 Moreover,  if $\varepsilon$ is sufficiently small, $\varepsilon\leq \varepsilon_0>0$, then the pair   $(U_\varepsilon ,Y_\varepsilon)$ can be joined with  $(U ,Y)$ be a diffeotopy in $X$ that is $C^0$-close  to the identity diffeomorphism.}  \vspace {1mm}

{\it Half-Proof}. The existence of $W_\varepsilon$ follows from the $C^0$-stability proven for individual manifolds in the previous section (where these were denoted $Y_\varepsilon$).

 Then the same argument  variation argument  delivers  $Y_\varepsilon$ with boundary
  in $W_\varepsilon$, since  $W_\varepsilon$ is Reifenberg flat and satisfies the above  filling inequality and this inequality does not allow escape of narrow fingers near $W_\varepsilon$.
  
  It is known  \cite {gruter}, \cite {nirnberg} that such a $Y$ has the same type of regularity at the boundary as at the interior points. In particular, these $Y_\varepsilon$ are smooth up to the boundary if  $n=dim X\leq 6$. 

One can not claim at this point  that the manifolds  $Y_\varepsilon$ are diffeomorphic to $Y$  for small $\varepsilon$ as it was done in the case of manifolds $Y$ without boundary, but the existence of such a diffeomorphism (and even  of a diffeotopy) is ensured by the {\it sharp filling inequality} that we prove in the next section.

\subsection { Sharp Filling under Reifenberg's Controle.}

Let us generalize the definition of \it $(\delta, \lambda)$-flatness} by replacing 
$\mathbb R^{n-1}  \subset \mathbb R^n$ by a more general the hypersurface in $\mathbb R^n$.

We limit ourselves to the case where this hypersurface serves as the boundary of a compact  convex domain $A\subset \mathbb R^n$  and define below {\it  $(\delta, \lambda)$-control of 
co-oriented hypersurfaces $W\subset X$ by $\partial A$}. Recall that  coorientation means that  we distinguished what is "locally inside" and what is "outside" $W$ and to simplify notation we assume that $W$ bounds a domain $U\subset X$. \vspace {1mm}

{\it Definition of    $(\delta, \lambda)$-Control.}  A  hypersurface $W=\partial U$ in  Riemannian manifold $X\supset U $ is called  {\it $(\delta, \lambda)$-controlled by $\partial A\subset \mathbb R^n$
 at a point $w\in W$ on the scale $\rho$},
   if there exists  a $\lambda$-bi-Lipschitz   "control" homeomorphism $L$ from the ball $B=B_w(\rho)\subset X$ to the Euclidean space  $\mathbb R^n$, $n=dim (X)$, such that 
\vspace {1mm}

$\bullet$ $L(U\cap B)\subset A$;

$\bullet$ the image $L(W\cap B)\subset \mathbb R^n $ is contained in the $\delta \rho$-neighbourhood of the boundary $\partial A\subset A$.

$\bullet$ the intersection of $\partial A$ with the Euclidean   ball $B_{L(w)}(\lambda^{-1}\rho)\subset   \mathbb R^n$ is contained in the $\delta \rho$-neighbourhood of the image $L(W\cap B)\subset A $.
\vspace {1mm}

We want to show that if $\lambda$ is close to $1$ and $\delta$ is small, then the {\it  relative $k$-filling profile} in $(U,\partial U)$ is almost the same as in $(A, partial A)$, where  such a  profile, 
 $ReFill^k_A(v)$, $ v>0$,  is defined as the infimum of the filling volumes of the relative  $(k-1)$-cycles  $C$ in $(A, \partial A)$ with $vol_{k-1}(C)\leq v$. Namely we have the following

\vspace {1mm}
{\it Sharp Filling Inequality.} Let a  hypersurface $W\partial U$ in  Riemannian manifold $X\supset U $   be   $(\delta, \lambda)$-controlled by $\partial A\subset \mathbb R^n$
 at all  point $w\in W$  and  on all scales $\rho<\rho_0$.\vspace {1mm}

{\it Then the $k$-filling profiles of  $U$ are bounded by those of $A$ as follows:
 $$ ReFill^k_U(v)\leq \alpha_A(\lambda, \delta)ReFill^k_A(v) \mbox {  for all $v\leq const_A  \lambda^{-k}\rho_0^{k-1}$ },$$
  where  $\alpha_A(\lambda, \delta)\to 1$ for $\lambda\to 1$ and   $\delta\to 0$.}

\vspace {1mm}
{\it Proof.} If a cycle $C$ is $d$-close to $W=\partial U$ for $d<<vol_{k-1}(C)^{1/(k-1)}$,
then  the rough filling inequality (2) from the previous section yields  a much stronger  filling bound  than that by  $ReFill^k_A(v)$.

On the other hand,  the part $C_{far}\subset C$ that is $d$-far from $W$
can be filled in $U$ as efficiently as in $A$. In fact,  if $\lambda$ is close to $1$, one can actually think of this  $C_{far}$
being the relative cycle in the subdomain $A_{-d}\subset A$ (that consists of the points that are $d$-far from the boundary $\partial A\subset A$) with $\partial C_{far}$ contained in the boundary of  $A_{-d}$.

The only remaining  problem is that the relative $k$-chain
$D^{\Rightcircle} \subset A_{-d}$ that fills  $\partial C_{far}$  modulo   $\partial A_{-d}$  may have large $(k-1)$-volume
of its intersection with $\partial A_{-d}$ that would make the $(k-1)$-volume of the cycle $ C^{\Rightcircle}= \partial D^{\Rightcircle}$ much greater than that of $C$,

However, the {\it coarea  inequality}, applied to $D^{\Rightcircle} $ intersected with the levels of the distance function $a\mapsto dist(a,\partial A)$, shows that some of these intersections for possibly larger but still controllably small $d' \geq d$ will have their volumes not much larger than that of $C$ unless all of $D^{\Rightcircle} $; hence, all of the $C$, lie within distance $<< vol_{k-1}(C)^{1/(k-1)}$ from the  boundary $\partial A$. This yields the proof of the sharp inequality  that completes the proof of the above $C^0$-stability of relative bubbles. \vspace {1mm}

{\it Remarks.} (a)  We formulated both,  rough and sharp filling inequalities for a specific purpose of proving the $C^0$-stability of the mean curvatures \& dihedral angles  of hypersurfaces. Probably, there is  a more 
general formulation of such an inequality and a more transparent   proof  that is  not  overburdened with trivial technicalities.

(b) If $\delta>0$ is sufficiently small $\delta\leq \delta_0 (A,\lambda)>0$, and   $W\subset X$ is  $(\delta, \lambda)$-controlled  by $\partial A\subset \mathbb R^n$, then, apparently, 
   $W$ is a  {\it H\"oder submanifold} in $X$. This   seems to follow by the argument from    \cite {reifenberg} and/or \cite{cheeger-colding}.

\subsection {Reifenberg and H\"oder at the Corners.}

Let us recall relative bubbles $Y_\varepsilon \subset  U_\varepsilon$ with  their boundaries 
$\partial Y_\varepsilon\subset W_\varepsilon =\partial U_\varepsilon$ in Riemannian
manifolds $X=(X,g_\varepsilon)$ with smooth Riemannian metrics that  $C^0$-approximate the original  metric $g_0$ on $X$. \vspace {1mm}

Let $U_{\varepsilon 1} \subset U$ be a cornered subdomain bounded in $U_{\varepsilon}$ by $Y_{\varepsilon}$ and let $A^\lefthalfcap\subset \mathbb R^n$ be the intersection of two  half spaces bounded by mutually orthogonal hyperplanes, i.e.  $A^\lefthalfcap=\mathbb R^2_+\times \mathbb R^{n-2}.$

Observe that the evaluation of the  filling profile of $A^\lefthalfcap$ reduces to that for  $\mathbb R^n \supset A^\lefthalfcap$, since every relative cycles $C\subset A^\lefthalfcap$
define an absolute cycle in $\mathbb R^n$, call $4\lefthalfcap C$ that is obtained  by reflecting $C$ 
 four times  around the   $(n-1)$-faces of $A^\lefthalfcap$.  Thus, the 
{\it extremal} relative  $(n-2)$-cycles $C$ (we care for the sharp inequality only for $k=n-1$) i.e. those maximizing  $$\frac {ReFill^k(C)^\frac {1}{k}}{vol_{k-1}(C)^\frac {1}{k-1}}$$
are intersections of  $A^\lefthalfcap$ with the round $(n-2)$-spheres that meet the two $(n-1)$-faces of 
 $A^\lefthalfcap$  at the $90^\circ$ angles. \vspace{1mm}

The limit  argument that was  used in section 4.6  for deriving the Reifenberg flatness of $Y_{\varepsilon\to 0}$  implies the following.

\vspace {1mm}

{\it There exist functions $\lambda =\lambda(\rho, \varepsilon )> 1$ and  $\delta=\delta(\rho, \varepsilon)> 0$, such that 
 $\lambda(\rho,\varepsilon)\to 1$ and   $\delta(\rho, \varepsilon)\to 0$  for $\rho\to 0$ and  $\varepsilon \to 0$  and  such that 
the  boundaries of the domains $U_{\varepsilon 1}$ are  $(\lambda, \delta)$-controlled by 
$A^\lefthalfcap$ at the   scales  $\leq \rho$ and at all  points $x\in \partial U_1$. }\vspace {1mm}.
 
It follows that the domain  $A^\lefthalfcap$ has almost the same  filling profile on the small scales $\rho$ as $A^\lefthalfcap$. Therefore, one is in a position to construct bubbles in  $U_{\varepsilon 1}$, call them  $Y_{\varepsilon 1} \subset U_{\varepsilon 1}$, with respect to the metric $g\varepsilon$  with boundaries in  
$\partial  Y_{\varepsilon 1}  \subset \partial U_{\varepsilon 1}$ with the same efficiency  and the same properties as we did it earlier   in $U_\varepsilon$, { \sc except} that we can not guarantee the smoothness of  $Y_{\varepsilon 1} $ at the points 
where $Y_{\varepsilon 1} $ meets the corner of $U_{\varepsilon 1}$.

However, all of the above seem to apply to general  to general $\mu$-bubble-hedra 
and yield 
 the  following solution of  the  {\it H\"older-regular  mean curvature $C^0$-stability problem}.
 \vspace {1mm}

Let $P$ be a compact strictly   preconvex cosimplicial polyhedral domain in a smooth Riemannian manifold $X=(X,g)$.

(Recall that "strictly   preconvex"  signifies the bound $\angle_{ij}<\pi $ for the dihedral angles of $P$; one  distinguishes the case   $\angle_{ij}=\pi/2 $
 where the  "$\mu$" in    "$\mu$-bubble" are measures  given by  continuos  density
 functions. 
 
 Also notice that an arbitrary $P$ becomes  cosimplicial, as defined in section 1   under a generic perturbation of its faces  and the general case can be probably    reduced to the cosimplicial one.)

 Let $X'$ be another Riemannian manifold with a family of  $C^2$-smooth metric $g'_\varepsilon$ and let 
 $f_\varepsilon:X' \to X$ be  an $e^\varepsilon$- {\it bi-Lipschitz} homeomorphisms. 
\vspace {1mm}

{\it Then there exist  $C^\alpha$-bi-H\"older homeomorphisms   $f'_\varepsilon : X'\to X$ for some $\alpha >0$  with the following properties.\vspace {1mm}

$\bullet_{\to0} $ The maps  $f'_\varepsilon$  converge to $f_\varepsilon$ in the $C^0$-topology for $\varepsilon \to 0$.\vspace {1mm}

$\bullet_{reg} $ the maps $f'_\varepsilon$ are $C^\infty$ diffeomorphisms  away from the corners, i.e. $(n-3)$-faces, $n=dim X$, of $P'_\varepsilon=(f_\varepsilon')^{-1}(P)\subset X'$.\vspace {1mm}

$\bullet_{stbl} $ 
$$ |mn.curv_{x'}(P'_\varepsilon)-mn.curv_x(P)| \leq \kappa(\varepsilon)\to 0 \mbox  {   for  $\varepsilon\to 0$, }$$
for  all $x'$ away from the corners   of $P'_\varepsilon$.}

(Recall that the mean curvature at the edges, i.e.  $(n-2)$-faces, $n=dim X$, is defined in terms of the dihedral angles as $\pi -\angle_{ij}$.) \vspace {1mm}

We do not go into the detailed proof since this result is neither  general enough to elucidate  the geometric meaning of $scal\geq 0$ 
nor is 
the H\"older regularity is sufficient for our applications to  positive scalar curvature. 

On the other hand, this  insufficient regularity can be bypassed in the  $C^0$-non-approximation  application as we shall see below.

\subsection {Proof of the $C^0$-Limit Theorem.}

We shall show in this section  that \vspace {1mm}

{\it  Smooth metrics  $g$ of negative scalar curvature  can not be $C^0$-approximated by metrics $g_\varepsilon$  with nonnegative scalar curvatures.} \vspace {1mm}

{\it Proof.} Let $g$ me a smooth  metric on an $n$-manifold $X$ with negative scalar curvature 
at a point $x_0\in X$. Then \vspace {1mm}

 ($\smallsquare $)   {\it all sufficiently small  neighbourhoods of $x_0$ contain  (tiny) mean curvature convex  cubical polyhedral domains $P\ni x_0$ with strictly  acute  (i.e. $<\pi/2$) dihedral angles. }\vspace {1mm}

{\it Proof of}  ($\smallsquare $) . We assume by induction that such domains exist in submanifolds $X'\subset X$ that contain $x_0$ that that have zero second fundamental form at $x_0$.

Clearly, there exists a codimension 1 submanifold in  $ X'\subset X$ that  contains $x_0$, such that

$(\bullet $)  submanifold  $X'$ is totally geodesic at $x_0$ i.e. its  second fundamental form vanishes at $x_0$; 

 $(\bullet $) the induced metric in $X'$  has  strictly negative  scalar  curvature;

  $(\bullet $)  the Ricci curvature  of $X$ at the normal   vector  $\nu_0$ to $X'$ at $x_0$  is strictly   negative.

  It is also clear that such an $X'$ admits an arbitrary small perturbation, call it $X'' \subset X$ near $x_0$ that still contains $x_0$ and such that  its mean curvature becomes zero near $x_0$ while keeping the second fundamental form  II zero at $x_0$. (Small  non-zero  II will also do.)

We may assume by induction on $dim(X)$ that  $X''$ contains required polyhedral  domains,
say $P''\subset X'$. Take such a $P''\ni x_0$ in $X''$ of very small diameter $\delta$ and let 
$P=P''\times [-\varepsilon,+\varepsilon]$ be the union of the  geodesic 
$2\varepsilon$-segments normal  to $P''$ and going by $\varepsilon$ in both normal directions at all points in $P''\ni x_0$.

If $\varepsilon$ is sufficiently small  that the "horizontal" faces $P''\times \pm \varepsilon \subset X$ have positive mean curvatures by the second variation formula. It is also clear that of the dihedral angles between these horizontal faces and the remaining "vertical" ones equal   $\pi/2$ while the angles between vertical angles are $<\pi/2$ for small $\varepsilon$. Also a simple computation shows that the mean curvatures of the vertical faces are positive.

Thus $P$ satisfies all requirements except for having some dihedral angles $=\pi/2$ but these
can be made acute by an arbitrary small perturbations of the two "horizontal" faces. \vspace {1mm}

If we could prove  the $C^0$-stability of the mean convexity for this $P$ with  $ P'_\varepsilon=(f'_\varepsilon)^{-1}(P)$  being $C^2$-smooth (rather than mere  
H\"older as in the previous section) everywhere including the corners, then   we  would  apply  "gluing around the edges"  to $ P'_\varepsilon$ as in section 2.1  and arrive at a  metric with positive scalar curvatures on the $n$-torus. 
But as we have not proved this stability, we need to combine the mean curvature stability for {\it individual hypersurfaces} (i.e. polyhedral domains of depth $d=1$) with the gluing   as follows.\vspace {1mm} 

Let $X=(X,g)$ be a smooth Riemannian manifold with compact strictly mean curvature convex boundary $Y=\partial X$. If  $scal(g)>0$ then  (see \cite {gromov-lawson1}, \cite {almeida})
the  double $2
\underset {Y}\diamond X$  admits a family of metrics, say  $2
\diamond g_\delta$, $\delta>0$, where  $g_\delta$ are metrics on $X$ such that: \vspace {1mm}

$\bullet_{\leq \delta}$ the metrics   $g_\delta$ are $\delta$-close to  $g$ in the $C^0$-topology; \vspace {1mm}

 $\bullet_{C^2}$ the double metrics $2
\diamond g_\delta$, the restriction of which on both copies of $X\subset 2
\underset {Y}\diamond X$ by definition  equal $g_\delta$, are $C^2$-smooth; \vspace {1mm}
 
 (Notice, that the $C^1$-smoothness of $2
\diamond g_\delta$ is equivalent to $Y$ being totally geodesic with respect to 
$g_\delta$, and  $C^2$-says something about the curvature tensor of  $g_\delta$
on $Y=\partial X$.) \vspace {1mm}

$\bullet_{sc>0}$ the metrics $g_\delta$, and hence  $2
\diamond g_\delta$,  have positive scalar curvatures.

  \vspace {3mm}
  
 Let  $P$ be an $n$-dimensional  rectangular reflection domain (see section 2), e.g. a cubical one, 
 and let $g$ be a smooth Riemannian metric on $P$ with respect to which the faces of $P$ have  strictly positive mean curvatures and the dihedral angles are all $\pi/2$.
 
 Let $g_\varepsilon$ be a family of smooth Riemannian metrics on $P$ such that:\vspace {1mm}
 
 $\bullet_{\leq \varepsilon}$ the metric  $g_\varepsilon$ is $\varepsilon$-close to $g$, in the $C^0$-topology for all $\varepsilon <0$;\vspace {1mm}
 
 $\bullet_{sc\geq 0}$  $ scal (g_\varepsilon)\geq 0$  all $\varepsilon <0$. 
 \vspace {1mm}

\textbf {Approximation/Reflection Lemma.} {\it There exists metrics $g_{\varepsilon, \delta}$ on $P$ for all  $\varepsilon, \delta>0$, and polyhedral subdomains $P'=P'_{\varepsilon, \delta}\subset P$   such that  \vspace {1mm}

 $\bullet_\sim$ the domains $P'$ are combinatorial equivalent to $P$, where such an equivalence is established by homeomorphisms $P\to P'$  that $C^0$-converge to the identity map for $\varepsilon, \delta\to 0$; \vspace {1mm}

 $\bullet_{\varepsilon+\delta}$. The metric  $g_{\varepsilon, \delta}$ is $\delta$-close to $g_\varepsilon$; thus, it is   $(\varepsilon+\delta)$-close to  $g$,  for all $\varepsilon, \delta>0$;  \vspace {1mm}

$\bullet_{sc\geq 0}$  $ scal (g_{\varepsilon, \delta})\geq 0$ for all $\varepsilon, \delta >0$;  \vspace {1mm}

 $\bullet_{reg}$ all  dihedral angles of $P'$ with respect to $g_{\varepsilon, \delta}$ equal $\pi/2$ and all 
faces of $P'$ are totally geodesic; moreover the canonical extension of $  g_{\varepsilon, \delta}$ to the metrics $\tilde g_{\varepsilon, \delta}$ on the manifold  $\tilde P'\supset P'$  where the corresponding reflection group $\Gamma$ acts
are $C^2$-smooth  for all $\varepsilon, \delta >0$.}

  \vspace {1mm}
  
  {\it Proof.}  Proceed by modifying faces of $P$ one by one and simultaneously  changing the metric.  
  
  Namely, assume by induction that  $\bullet_{reg}$ is satisfied by some faces, call them $W'_{reg}$  and  dihedral angles between them for some metric $g_{reg}$.
  Then "regularize" an extra face, say $W_i$  by
 solving the corresponding  Plateau $\mu$-bubble problem (see the  previous section) where the boundary of the new face $W'_i$ --  solution of this Plateau problem,  is contained in the union of the faces $W_{reg}$ and of all not yet regularized faces $W_j$  except for  $W_i$ itself.
 
  Observe that

{\it  the  regularity of the extremal $W'_i$ at the points in the union of $W_{reg} $ follows 

from the interior regularity by the standard and obvious reflection argument. }
 
 Finally,  replace  $g_{reg} $ by  the above  $g_{reg, \delta} $ that  makes $W'_i$ regular as well; thus,  the inductive step is  accomplished.

 An essential point here is that the condition $\bullet_{\leq \delta}$  implies that the metric $g_{reg, \delta} $ may be assumed arbitrarily $C^0$ close to $g_{reg}$ and so the relevant   filling profiles  essentially do not change as we pass  from  $g_{reg}$  to $g_{reg, \delta} $; thus, the process of consecutive "regularization" of faces and metrics goes unobstructed. \vspace {1mm}
 
{ \it Warning.}      The metrics $g_{\varepsilon, \delta}$   are $C^0$-close to $g$ but, in general, this closeness, unlike that between $g_{reg}$ and   $g_{reg,\delta}$,   "does not respect" the boundary   $Y=\cup_iW_i$ of $P$: there is {\it no Lipschitz control} over the  homeomorphism that
moves $P$ to $P'$.

  \vspace {1mm}
 
 {\it Conclusion of the Proof of the  $C^0$-Limit Theorem.}
 The  cubical    domain  from  the above  ($\smallsquare $)     can be trivially made strictly mean convex with
 all dihedral angles $\pi/2$ and, by the Lemma, the   the solution of the Geroch conjecture for tori applies.
 
 This settles the case $\kappa=0$  that extends to $\kappa\neq 0$ by (a)
 of  1.8.

\vspace {1mm}

    \vspace {1mm}

\section { Conjectures and Problems.}

 \subsection { On  Topography of Plateau Wells.} 
 
Existence/non-existence of positive scalar curvature on a manifold $X$  is invariant under the codimension $2$-surgery of $X$, i.e. adding $m$-handles with $m\leq n-2$; accordingly, one wishes to have a counterpart of the $C^0$-non-approximation  property in the category of manifolds taken modulo such surgery whatever this means. 
 
  This agrees with the observation 
that  attaching 
"thin $m$-handles"  to $X$ with $ m\leq n-2$   does not significantly change the topography of
 Plateau traps  in $X$ that can be seen as    "wells" in the $(n-1)$-volume   landscape in the space of hypersurfaces in $X$. 
 
 A representative example of such "insignificant" change  is as follows.

\vspace {1mm}

Let $Y\subset X$ be a  locally $vol_{n-1}$-minimizing hypersurface which is trapped in a small 
neighbourhood $U \supset Y$. 

Let $C$ be  a smooth  curve joining two points
 $x_1,x_2 \in X\setminus U$ which are
positioned "relatively far from" $U$ and  such that  $C$ itself as well as  all "moderately large"  perturbations"  $C'$ of $C$
intersect $Y$. 

Take an  $\varepsilon$-thin normal neighbourhood 
$D\supset C$ and modify the original metric $g$ of $X$
on $D$ by enlarging the lengths  of the geodesic segments in $D$ normal to $C$  by the factor    $\varepsilon^{-1-\alpha}$, for a small $\alpha>0$, e.g. $\alpha=0.1$.

 If $n=dim (X) \geq 3$, then the resulting  enlarged metric $g' \geq g $ on $X$ is,
   "on the average",
 $\varepsilon^\beta$-close
 to $g$ for $\beta >0$
and, in many cases, there is a minimal $Y'\subset X'$ corresponding to $Y$, 
 which, for $n-1=dim(Y)\geq 2$, is "on the average"
$\approx \varepsilon^\gamma$-close to $Y$ for a $\gamma>0$.

 But since 
such a $Y'$ must intersect $C$, it is obtained from (a slightly perturbed) $Y$ by 
attaching  several thin and {\it $\approx\varepsilon^{-\alpha}$-long} "fingers"
 corresponding to (sufficiently stable) intersection 
points of $C$ with $Y$. Since length$_{fing} \to \infty$ for $\varepsilon \to 0$ 
one can not avoid "on the average".

\vspace {1mm}

Similarly, one can  compare traps in {\it non}-equidimensional 
manifolds.  The simplest examples are Riemannian product manifolds, where, e.g. 
$X' =X\times \mathbb R^k$
  has essentially the same minimal subvarieties and Plateau wells as $X$.
 \vspace {1mm}

Besides looking at what happens  to $X'$ in the "immediate neighbourhood" of an $X$ 
we  want to keep track of  traps/wells in manifolds $X'$ that are  only
 moderately close to $X$  and where only relatively
 deep/wide  Plateau wells in $X$ have a chance of being shadowed by wells in $X'$.

\vspace {1mm}

\subsection {Webs and Honeycombs.}

A {\it Plateau  $m$-web}  $\cal M$ in  a Riemannian $n$-manifold   $X$ 
is an $m$-tuple of
 foliations  ${\cal M}_i$, $ i=1,2,...m$  by minimal subvarieties 
of codimension $1$ such that
 the subvarieties from different foliations are mutually transversal and make constant
(dihedral)  angles, say $\angle_{ij}$, $ i,j=1,2,...m$. 

\vspace {1mm}

{\it  Local Web Conjecture}.  If $m>n$ then such a web is locally isometric 
to the flat one, i.e. to $\mathbb R^n$ with $m$
families of parallel hyperplanes. Furthermore, if $m=n$ and the web is {\it normal}, i.e. $\angle_{ij}=\pi/2$
for all $j\neq i=1,2,...n$, then it is also flat, provided $X$ has $Sc\geq 0$.

\vspace {1mm} Let us look closer at  the normal webs.

 First,  every transversal  $n$-web 
 locally equals   a {\it coordinate web}:  there are
local coordinates $x_1,..., x_n$ in $X$, where  $\cal M$  identifies with
 the $n$ families of the coordinate
 hypersurfaces $x_i=cost$. 
 
  The normality condition signifies that the 
  $g_{ij}$-terms of the Riemannian metric vanish for $i\neq j$, while the minimality implies that the
 products $ G_i=\prod_{j\neq i} g_{jj}$ are invariant under the flows by the  (coordinate) 
 vector fields  $\partial_i=\partial/dx_i$. 
 
It follows, that each $G_i$ is, in fact, a (positive)  function in $n-1$  (rather than $n$)  variables, namely
in $x_j$ for $j\neq i$, and every $n$-tuple of such functions defines a normal Plateau web, where
 $g_{ii}$  are uniquely determined by the equations $\prod_{j\neq i} g_{jj}= G_i$.

It is seems not hard to show by a direct computation that the scalar curvature of such metric
$\sum_i g_{ii}dx_i^2$ is  {\it strictly} negative unless the web is flat, but I did not hat this. 
On the other hand   the inequality $scal\leq 0$ follows, as we know, from the Geroch conjecture for tori.

Apparently, the above representation of metrics generalizes  to non-normal Plateau $n$-webs
 (with constant but not normal  $\angle_{ij}$); this would imply the  local  conjecture for $(m>n)$-webs.

\vspace  {3mm}

Let us generalize the concept of normal web as follows.

 A cubical polyhedron $P_0$ is called {\it a  normal  Plateau honeycomb} if  a Riemannian $n$-manifold $X$ is a collection $\cal P$
of cubical domain $P\subset X$ such that $P_0\in \cal P$ and  \hspace {1mm}

 $\bullet $  the faces of all  $P\in \cal P$ are  minimal (possibly singular?)  hypersurfaces in $X$ where  all dihedral angles between the  faces equal $\pi/2$;\hspace {1mm}

 $\bullet $  there exists a face preserving continuous map of every $P\in \cal P$ onto the $n$-cube $[0,1]^n$ 
 such that the pullbacks of all $(n-1)$-subcubes in it parallel to the faces, $[0,1]\times [0,1]\times ... \times t\times ... \times [0,1]$, $t\in [0,1]$ are minimal hupersurfaces, call them $Q\subset P$  normal to the boundary of $P$;

 $\bullet $ the two parts into which  such a $Q$ divides $P$ are cubical domains that are
 elements of $\cal P$. 

\hspace {1mm}

{\it Motivation.}   Let  $P_0$ be a non strictly  mean curvature convex  cubical polyhedral domain,  say $P_+ \subset P$,
with all dihedral angles $\leq \pi/2$. Then either $P_0$ contains a   {\it strictly}  mean curvature convex  cubical polyhedral domain
with all dihedral angles $\leq \pi/2$ or 
 $P_0$ is  a    normal  Plateau honeycomb. \vspace {1mm}

{\it Justification.} If there is a regular point in a face where the mean curvature is strictly
positive, or if some angle $\angle_{ij}$ is somewhere $,\pi/2$, then $ P'$ can be 
obtained by  smoothing of an arbitrary small perturbation of $P$ which is achieved by
 an elementary linear(ization) argument.

Thus, we may assume that $P$ is normal Plateau.
If one of faces, say $Y$ of $P$ is not locally minimizing, one can cut $P$ by a minimizing
face and thus one may assume that every pair of opposite faces in $P$ has at least one of
 them, say $Y$, 
 being locally minimizing. 
 
 If such a $Y \subset P$ is not isolated, it serves as a leaf of a Plateau foliation and if assume
 non-existence of a Plateau honeycomb, we conclude that some of these faces, 
 let it be $Y$, is isolated. 

Then, for a sufficiently small $\varepsilon$, there exists an $\varepsilon$ bubble
$Y_\varepsilon \subset P$ and then the "band" $P'$ between  this  $Y_\varepsilon$
and $Y$ in $P$ does the job (after a small perturbation and smoothing making all faces
 of  $ P' $ strictly mean convex).\vspace {1mm}
 
{\it Remark.}  A technical difficulty in this argument resides in the regularity 
of our bubbles, especially  at the corners.

{\it Questions.} Do all  normal  Plateau honeycombs  are isometric to Euclidean solids?

Does it help to assume that $scal(P)=0$?

How much does  a presence of  singularities in  minimal hypersurfaces $Q$  
complicate the geometry of $P$?

\subsection { Nested Cubes and  Small Diameter Conjecture.}

Let $P\subset X$ be a
 normal (i.e. with mutually normal faces) mean curvature convex cubical  polyhedral domain 
  of depth $d=n=dim(X)$ in a Riemannian manifold $X$.
 
 Does $P$ contains a normal $n$-cubical subdomain  $P_\smallsquare \subset P$ of an arbitrary small $diam_X(P_{\smallsquare})\leq \delta$ for a given $\delta>0$.  \vspace {1mm}

Notice that  according to  {\tiny  ($\smallsquare $)}  from  section 4.9  a presence of a point in $P$ where scalar curvature $<0$ implies the  existence of $P_\smallsquare$  and the solution of the $P_\smallsquare$-problem  follows in some cases from the solution of the Geroch conjecture, e.g. where the faces of $P$ are smooth with {\it strictly} positive mean curvatures.

But  it is instructive to construct $P_\smallsquare $ by a  direct  argument keeping an eye on singular spaces, where the simplest case is that of a manifold $X$  with a $C^1$-smooth metric where
curvature is defined only in distribution sense. Moreover, one can formulate this problem even for $C^0$-metrics  in terms of  $\Delta$-stable mean convexity (see section 4.4) where no curvature exists even in a weak sense. 

More specifically, say that a cubical  subdomain $P' \subset P$ is a {\it sandwich} in $P$ if, 
combinatorially,   $P'$ looks as a "rectangular slice" of a cube, namely as  $ [a,b]\times [0,1]^{n-1} \subset  [0,1]^{n} $ for $0\leq a<b\leq 1$.

In other words,  $  P'$  is bounded in $P$ by  a pair of its mutually disjoint "new" faces that are  hypersurfaces, say  $Q'_($ and $Q'_)$ in $P$, both  separating  a pair of opposite faces, say $Q_($ and $Q_)$ of $P$, where "separating" here also signifies "being homologous to" in the the obvious sense. Thus, 
 the boundaries  $\partial Q'$ and $\partial Q'_-$ are contained in the boundary $\partial P$ where they intersect all  faces,  but  $Q$ and $Q_-$ of $P$, unless one (or both)  of these new faces equals the "old"
$Q$ or (and) $Q_-$.

Consider decreasing  sequences   $P=P_0\supset P_1\supset P_2\supset ...\supset P_k...$  of normal cubical domains  where
  $P_k$
is a sandwich in $ P_{k-1}$ for all $k=1,2,...$.

Call a closed subset $P_\infty\subset P$ a  {\it normal micro-cube} if it equals the intersection of $P_k$ in  such a sequence.\vspace {1mm}

{\it Example.} If $P$ equals the ordinary cube $[0,1]^n$ then (connected) micro-cubes are exactly Cartesian products of $n$ subintervals in  $[0,1]$, possibly, some reduced to points.

\vspace {1mm}

{\it Question.} What are geometries and topologies of these micro-cubes? \vspace {1mm}

{\it Conjecture.} Every normal cubical domain $P$ of depth $d=n=dim X$ contains a {\it zero dimensional} micro-cube in it.\vspace {2mm}

If $Y$ and, hence $Y_\varepsilon$ are regular, then such band can be regarded as an "infinitesimally 
thickened" $(n-1)$ dimensional cubical polyhedron and, after $n$-steps, one arrives at a 
zero dimensional 
 "cube" $P_{\smallsquare} \subset P$. 
 
 But this dimension reduction process does not a priori apply to  quasiconical singularities  
 where the dimension reduction is not apparent. (We temporarily disregard singularities at the corners at this point.) 
 
 Thus all we can claim is that  some $P_{\smallsquare}$ is no greater than the singularity;
  thus, $dim(P_{\smallsquare})\leq n-8$.

Then, if $dim X\leq 9$, some cubical domain $P_k$  approximating $P_{\smallsquare}$
is {\it spin}, since the the {\it $2$-dimensional} Stiefel-Whitney class $w_2\in H^2(P;
\mathbb Z_2)$  vanishes on $P_{\smallsquare}$ that have dimension $<2$ and  the Dirac operator method applies.

\vspace {1mm}

{\it  Corners do not Matter.}  Indeed,we  can arrange   $P_k$ such that $P_{k+n}$ is contained  in the {\it interior} of $P_k$ for all $k$. Thus, any nuisance at the boundary will be eventually forgotten. Besides we  can apply arbitrarily small $C^0$-perturbations to the Riemannian metrics  $g_k$ in $P_k$, that  smoothes the natural extensions  $\tilde g_k$ of  of $g_k$ to the orbicoverings (reflection developments) $\tilde P_k$ of $P_k$ as in section 2.

\vspace {1mm}

{\it Further Remarks and  Questions.} (a) There is  counterpart of the small  diameter problem that
appeals to  "billiard minimal"  hypersurfaces that do not have to be "parallel" to the original faces of $P$ and that can be formulated in terms of the $\mathbb T^n$-essential manifold
$X  =\tilde P/\mathbb Z^n$ associated to $P$ as follows.

Take a minimal $(n-1)$-cycle $Y(1) \subset X(1) =X$
 which may be different from those corresponding to faces of $P$. Then take an infinitesimally  narrow "bubble-band"
  $X'(1) \subset X(1)$ around $Y(1)$ and that is bounded by two smooth hypersurfaces with positive mean curvatures as we did in section 4.3.  Let   $X(2) =2X'(1)$ be the double  of  $X'(1)$.
  
  Continue similarly with some 
  $Y(2) \subset X(2)$ thus obtaining  $X(3)$, etc. 
  
   $X(1), X(2), X(3),...$.
  
  Do the  diameters of $X(i)$ (as well as of the Plateau-hedra  obtained by "cutting these $X(i)$ open")  
  converge to zero for $i\to \infty$ for suitable $X(i)$?\vspace {1mm}

(b) Can one 
move several faces  of a $P$ inward {\it simultaneously}, e.g.  near a corner (vertex) of $P$ by solving the corresponding 
 linearized  equations  and applying an implicit function theorem?

\vspace {1mm}

(c) Let  $X$ be an dimensional Riemannian manifold and  $Z \subset X$  be a compact piecewise smooth 
 submanifold   of dimension $\leq n-2$, e.g a curve in $\mathbb R^3$.

Suppose that $Z$  admits an  a mean curvature convex normal cubical approximation i.e. it  equals the
 intersection of a decreasing family of normal  $n$-cubical  m.c.c. polyhedra
$P(k) \subset X$, $k=1,2,...$.

Is then $Z$ necessarily smooth away from its boundary? Is it, moreover, totally geodesic?

We are especially concerned with the possibility of a "bad" 
(but, potentially most interesting) approximation,
where { \it all  faces} $W(k)=W^{n-1}(k)$ of $P(i)$ (there are $2n$ of these faces) and, consequently, the  $n-2$-faces,    are eventually
dense in $Z$, i.e.
 $dist_{Hau}(Y(i), Z)\to 0$,  for $i \to \infty$ and    for every sequence of $(n-1)$-faces $Y(i)$ of $P(i)$.

\vspace{1mm}

Does a $Z$ which admits a "bad" approximation by normal  cubical m.c.c. (or Plateau) $P$ 
necessarily consist of a single point? 

 This seems easy for $n\leq 7$  but a similar question   for non-cubical $P$
  appears non-trivial
even for mean curvature convex    $ P\subset \mathbb R^3$.

(The worst
 scenario  from the perspective of the $scal>0$-problems is when not only the families of faces $W(i)$
of $X(i)$ are eventually dense in $Z$, but every family of stable minimal "billiard subvarieties" 
$Y'(i) \subset X(i)$  is also eventually dense in $Z$.

\vspace {1mm}

(d) Possible quasi-conical singularities of minimal varieties is the apparent source of our problems, but, eventually,
 singularities should serve in our favor: they significantly constrain the shape of minimal varieties.

 For example, let $X$ be a compact Riemannian manifold. Then, probably, 
  every  "small" cubical
 Plateau-hedron  $ P\subset X$ with
$\angle_{ij}\leq \pi/2$ and 
 with locally  minimizing $(n-1)$-faces has all these faces non-singular, say in the interior points, 
where "small" means $diam (P) \leq const_X$.

\vspace {1mm}

 (e)  Observe the following converse to the existence of mean convex neighbourhoods of  isolated minimizing hypersurfaces.    

Let  a connected Riemannian $n$-manifold $X$ be 
divided  by  be a compact connected  subset $Y$ in into two connected domains $X_+$ and $X_-$ 
 with common boundary $\partial X_+=\partial X_-=Y.$ \vspace {1mm}

Let  $Y_0$ equals the intersection of compact domains $W\supset Y_0$ in $X$ with
  smooth mean curvature convex
 boundaries in $X$. \vspace {1mm}

 {\it If  $Y$ has finite $(n-1)$-dimensional Hausdorff measure,
 $vol_{n-1}(Y)<\infty$,  then $Y$ is  locally (in the space of $\mathbb Z$-currents)  $vol_{n-1}$-minimizing. In particular, $Y$ 
 is smooth of dimension $n-1$ away
 from a compact subset $\Sigma \subset Y_0$ with  $dim_{Hau}(\Sigma)\leq n-8$.} \vspace {1mm}
 
 (I am uncertain  on what happens if $vol_{n-1}(Y)=\infty$ that
  pertains to the question of
   possible topologies of Hausdorff limits of stable minimal hypersurfaces with volumes $\to \infty$  in  $X$.)
 \vspace {1mm}

(f)
Let  $Z=\cap_kP(k)$ for a sequence of  m.c.c. polyhedra $P(k)$ of  same combinatorial types and  with  given bounds in their dihedral angles, where one distinguishes the case  where all $P(k)$ are Plateau-hedra.

Consider ("informative") sequences of points
$x_i\subset X$ and sequences of numbers $\lambda_i\to \infty$, and 
take the (sub)limits of $ \lambda_i P(i) \subset \lambda_i(X,x_i)$.

How much of geometry of $Z$, e.g. in the case of "bad approximation", can be extracted 
from the resulting Euclidean picture(s)?

\subsection{Gauss Bonnet  Prism Inequalities and  the Extremal  Model Problem.}

Let $P$ be a mean curvature convex  $3$-dimensional Riemannian manifold with corners that is combinatorially equivalent to a prism, that is a product of a $k$-gon by  a line segment.

Let the dihedral angles at the top and at the bottom of $P$ be $\leq \pi/2$,
call such prisms {\it normal}, and let
 the dihedral angles between the  remaining faces (sides) of $P$ be bounded by some numbers $\alpha_1,...,\alpha_i,...,\alpha_k$.

\textbf {$3D$ Gauss Bonnet  Prism Inequality. }  {\it  If $P$ has non-negative scalar curvature, then the numbers  $\alpha_i$
are bounded from below by 
  $$\sum_{i=1,...,k} (\pi-\alpha_i)\leq 2\pi.$$}

{\it  Proof.} Let $Y_{min} \subset P$ be an area minimizing surface separating the top of $P$ from the bottom, that, observe, is  {\it normal} to $\partial P$. 

Temporally assume  that $Y_{min}$ is $C^2$-smooth including the corner points where
$Y_{min}$ meets the "vertical" edges of $P$ and recall that  the second variation of $area(Y_{min})$ is 
  $$area''=-\int_{Y_{min}}Ricci_P(\nu)dy -\int_{\partial Y_{min}}curv_\nu(\partial P)ds$$
for $s$ being the length parameter in $\partial Y_{min}$.

Observe 
following \cite {schoen-yau1} (also compare \cite {burago-toponogov})
that
$$Ricci_P(\nu)=\frac{1}{2}(scal(P) +\lambda_1^2+\lambda_2^2)-K(Y_{min})$$
for $\lambda_i$ denoting the principal curvatures of $Y_{min}$
and $K$ being the sectional curvature of $Y_{min}$, while
$$curv_\nu(\partial P)= mn.curv(\partial P)-curv(\partial Y_{min}).$$

Thus, we conclude as Schoen and Yau do  in \cite {schoen-yau1} that

$$0\leq area'' \leq \int_{Y_{min}}K(Y_{min})dy+\int_{\partial Y_{min}}curv(\partial Y_{min})ds,$$
where

 $$\int_{Y_{min}}K(Y)dy+\int_{\partial Y_{min}}curv(\partial Y_{min})ds +\sum_{i} (\pi-\alpha_i)=2\pi \chi(Y_{min})$$
 by the Gauss Bonnet formula. Hence,
 $$\sum_{i} (\pi-\alpha_i)\leq 2\pi \chi(Y_{min})\leq 2\pi.$$

{\it On Regularity at the Corners.} The  above $C^2$-smoothness assumption may be violated at the corners of $Y_{min}$
but it is easy to see by looking at how $Y_{min}$ is approximated at the vertices by the {\it necessarily unique}  tangent cones, that $Y_{min}$ admits tangent planes (or  rather cones) at the corners. 

Also one sees  in this limit cone picture  that the  $length(\partial Y_{min})<\infty$; hence $\int_{\partial Y_{min}}|curv(\partial Y_{min})|ds<\infty$ and 
since by minimality of  $Y_{min}$  the curvature $K(V)$ is bounded from  above,
 $\int_{Y_{min}}|K(Y_{min})|dy<\infty$. This is sufficient to justify the above computation.
\vspace {1mm}

{\it Extremality and Rigidity.}  The above argument shows that convex Euclidean prisms are extremal for the prism inequality, moreover they are dihedrally rigid: 

{\it If a  normal mean curvature convex prism $P$ with $scal(P)\geq 0$  satisfies $\sum_{i} (\pi-\alpha_i)= 2\pi$ then it is isometric to a convex  Euclidean prism.}\vspace {2mm}

{\it  Non-zero Bounds on  Scalar Curvature}. If $scal(P)\geq 2\kappa$ then the inequality $$0\leq area'' \leq \int_{Y_{min}}K(Y_{min})dy+\int_{\partial Y_{min}}curv(\partial Y_{min})ds,$$
becomes   
$$0\leq area'' \leq \int_{Y_{min}}K(Y_{min})dy+\int_{\partial Y_{min}}curv(\partial Y_{min})ds- \kappa \cdot area(Y_{min})$$
and 
 $$\sum_{i} (\pi-\alpha_i)+\kappa \cdot area(Y_{min})\leq 2\pi\chi(Y_{min}).$$
Notice that this is sharp (i.e. turns into an equality) for $P$ being the product of a $k$-gonal surface $Q^2_\kappa$ of constant curvature $\kappa$ by a line segment.

Also observe that the area of $Y_{min}$ is bounded from below by an $A$, (that is of use for $\kappa>0$)  if $P$ admits a $1$-Lipschitz map $f$ onto a disk $D$
of area $A$, such that $f$ sends the side-boundary of $P$ on $\partial D$
with the top and the bottom of $P$   being sent to $D$ with degree 1.
 
On the other hand,  $area(Y_{min})\leq A$ (that may be used for $\kappa<0$) if
the distance $d$ between the top and bottom in $P$ is related to the volume of $P$ by 
 $vol(P)/d\leq A$. 
 
 \vspace {1mm}

\textbf {Semi-Integral Inequality.} If $\kappa>0$, then the inequality  $scal (P)\geq \kappa$  can be significantly relaxed by requiring that the  {\it integral $\int_{Y}scal (P)dy$ } is bounded from below by $2\kappa\cdot area(Q^2_\kappa)$ for the above  $k$-gon $Q^2_\kappa$ and  {\it all} surfaces $Y\subset P$ separating the top from the bottom.  

The resulting inequality  (compare \cite {petrunin} is most informative if $scal(P)\geq 0$ and it can be also meaningfully used if $scal(P)\geq \kappa_-$for some $\kappa_-<0$
if one also has an upper bound on $area(Y_{min})$. For example, this shows the following.
\vspace {1mm}
 
  {\it Given an  Riemannian metric on $P$, it can not  modified  with an uncontrolled  enlargement of  its scalar  curvature along $N$ line segments joining the top to the bottom of $P$, e.g. by inserting $N$ tubes  of fixed thickness, say  $\delta$ independent of $N$, with scalar curvatures  $\geq \varepsilon = \frac {100}{N\delta^2}$ and large $N$ without generating a proportional amount of negative scalar curvature.} \vspace {1mm}

 {\it On Non-zero bounds on the  Mean Curvatures of the Faces of $P$.} Besides allowing $\kappa\neq 0$  one may similarly allow non-zero lower
bounds on the  mean curvature by some numbers  of the side faces of $P$ by some $m_i$, $i=1,2,...,k$.

If  $m_i>0$, this may be used   together with a {\it lower bound} on the "widths" of these faces, but it is less clear what kind of  {\it upper bound} on the size of these faces (and/or on all of $P$) may serve along with some $m_i<0$. \vspace {1mm}

\textbf { $K$-Area and Semi-integral Inequalities for  Closed and for Cornered  $n$-Manifolds}.  Let $X$ be a closed oriented  $n$-dimensional Riemannian manifold,
where the fundamental cohomology class $[X]^n\in H^n(X;\mathbb Z)$ equals the $\smile$-product of $2$-dimensional classes and a  class $h^k\in H^k(X;\mathbb Z)$ coming from $H^k(\Gamma;\mathbb Z)$ under the classifying map $X\to K(\Gamma;1)$.

An instance of such an $X$ is a Cartesian product of complex projective spaces, e.g. of $2$-spheres,  and  of a closed $k$-manifold  $Z$ of non-positive sectional curvature, e.g. $Z=\mathbb T^k$.

Denote by  $sc.ar(X)$  the infimum of the numbers $A$, such that the $N$-th multiples of all non-zero integer  $2$-dimensional homology classes in $X$ are representable by  surfaces  $Y$  in $X$ (possibly, with self intersections),  such that
    $$\int_Y scal(X)dy\leq N\cdot A  \mbox  { for all }  N=1,2,3,... .$$

Notice that this definition makes sense only if $scal(X)\geq 0$; otherwise,  $sc.ar(X)=-\infty.$ \vspace {1mm}

 Let  the universal covering $\tilde X$ of the  manifold $X$ be  spin, e.g. homeomorphic to a product of  $2$-spheres, or more generally, of   complex projective spaces $\mathbb CP^l$ for {\it odd $ l$}, and of $\mathbb R^k$. 
 
 Also assume that
 the class $h^k$ has  {\it infinite $K$-area} in the sense of \cite {positive}, e.g. this class   comes from (a cohomology class  of)  a complete  manifold $Z'$  of non-positive curvature by a map $X\to Z'$.

 \vspace {1mm}

{\it Then $sc.ar(X)\leq const_{top}$ where $const_{top}<\infty$ depends  on  the topology of $X$. }

 \vspace {1mm}
 
 {\it Proof}. By the Whitney-Hahn-Banach  (duality) theorem, the fundamental class of $X$ is representable by  a product of smooth closed  $2$-forms $\omega_i$
 with their sup-norms bounded by 
 $$||\omega_i||_{sup}\leq \kappa+\varepsilon \mbox  { for } \kappa= (sc.ar(X))^{-1}  \mbox { and all } \varepsilon>0.$$

 Since all $\omega_i$ are representable as curvature forms of  complex line bundles 
 over $X$,  theorem 
 5$\frac{1}{4}$ from \cite {positive} applies and the proof follows.
 \vspace {1mm}

{\it Remark/Question}. the Riemannian products of $2$-spheres by the torus, are probably extremal  for this inequality, i.e.
 $sc.ar(X)= const_{top}$ for these $X$.
 
 This seems to follow from (a suitable form of)  {\it the area extremality} of the product of $2$-spheres with arbitrary metrics of positive curvatures with the flat torus  $\mathbb T^k$,  see 
  \cite {min-oo2},  \cite {listing},  \cite {listing2},  \cite {goette},  but I did not hat it carefully.

More general products of complex projective spaces with K\"ahler metrics, 
  $\left (\times_i \mathbb CP^{l_i}\right)\times \mathbb T^k$
 are  also extremal \cite {goette2} and this probably, extends to singular metrics  with curvatures concentrated along divisors.   \vspace {1mm}

 {\it Potential  Corollary.}  Let $P$ be an $n$-dimensional mean curvature convex  cornered spin  manifold that is combinatorially equivalent to   a Cartesian product $P_0$ of reflection domains in the spheres $S^2_i$ and in $\mathbb R^k$.

Observe that such a  $P$ serves as a reflection domains in a closed manifold $\hat X$  that  admits a 
 map with degree one onto a product of $2$-spheres  $S^2_i$ and the torus $\mathbb T^k$.

  Let  $scal(P)=0$, e.g. $P$ is isometric to a cornered domain in $\mathbb R^n$.  
 
 Then the above implies a  certain lower 
 bound on the dihedral angles  $\alpha_j$ of $P$ as follows.

 Assume, to simplify the  notation, that the above reflection (Cartesian product)   domain $P_0$  has its all dihedral angles $\pi/2$. Thus all all spherical reflection domains $P_i\subset S^2_i$
 in this case are spherical triangles with $90^\circ$ angles.

Let us  endow  the above closed manifold $\hat X\supset P$  (obtained by reflecting $P$)  with the singular Riemannian metric  $\hat g$  that extends that on $P$ and observe that 
 the essential contribution to (non-negative!) scalar curvature of this metric comes from the edges, i.e. codimension $2$ faces of $P$. 
 
 Namely, the contribution of such an edge with the dihedral angle $\alpha_j$  to  a surface $Y$ transversally intersecting the "descendent" of this edge in $X$ equals $2\pi-4\alpha_j$. 
 
 Thus,  assuming the truth of the above extremality statement for  of $\times S^2_i\times \mathbb T^k$,  \vspace {1mm}

 {\it the sums   $\Sigma_3$  of all triples of dihedral angles in $P$ corresponding to the  triples 
 of  vertices of the reflection triangles  in the spheres  $S^2_i$ are bounded from below  by 
 $\Sigma_3 \geq \pi$.} \vspace {1mm}

Notice that we do not have to take $N$-multiples of our $2$-dimensional  homology classes in this case.  

Also observe that if  $n=3$, then the above  reduces to a special case of the  Gauss-Bonnet prism inequality, namely for   "triangular"  prisms $P^3=\triangle \times [0,1]$.
 
\vspace {1mm}

{\it Remarks on   Pure Edge Singularities.} The $K$-area inequalities and the related extremality/rigidity results for closed  for Riemannian manifolds 
$(X,g)$ can be expressed in terms of the size/shape of this $X$ with the metric
$scal(g)\cdot g$, where, observe this metric is invariant under scaling of $g$. (Compare  \cite {min-oo2},  \cite {listing},  \cite {listing2},  \cite {goette}.)

Namely, the suitable for the present purpose area extremality of an $X$ says that if another manifold, say $(X',g')$  admits a map $f: X\to X'$ of positive degree than this $f$ can not be strictly  area decreasing with respect to  the metrics $scal(g)\cdot g$ and  $scal(g')\cdot g'$.

The case that is relevant in the present context is that of singular metrics $g$ with the singularity, say $\Sigma \subset X$ being of codimension $2$ and with the main contribution to the scalar curvature being supported on $\Sigma$ transversally to  $\Sigma$.  

One may think of this $\Sigma \subset X$ as a "divisor" in $X$ (divisors in complex manifolds $X$ provide a pool of interesting examples) and  the $K$-area formulas apply to the Dirac operator twisted with the line bundle associated to this "divisor" or to
Whitney sums of such bundles.

This is what we have actually done for the above "prisms" and it would be interesting to look more systematically on more general  such $(X,\Sigma)$.

\vspace {1mm}

{\it Question.}  Can the $K$-area be used for   bounding  from below  some  "combinatorial sizes/shapes"  of mean curvature convex polyhedra $P$ in $\mathbb R^n$, say for $n=3$, with all their dihedral angles bounded by  $\pi-\varepsilon$ for a given $\varepsilon>0$?

 \vspace {2mm}

 \textbf{Hyperbolic Prism Inequality.} The  $3D$ Gauss-Bonnet prism inequality  was  "modeled" on the prisms in the product spaces $V^2_\kappa\times \mathbb R$ for surfaces $V^2_\kappa$ of constant curvature $\kappa$.
 
Now let us look at prisms $P_{horo}$ in the hyperbolic space $H^3$ of curvature $-1$ where the top and the bottom are contained in two parallel horospheres. Observe that the side faces of $P_{horo}$ are totally geodesic, and the mean curvatures of the top and the bottom equal $+1$ and $-1$ respectively.

Now, let $P$ be a normal prism with $scal(P)\geq -6=scal(H^3)$ where the side faces are mean curvature convex while the mean curvatures of the top
and of the bottom are bounded from below by $+1$ and $-1$. 

The above argument, applied to constant mean curvature bubbles $Y_{bbl} \subset P$  (compare \S 5$\frac{5}{6}$ in \cite {positive}), implies  that 
{\it $$\sum_{i=1,...,k} (\pi-\alpha_i)\leq 2\pi,$$
where the equality holds if and only if $P$ is isometric to some  $P_{horo}$ in the hyperbolic space $H^3$.}

\vspace {1mm}

{\it On Spherical Prisms.}  Let is indicate a similar extremality/rigidity  property of  spherical rather
than horospherical prisms that applies to   $3$-manifolds $V^3_\kappa$  of constant curvature $\kappa$ for all $-\infty <\kappa< +\infty$.

Let $ P_{sph}\subset V^3_\kappa$ be a normal prism where the top and the bottom lie in two concentric spheres  in $V^3_\kappa$ (or two concentric umbilical surfaces if $\kappa<0$) and where the side faces of $ P_{sph}$ are totally geodesic.
Denote by $Y_{sph}(d)\subset  P_{sph}$ the $d$-level of the  distance function to the (concave) bottom of $P_{sph}$ and let $m(d)$ denotes the mean curvature of  $Y_{sph}(d)$.

Observe that the scalar curvature of $Y_{sph}(d)$ equals $scal (V^3_\kappa) + \frac{3}{4}m(d)^2$, compare 
\S 5$\frac{5}{6}$ in \cite {positive}, and, thus, the number 
$$S=_{def}area(Y_{sph})(d)\cdot (scal (V^3_\kappa) + \frac{3}{4}m(d)^2)$$
is independent of $d$ and where, clearly, $S\geq 0$  for all $\kappa$.

Let $a_i$, $ i=1,...,k$,
denote  the dihedral angles between the side faces of  $P_{sph}$.

\vspace {1mm}
 Let $P$ be a normal Riemannian prism combinatorially equivalent to $P_{sph}$ (with not necessarily  $scal(P)\geq scal(V^3_\kappa)$) where
the side faces of $P$ are mean curvature convex while the mean curvatures of the top and the bottom are bounded from below by those of $ P_{sph}$. Thus top of $P$ is mean convex while the bottom may be concave. 

Let $\delta =\delta(p)$ denotes the distance from $p$ to the bottom of $P$.

 \vspace {1mm}

{\it Let  the distance between the top and bottom in $P$ is (non-strictly)   greater
than that in $P_{sph}$ and let
$$\int_Y (scal (P) + \frac{3}{4}m(\delta)^2)dy\geq S$$
for all surfaces $Y\in P$ that separate the top from the bottom.} 

 {\it Then the upper bounds $\alpha_i$
on the dihedral angles between the side faces of $P$ satisfy
$$\sum_i\alpha_i\geq \sum_i a_i.$$}\vspace {1mm}

The above spherical  prism inequality can be applied to particular  small domains $P$ in $3$-manifolds $X$ delivered by an argument similar for proving $(\smallsquare)$ from 5.9, where it shows that 

 {\it If a  $C^2$-metric  $g$ on a $3$-manifold equals a $C^0$-limit of  $C^2$-metrics $g_i$ with $scal(g_i)\geq \kappa$ for a given $-\infty<\kappa<+\infty$, then  $scal(g)\geq \kappa$}.

In fact, one   needs for this purpose only  polyhedral domains of depth $2$, i.e. {\it without corners}, namely, those combinatorially equivalent to round $3$-cylinders, where the relevant bubbles, that separate top from the bottom, are $C^2$-smooth at the boundary. \vspace {1mm}

{\it Question.} Can one prove all of the   above   "prism  inequalities"  by a purely Dirac-type argument (e.g.  elaborating on that  in \cite{min-oo1} and/or in \cite{davaux}  and where \cite{lusztig} may be relevant)  with no localization to  minimal hypersurfaces? \vspace {1mm}

\vspace {1mm}

{\it  More on Pure Edges.} The study of the above $P$ can be reduced  by (multi)-doubling  to the "pure edge" picture (see 1.1., 2.3). 

Then one observes that the semi-integral  inequalities are stable under smoothing of the edges since these smoothings   that do not significantly  change the integrals of the scalar curvatures over the relevant  surfaces; thus, one can fully  avoid the  "singularities at the corners" problem. \vspace{1mm}

It seems  that the geometry of  manifolds  $X$ of all dimensions   with $scal\geq 0$ (and, possibly, with $scal\geq \kappa$) is governed by integrals  of $scal(X)$ over surfaces  in $X$, or, rather   by 
infima of such integrals over surfaces $A$ with (small?)   boundaries $B$ and  with $A$ being "not-too-far" from $B=\partial A$.

However,  this is insufficient, for instance,  for extending the  (dual)  Kirszbraun theorem  to  mean curvature convex     simplices $P=\Delta_{curv}^3$  with $scal (P)\geq 0$   that is  needed for  proving the dihedral extremality  of the ordinary simplices  $\Delta^3 \subset \mathbb R^3$   (see section 1.6).

But, possibly, this may be achieved  with  $\mu$-bubbles  $Y\subset  P$,  $\partial Y\subset  \partial P$, that are  obtained by  inward deformations of   the surface $Y_0 \subset \partial P=\partial  \Delta_{curv}^3$, for $Y_0$ being   $\partial \Delta_{curv}^3$  minus four small discs around   four vertices of 
 $\Delta_{curv}^3$, and where $\mu$ are   suitable   measures  supported on $\partial \Delta^3_{curv}$.
 
 Such $Y$ do not have to be normal to the edges that  may increase the integral of the scalar curvatures of  (regularized meteoric on) $\Delta^3$ over $Y$.

 \vspace {1mm}

{\it  More  on $n$-Dimensional "Prisms", etc.} 
Application of minimal hypersurfaces (see \cite  {schoen-yau2}, \cite  
{gromov-lawson3}   depends on their  regularity. This limits one to $n\leq 7$, where the case $n=8$ follows from the removal of isolated  singularities \cite {smale} with the general case, probably, amenable to Lohkamp's  "around singularities" method.  (Removal and/or "going around" singularities is good for proving extremality of cornered domains but it is poorly adjusted for   rigidity.) 

But even disregarding the singularity problems, (including that at the corners) one needs to adjust the warping  from  \cite  {gromov-lawson3} to the above Gauss-Bonnet argument.

On the other hand, as we saw above, the Dirac operator/$K$-area method applies in some cases.

Furthermore, let  for example,  $Y_{horo}^{n-1}$ be a Euclidean reflection domain in a horosphere in the hyperbolic $n$-space, then the study of the corresponding prisms
reduces to that of manifolds homeomorphic to $\mathbb T^{n-1}\times \mathbb R$, where the corresponding result is proven in \cite {positive}
for $n\leq 7$. This sharpens/quantifies  Min-Oo's rigidity theorem \cite {min-oo1}\cite {min-oo2} for hyperbolic $n$-spaces in the case  $n\leq 7$ (with \cite {smale} one gets extremality but not rigidity for $n=8$) but leaves 
$n\geq 9$ open. (Possibly, Lohkamp's method  applies here.)\vspace {1mm}

{\it What are, in general, extremal polyhedral  objects in symmetric spaces $X$ of non-constant curvatures? }

If $X$ is a Cartesian product of constant curvature spaces, then natural candidates are Cartesian  products of poly-bubble-hedra  but it is unclear 
 waht happens in general. On the other hand, there are several rigidity (and non-rigidity)  results for  compact and complete  manifolds without boundary \cite {min-oo1}  \cite {min-oo2}  \cite {llarull}  \cite {listing}\cite {brendel}.

 \vspace {1mm}

 {\it Model Problem.} Is there a single theory of manifolds $X$  with $scal(X)\geq \kappa$
or there may be several theories associated with different (homogeneous? symmetric?) model spaces that serve as extrema for geometric/topological inequalities  in spaces with $scal\geq \kappa$?
 
 \vspace {1mm}

{\it Rigidity and Stability.}  All(?)   geometric/topological inequalities, in particular, those  concerning (smooth as well as cornered)  manifolds $X$ with $scal(X)\geq \kappa$, whenever these are sharp,  are accompanied by rigidity
problems where one asks for a  description of  $X$ where such inequalities become equalities.

 But even when  such  rigidity of an  $X$ is known, e.g. for $X=\mathbb T^n$ with flat metric, it is not quite clear in what sense such an $X$ is stable. For example, one can (almost) unrestrictedly blow  huge  "bubbles"  with positive scalar curvature 
  "grounded" in codimension $k$ subsets in a Riemannian manifold $X$ for $k\geq 3$. But the above semi-integral inequalities indicate that this is impossible with $k\geq 2$.

 This, probably,  can be  interpreted as follows:  \vspace {1mm}
 
 there is a particular "Sobolev type  weak metric" $dist_{wea}$  in the space of $n$-manifolds $X$, such that,  for example,  tori $X$ with $scal(X)\geq-\varepsilon$, when properly normalized,  (sub)converge to 
  flat tori,    $dist_{wea} (X,\mathbb T^n_{flat})\to 0$ for $\varepsilon\to 0$, but  these $X$ may, in general, diverge in   stronger metrics.    \vspace {1mm}

 (A preliminary step toward  construction of  such  a  metric can be seen in  \cite {sormani}.)

\subsection{Spaces or Objects?}

The Dirac operator $D$ and minimal hypersurfaces seem to belong to different wolds but they unexpectedly meet in
spaces with lower bounds on their scalar curvatures. \vspace {1mm}

{\it Is there a deeper link between "Dirac" and "minimal".} \vspace {1mm}

It should be noted that the  constrains on the {\it geometry} of  manifolds $X$  with $scal \geq 0$ (and $scal(X) \geq \kappa$, in general), and the issuing constrains on the topology of $X$, that obtained with the Dirac operator are significantly different from those obtained with minimal hypersurfces, even for  spin manifolds of dimension $\leq 7$ where both methods apply. 

An essential drawback of  a {\it direct} application of   Dirac operators  is the requirement of {\it completeness} of $X$, while the  use of minimal surfaces in dimensions $n\geq 4$ delivers upper {\it distance bounds} rather  than  {\it area bounds} as the Dirac operators do. (Particularly fine results can be obtained for $dim(X)=4$ with a use of the Seiberg-Witten equation but this is fully beyond the scope of the present paper.)

 \vspace {1mm}

{\it Can one unify the two methods and, thus,  obtains more precise  results?}\vspace {1mm}

The key to the application of Dirac to scalar curvature is the Lichnerowicz formula $D^2=\nabla^2+scal/4$ where $\nabla^2$ is more positive than an ordinary Laplace operator, while the application of minimal (hyper)surfaces depends on the second variation formula in the Schoen-Yau form that makes it similar to Lichnerowicz formula. \vspace {1mm}

{\it Are  there further formulas of this kind mediating between these two?} \vspace {1mm}

The next step in   approaching the problem is to change our perspective on manifolds
in the spirit of ideas of Fedia Bogomolov and Maxim Kontsevich.

Fedia suggested looking at stable vector bundles with zero $c_1$ on an algebraic variety $X$  as  at  "coherent families" of flat bundles over the curves $C\subset X$, while Maxim's idea was to regard Riemannian manifolds as (special cases of)  functors from the "category of graphs" to the category of measure spaces: the value of such a functor $F$ on a graph
$G$ with given edge lengths is the space of maps $G\to X$ with the Wiener measure on it.

If we want to see  the Dirac operator  in Maxim's picture, we have  to consider graphs  $G$
along with flat  $O(n)$-bundles $V$, $n=dim X$, over them and with an embedding of the "tangent spaces"
of the graphs at all points $g\in G$ into $V_g$ (that is most informative at the vertices  $g\in G$ of high valency). Then one needs to define an appropriate structures, including measures,
on the the spaces of maps $(G,V)\to  (X,T(X))$ for the tangent bundle $T(X)$.

The flat structures in $V$  encode the parallel transport in $X$ that allows one to speak of "Dirac" as well as of other geometric differential operators. On the other hand, a measure
on the space of maps $G\to X$ allows an integration of the numbers of intersections of
these graphs with hypersurfaces $Y\in X$, thus keeping track of the volumes of these $X$.

\vspace {1mm}
Besides all this, one can associate with every cycle in $G$ the range of values of the 
areas of surfaces $A$  fillings this cycles in  $X$ along with the  integrals $\int_A scal(X) da$, where, possibly, the values of the genera and Euler characteristics  of $A$ may be also  relevant.

Possibly,  the spirit of  semi-integral inequalities suggests  that  we may actually forfeit graphs  and think of $X $  as a contravariant  functor from  a suitable  category of surfaces  to some set category.\vspace {1mm}

{\sc Final Questions.} Can one make a mathematical theory  along these lines  with manifolds  being replaced  by objects of a more abstract and more flexible category  of  functors from a category of   "extended objects" to sets?

Can one, thus, extend basic results on positive scalar curvature to  {\it singular}  and/or {\it infinite dimensional} spaces?

\end{document}